\documentclass[11pt,a4paper,usenames,dvipsnames]{article}

\usepackage{url,amsmath,amssymb,latexsym,mathrsfs,comment,amsthm,enumerate,geometry,calc,ifthen,mathdots,textcomp,multicol}

\usepackage{cases}

\usepackage[noadjust]{cite}

\usepackage[pagebackref,colorlinks]{hyperref}

\hypersetup{pdfborder={0 0 0}}

\renewcommand*{\backref}[1]{}
\renewcommand*{\backrefalt}[4]{%
      \ifcase #1 %
      \relax
      \or
      (Page #2).%
      \else
      (Pages #2).%
      \fi%
}

\usepackage[margin=10pt,font=small,labelfont=bf, labelsep=period]{caption}

\usepackage{makecell}

\usepackage[all,cmtip,2cell]{xy}

\usepackage[x11names, svgnames, rgb,table]{xcolor}
\usepackage{hhline}
\usepackage{multirow}

\usepackage{pst-node}
\usepackage{tikz-cd}

\usepackage{enumitem}

\usepackage[utf8]{inputenc}
\usepackage[T1]{fontenc}

\usepackage{tikz}
\usepgflibrary{snakes,arrows,shapes}
\pgfdeclarelayer{background layer}
\pgfsetlayers{background layer,main}

\usetikzlibrary{decorations.markings}
\usetikzlibrary{arrows,matrix}
\usepgflibrary{arrows}
  
\tikzset{->-/.style={decoration={
  markings,
  mark=at position #1 with {\arrow{{latex}}}},postaction={decorate}}}
\newcommand\nc\newcommand
\nc\rnc\renewcommand

\nc\trans[1]{\left(\begin{smallmatrix}#1\end{smallmatrix}\right)}

\nc\FIG{\operatorname{\mathsf{IG}}}
\nc\FPG{\operatorname{\mathsf{PG}}}
\nc\bbB{\mathbb B}
\nc\N{\mathbb N}
\nc\Z{\mathbb Z}
\nc\TL{\mathcal T\!\mathcal L}
\rnc\S{\mathcal S}

\newcommand{\uv}[1]{\fill (#1,2)circle(.17);}
\newcommand{\lv}[1]{\fill (#1,0)circle(.17);}
\newcommand{\tuv}[1]{\fill (#1,2)circle(.2);}
\newcommand{\tlv}[1]{\fill (#1,0)circle(.2);}
\newcommand{\uvs}[1]{{\foreach \x in {#1} { \uv{\x}}}}
\newcommand{\lvs}[1]{{\foreach \x in {#1} { \lv{\x}}}}
\newcommand{\darcx}[3]{\draw(#1,0)arc(180:90:#3) (#1+#3,#3)--(#2-#3,#3) (#2-#3,#3) arc(90:0:#3);}
\newcommand{\darc}[2]{\darcx{#1}{#2}{.4}}
\newcommand{\uarcx}[3]{\draw(#1,2)arc(180:270:#3) (#1+#3,2-#3)--(#2-#3,2-#3) (#2-#3,2-#3) arc(270:360:#3);}
\newcommand{\uarc}[2]{\uarcx{#1}{#2}{.4}}
\newcommand{\stline}[2]{\draw(#1,2)--(#2,0);}
\nc{\uarcs}[1]{{\foreach \x/\y in {#1}{ \uarc{\x}{\y} }}}
\nc{\darcs}[1]{{\foreach \x/\y in {#1}{ \darc{\x}{\y} }}}
\newcommand{\stlines}[1]{{\foreach \x/\y in {#1} { \stline{\x}{\y} }}}

\newcommand{\darcxcol}[4]{\draw[#4](#1,0)arc(180:90:#3) (#1+#3,#3)--(#2-#3,#3) (#2-#3,#3) arc(90:0:#3);}
\newcommand{\darccol}[3]{\darcxcol{#1}{#2}{.4}{#3}}
\nc{\darccols}[2]{{\foreach \x/\y in {#1}{ \darccol{\x}{\y}{#2} }}}
\newcommand{\uarcxcol}[4]{\draw[#4](#1,2)arc(180:270:#3) (#1+#3,2-#3)--(#2-#3,2-#3) (#2-#3,2-#3) arc(270:360:#3);}
\newcommand{\uarccol}[3]{\uarcxcol{#1}{#2}{.4}{#3}}
\nc{\uarccols}[2]{{\foreach \x/\y in {#1}{ \uarccol{\x}{\y}{#2} }}}
\newcommand{\stlinecol}[3]{\draw[#3](#1,2)--(#2,0);}
\newcommand{\stlinecols}[2]{{\foreach \x/\y in {#1} { \stlinecol{\x}{\y}{#2} }}}

\nc{\uvert}[1]{\fill (#1,2)circle(.2);}
\rnc{\lvert}[1]{\fill (#1,0)circle(.2);}
\nc{\custpartn}[3]{{\lower1.4 ex\hbox{
\begin{tikzpicture}[scale=.3]
\foreach \x in {#1}
{ \uvert{\x}  }
\foreach \x in {#2}
{ \lvert{\x}  }
#3 \end{tikzpicture}
}}}

\newcounter{ncols}
\newcounter{incols}
\newenvironment{partn}[1]{
  \setcounter{ncols}{#1} \setcounter{incols}{\thencols - 1}\setlength{\arraycolsep}{1pt}
  \Bigl( \hspace{-1.5truemm}\scriptsize 
    \begin{array}{@{\hskip 3pt} c *{\theincols}{|c} @{\hskip 3pt}  }
}{
     \end{array}
     \normalsize \hspace{-1.5truemm}\Bigr)\setlength{\arraycolsep}{6pt}
}

\nc\la\langle
\nc\ra\rangle

\nc\ol\overline
\nc\ul\underline
\nc\da\downarrow
\nc\pr\bullet
\nc\ccirc{\circ'}
\nc\bc\odot
\nc\sm\setminus

\nc\Grp{{\bf Grp}}
\nc\SL{{\bf SL}}
\nc\RSS{{\bf RSS}}
\nc\OG{{\bf OG}}
\nc\IGRSS{{\bf IGRSS}}
\nc\IGRSSP{\IGRSS_P}
\nc\PA{{\bf PA}}
\nc\PG{{\bf PG}}
\nc\WCPG{{\bf WCPG}}
\nc\TCPG{{\bf TCPG}}
\nc\CPG{{\bf CPG}}
\nc\Set{{\bf Set}}
\nc\IS{{\bf IS}}
\nc\IG{{\bf IG}}
\nc\bG{{\bf G}}
\nc\bS{{\bf S}}
\nc\bC{{\bf C}}
\nc\bD{{\bf D}}
\nc\bI{{\bf I}}
\nc\bT{{\bf T}}
\nc\bbG{\mathbb G}
\nc\bbS{\mathbb S}

\nc\JErev[1]{\textcolor{magenta}{#1}}

\nc\ds\displaystyle
\nc\ts\textstyle

\nc\bn{{\bf n}}

\nc\al\alpha
\nc\be\beta
\nc\ga\gamma
\nc\de\delta
\nc\ve\varepsilon
\nc\lam\lambda
\nc\om\omega
\nc\ze\zeta
\nc\si\sigma
\nc\Si\Sigma
\nc\Ga\Gamma
\nc\De\Delta
\nc\Om\Omega
\nc\io\iota

\nc\nab\nabla

\nc\vt\vartheta
\rnc\th\theta
\nc\Th\Theta

\nc\BY{\qquad\text{by}\qquad}
\nc\GIVENBY{\qquad\text{given by}\qquad}
\nc\AND{\qquad\text{and}\qquad}
\nc\ANDSIM{\qquad\text{and similarly}\qquad}
\nc\ANDSO{\qquad\text{and so}\qquad}
\nc\ANd{\quad\text{and}\quad}
\nc\COMMA{,\qquad}
\nc\COMMa{,\quad}
\nc\WHERE{\qquad\text{where}\qquad}

\rnc\iff{\ \Leftrightarrow\ }
\nc\IFf{\quad \Leftrightarrow\quad }
\nc\Iff{\ \ \Leftrightarrow\ \ }
\nc\IFF{\qquad \Leftrightarrow\qquad }
\rnc\implies{\ \Rightarrow\ }
\nc\IMPLIES{\qquad \Rightarrow\qquad }

\nc\set[2]{\{#1:#2\}}
\nc\bigset[2]{\big\{#1:#2\big\}}
\nc\pres[2]{\la#1:#2\ra}

\nc\bit{\begin{itemize}[label=\textbullet, leftmargin=5mm]}
\nc\eit{\end{itemize}}
\nc\ben{\begin{enumerate}[label=\textup{(\roman*)},leftmargin=10mm]}
\nc\bena{\begin{enumerate}[label=\textup{(\alph*)},leftmargin=10mm]}
\nc\een{\end{enumerate}}
\nc\bmc{\begin{multicols}}
\nc\emc{\end{multicols}}

\nc\pf{\begin{proof}}
\nc\epf{\end{proof}}
\nc\pfclaim{\begin{quote}\begin{proof}}
\nc\epfclaim{\end{proof}\end{quote}}
\nc\epfres{\hfill\qed}
\nc\epfreseq{\tag*{\qed}}
\let\oldproofname=\proofname
\renewcommand{\proofname}{\rm\bf{\oldproofname}}
\nc{\pfitem}[1]{\medskip \noindent #1.}
\nc{\firstpfitem}[1]{#1.}
\nc{\pfcase}[1]{\medskip\noindent {\bf Case #1.}}
\nc\aftercases{\medskip\noindent}

\renewcommand{\H}{\mathrel{\mathscr H}}
\renewcommand{\L}{\mathrel{\mathscr L}}
\newcommand{\R}{\mathrel{\mathscr R}}
\newcommand{\D}{\mathrel{\mathscr D}}
\newcommand{\J}{\mathrel{\mathscr J}}
\newcommand{\K}{\mathbb K}
\nc\rH{\mathrel{\H}}
\nc\rL{\mathrel{\L}}
\nc\rR{\mathrel{\R}}
\nc\rD{\mathrel{\D}}
\nc\rJ{\mathrel{\J}}
\nc\rK{\mathrel{\K}}
\nc\rsi{\mathrel{\si}}

\nc\leqL{\leq_\L}
\nc\leqR{\leq_\R}
\nc\leqJ{\leq_\J}
\nc\leqH{\leq_\H}
\nc\leqF{\leq_{\F}}
\nc\geqF{\geq_{\F}}

\newcommand{\Sing}{\operatorname{Sing}}

\newcommand{\opp}{{\operatorname{op}}}
\newcommand{\id}{\operatorname{id}}
\newcommand{\dom}{\operatorname{dom}}
\newcommand{\codom}{\operatorname{codom}}
\newcommand{\coker}{\operatorname{coker}}
\newcommand{\rank}{\operatorname{rank}}

\nc\mr\mathrel
\nc\pc[2]{(#1,#2)^\sharp}

\nc\U{\mathcal U}
\nc\V{\mathcal V}
\nc\G{\mathcal G}
\rnc\iff{\ \Leftrightarrow\ }
\rnc\implies{\ \Rightarrow\ }
\nc\Implies{\quad \Rightarrow\quad }
\nc\F{\mathrel{\mathscr F}}
\nc\C{\mathscr C}
\nc\M{\mathcal M}
\nc\CC{\mathcal C}
\nc\DD{\mathcal D}
\nc\FF{\mathcal F}
\nc\I{\mathcal I}
\rnc\O{\mathcal O}
\rnc\P{\mathscr P}
\nc\PP{\mathcal P}
\nc\T{\mathcal T}
\nc\p{\mathfrak p}
\nc\q{\mathfrak q}
\rnc\r{\mathfrak r}
\nc\s{\mathfrak s}
\rnc\t{\mathfrak t}
\nc\bd{{\bf d}}
\nc\br{{\bf r}}
\nc\op\oplus
\nc\lra{\mathrel\leftrightarrow}
\nc\es\varnothing
\nc\rev{\textup{rev}}
\nc\corestt{{\upharpoonleft}}
\nc\restt{{\upharpoonright}}
\nc\corest{{\downharpoonleft}}
\nc\rest{{\downharpoonright}}
\nc\WHERe{\quad\text{where}\quad}
\nc\ot\leftarrow
\rnc\a{\mathfrak a}
\rnc\b{\mathfrak b}
\rnc\c{\mathfrak c}
\rnc\d{\mathfrak d}
\nc\sub\subseteq
\nc\mt\mapsto
\nc\im{\operatorname{im}}
\nc\B{\mathcal B}
\nc\E{\mathbb E}

\numberwithin{equation}{section}

\newtheorem{thm}[equation]{Theorem}
\newtheorem{lemma}[equation]{Lemma}

\newtheorem{prop}[equation]{Proposition}
\newtheorem*{ESN}{Ehresmann--Schein--Nambooripad Theorem}
\newtheorem*{MT}{Main Theorem}

\theoremstyle{definition}

\newtheorem{defn}[equation]{Definition}
\newtheorem{rem}[equation]{Remark}
\newtheorem{eg}[equation]{Example}

\newcounter{caseco}

\newcounter{subcaseco}

\newcounter{stepco}

\newcounter{stageco}

\begin{document}

\title{A groupoid approach to regular $*$-semigroups}
\date{}

\author{James East\footnote{Supported by ARC Future Fellowship FT190100632.} \ and P.A.~Azeef Muhammed\\[3mm]
{\it\small Centre for Research in Mathematics and Data Science,}\\
{\it\small Western Sydney University, Locked Bag 1797, Penrith NSW 2751, Australia.}\\[3mm]
{\tt\small J.East@WesternSydney.edu.au}, {\tt\small A.ParayilAjmal@WesternSydney.edu.au}}

\maketitle

\begin{abstract}
%\noindent
In this paper we develop a new groupoid-based structure theory for the class of regular $*$-semigroups.  This class occupies something of a `sweet spot' between the important classes of inverse and regular semigroups, and contains many natural examples.  Some of the most significant families include the partition, Brauer and Temperley-Lieb monoids, among other diagram monoids.

Our main result is that the category of regular $*$-semigroups is isomorphic to the category of so-called `chained projection groupoids'.  Such a groupoid is in fact a triple $(P,\mathcal G,\varepsilon)$, where:
\bit
\item $P$ is a projection algebra (in the sense of Imaoka and Jones),
\item $\mathcal G$ is an ordered groupoid with object set $P$, and 
\item $\varepsilon:\mathscr C\to\mathcal G$ is a special functor, where $\mathscr C$ is a certain natural `chain groupoid' constructed from $P$.  
\eit
Roughly speaking: the groupoid $\mathcal G=\mathcal G(S)$ remembers only the `easy' products in a regular $*$-semigroup $S$; the projection algebra $P=P(S)$ remembers only the `conjugation action' of the projections of $S$; and the functor $\varepsilon=\varepsilon(S)$ tells us how $\mathcal G$ and $P$ `fit together' in order to recover the entire structure of $S$.  In this way, we obtain the first completely general structure theorem for regular $*$-semigroups.

%Among other applications, we use our structure theorem to obtain new, and perhaps more transparent, constructions of fundamental regular $*$-semigroups.  We also use the chain groupoids to demonstrate the existence of free (idempotent-generated) regular $*$-semigroups associated to arbitrary projection algebras.  

As a consequence of our main result, we give a new proof of the celebrated Ehresmann--Schein--Nambooripad Theorem, which establishes an isomorphism between the categories of inverse semigroups and inductive groupoids.  Other applications will be given in future works.

We consider several examples along the way, and pose a number of problems that we believe are worthy of further attention.

\medskip

\medskip


\end{abstract}

\newpage

\tableofcontents

\newpage

\section{Introduction}\label{sect:intro}

The \emph{ESN Theorem} is arguably the most important result in inverse semigroup theory:

\begin{ESN}
%The categories of inverse semigroups (with semigroup homomorphisms) and inductive groupoids (with inductive functors) are isomorphic.
The category $\IS$ of inverse semigroups (with semigroup homomorphisms) is isomorphic to the category $\IG$ of inductive groupoids (with inductive functors).
\end{ESN}

This result was first formulated as above by Lawson in \cite[Theorem 4.1.8]{Lawson1998}, who named it after the three mathematicians who had contributed most to its development.  The key motivation for the theorem is that it ties together the two main approaches to partial symmetries 
%in non-Euclidean geometries 
that were not captured by groups: the \emph{inductive groupoids} of Charles Ehresmann \cite{Ehresmann1957,Ehresmann1960}, and the \emph{inverse semigroups} of Viktor Wagner~\cite{Wagner1952,Wagner1953} and Gordon Preston~\cite{Preston1954,Preston1954b}.  A fuller account of the historical background can be found in \cite{Lawson1998,HL2017,Hollings2012}; see also \cite{MPR2021}.

The exact nature of the isomorphism in the ESN Theorem reveals a lot more than is evident in the above statement.  Not only are the \emph{categories} $\IS$ and $\IG$ `the same', but in fact the \emph{objects} of these categories are, too:  inverse semigroups \emph{are} inductive groupoids, and conversely.  More formally, an inverse semigroup~$S$ can be identified with the groupoid~$\G(S)$ with:
\bit
\item object/identity set $E(S)$, the semilattice of idempotents of $S$, and
\item morphism set $S$, where a semigroup element $a$ is thought of as a morphism $aa^{-1}\to a^{-1}a$.
\eit
Together with a natural order, this is enough to recover the entire structure of $S$.  At the morphism level, inverse semigroup homomorphisms $S\to S'$ and inductive functors $\G(S)\to\G(S')$ are precisely the same maps.
%The key step here is to show that the inductive groupoid axioms are enough to extend the (partial) composition to the (total) product in an associated inverse semigroup.  
%We will say more about this in Section \ref{sect:I}.
We will give more detailed information (including a new proof of the theorem) in Section \ref{sect:I}, but for now it is worth noting that the (partial) composition in $\G(S)$ is a restriction of the (total) product in~$S$.  Roughly speaking, the groupoid `remembers' only the `easy' products in the semigroup.  This is particularly useful in defining inverse semigroups (from inductive groupoids) or inverse semigroup homomorphisms (from inductive functors).

The ESN Theorem has had a major impact outside of (inverse) semigroup theory, with one particularly fruitful direction being the application of inverse semigroups to C$^*$-algebras.  For example, the articles \cite{Paterson2002,Li2012,Li2017,LRRW2018, BLS2017,MS2014,KS2002, FMY2005,Lawson2012,BEM2017} display a mixture of semigroup-theoretic and groupoid-based approaches.  For more extensive discussions, and many more references, see Lawson's survey \cite{Lawson2020} and Paterson's monograph~\cite{Paterson1999}.  See also 
%\cite{Steinberg2008,Steinberg2010,Steinberg2006, Stein2016,Stein2017,Stein2019,Stein2020,MS2021,Stein2022} 
\cite{Steinberg2008,Steinberg2010,Steinberg2006} 
for applications to representation theory.

Another way to measure the influence of the ESN Theorem is to consider its extensions and generalisations.  When faced with an important class of semigroups, one naturally looks for an `ESN-type link' with an equally-natural class of categories, and vice versa.  Thus, we have ESN-type theorems for regular semigroups, restriction semigroups, Ehresmann semigroups, concordant semigroups, DRC semigroups, and others; see for example \cite{Wang2016,Lawson2004,DP2018,FitzGerald2019,GW2012,FitzGerald2010,Armstrong1988,Wang2020,Wang2019 ,Nambooripad1979 ,Lawson1991,Lawson2021 ,Hollings2012,GH2010,Gould2012,Hollings2010,
Wang2022,FK2021,Stokes2017,Stokes2022,Stokes2022b
}.  As with several `dualities' in mathematics \cite{CD1998,Stone1936,Birkhoff1937,Pontrjagin1934a,Pontrjagin1934b,Priestley1970,Priestley1972}, such correspondences allow problems in one field to be translated into another, where they can hopefully be solved, and the solution then reinterpreted in the original context.\footnote{We note, however, that in these dualities morphisms are reversed when moving from one category to the other, whereas the directions are consistent in the ESN Theorem and its variations just cited.}

Arguably the most significant of the papers just cited was Nambooripad's~1979 memoir~\cite{Nambooripad1979}.  As well as pioneering the categorical approach (which led to the `N' in `ESN'), this paper majorly extended the ESN Theorem to the class of \emph{regular} semigroups.  These are semigroups in which every element has \emph{at least one} inverse, rather than a unique one as for inverse semigroups.  Dropping the uniqueness restriction on inverses leads to a far more general class of semigroups, which contains many additional natural examples, including semigroups of mappings, linear transformations, and more.  The increase in generality of the semigroups led inevitably to an increase in the complexity of the categorical structures modelling them.  Although Nambooripad still represents a regular semigroup by an ordered groupoid, the groupoid is not enough to completely recover the semigroup; an additional layer of structure is required.  The need for this extra data is due to the fact that the object set of such a groupoid is not a semilattice (as for inverse semigroups), but instead a \emph{(regular) biordered set}: a partial algebra with a pair of intertwined pre-orders satisfying a fairly complex set of axioms.  The increase in generality also led to the sacrifice of the category \emph{isomorphism}.  Nambooripad's main result, \cite[Theorem~4.14]{Nambooripad1979}, is that the category of regular semigroups is \emph{equivalent} to a certain category of groupoids (which Nambooripad also called \emph{inductive}).  

One way to understand this sacrifice is via the \emph{loss of symmetry} when moving from inverse to regular semigroups.  As noted above, inverse semigroups were devised to model partial symmetries that were unrecognisable by groups.  This is formalised in the Wagner--Preston Theorem, which states that any inverse semigroup is (isomorphic to) a semigroup of partial symmetries of some mathematical structure; see \cite[Theorem 1.5.1]{Lawson1998}.  But an inverse semigroup is itself an extremely symmetrical mathematical structure in its own right, so much so that the canonical proof of the Wagner--Preston Theorem involves showing that an inverse semigroup is (isomorphic to) a semigroup of partial symmetries of \emph{itself}.

The above is not true of arbitrary regular semigroups, however, as any regular semigroup representable by partial symmetries is necessarily inverse.  In a certain sense, an over-abundance of inverses leads to an over-abundance of idempotents, and then to ambiguity in identifying the domain and range of a semigroup element/morphism, hence the necessity of Nambooripad's more elaborate groupoids.

The current article concerns a class of semigroups occupying something of a `sweet spot' between inverse and regular semigroups, the so-called \emph{regular $*$-semigroups} of Nordahl and Scheiblich \cite{NS1978}.  These are semigroups with an additional unary operation $a\mt a^*$ satisfying the laws:
\[
(a^*)^* = a = aa^*a \AND (ab)^* = b^*a^*.
\footnote{Other names exist in the literature for such semigroups.  For example, they are sometimes called \emph{$*$-regular semigroups} \cite{Adair1982,Polak2001}, even though the latter typically refers to a different (larger) class of semigroups \cite{Drazin1979}.  Nambooripad and Pastijn used the term \emph{special $*$-semigroups} \cite{NP1980}.}
\]
The class of regular $*$-semigroups has received a great deal of attention in recent years, partly because of the prototypical examples of the so-called \emph{diagram semigroups}.  These semigroups include families such as the Brauer, Temperley-Lieb, Kauffman and partition monoids, and have their origins and applications in a wide range of mathematical and scientific disciplines, from representation theory, low-dimensional topology, statistical mechanics, and many more; see for example \cite{ER2022,LF2006,FL2011,KV2019,ACHLV2015,Auinger2014,Auinger2012,KM2006,
BDP2002,Jones1994_2,Jones1987,Jones1994_a,Martin1996,Martin1994,Martin2015,EMRT2018,ER2020,ER2022b,ER2022c, MM2007,Maz2002,Wenzl1988,LZ2015,HR2005,Wilcox2007,BH2014,BH2019,LZ2012,DEEFHHLM2019,DEEFHHL2015, MR1998,MM2014,Brauer1937,TL1971,Kauffman1987,Kauffman1990,Jones1983_2,EG2017}.  These monoids have provided a strong bridge between semigroup theory and these other disciplines; for a fuller discussion of this fruitful dialogue, and for many more references, see the introduction to~\cite{EG2017}.

While an element $a$ of a regular $*$-semigroup $S$ might not have a \emph{unique} inverse, it certainly has a \emph{distinguished} inverse, namely $a^*$.  This resolves the above `ambiguity', and allows one to define a groupoid $\G=\G(S)$ in entirely the same way as for inverse semigroups, where now we think of $a$ as a morphism $aa^*\to a^*a$.  While these domain/range idempotents do not generally account for \emph{all} the idempotents of $S$, they are precisely the so-called \emph{projections}, i.e.~the elements $p\in S$ satisfying $p^2=p=p^*$.  The set $P=P(S)$ of all such projections forms a (unary) \emph{projection algebra}, as originally formulated by Imaoka \cite{Imaoka1983}, albeit under a different name.  As we will see, there are important algebraic and order-theoretic ties between the groupoid $\G$ and the algebra~$P$, resulting in what we will call an \emph{(abstract) projection groupoid} $(P,\G)$.  An additional link comes in the form of a certain functor $\ve=\ve(S):\C\to\G$, where $\C=\C(P)$ is the so-called \emph{chain groupoid} of~$P$; this groupoid, which will be defined later, bears some resemblance to Nambooripad's groupoid of \emph{$E$-chains} \cite{Nambooripad1979}.  One consequence of the current work is that the triple $(P,\G,\ve)$, which we will call a \emph{chained projection groupoid}, contains enough information to completely recover the structure of $S$.  
We go much further than this, however, and prove the~following:

\begin{MT}[see Theorem \ref{thm:iso}]
The category $\RSS$ of regular $*$-semigroups (with $*$-semigroup homomorphisms) is isomorphic to the category $\CPG$ of chained projection groupoids (with chained projection functors).
\end{MT}

The definition of some of the terms in the above theorem will be deferred until later sections.  For now we will give a brief outline of the structure of the paper.  
Section \ref{sect:prelim} contains preliminary material and results on semigroups and (small) categories.  
Regular $*$-semigroups are introduced in Section \ref{sect:RSS}, where among other things we construct a functor $\RSS\to\OG$ into the category of ordered groupoids (see Proposition \ref{prop:calGfunctor}).  
Section~\ref{sect:P} concerns the category $\PA$ of projection algebras; we again construct a functor $\RSS\to\PA$ (Proposition \ref{prop:Pfunctor}), and give a new proof of the known fact \cite{Imaoka1983} that every (abstract) projection algebra is the algebra of projections of some regular $*$-semigroup (Proposition \ref{prop:PtoS}).  
Section \ref{sect:CP} builds on Section \ref{sect:P}, and introduces the chain groupoid $\C(P)$ associated to a projection algebra~$P$; this again leads to a functor $\PA\to\OG$ (Remark \ref{rem:Cfunctor}).  
The above ideas are all brought together in Section \ref{sect:G}, where we introduce the category $\CPG$ of chained projection groupoids, and construct a functor $\bG:\RSS\to\CPG$ (Theorem \ref{thm:Gfunctor}).
In Section \ref{sect:GtoS} we construct a functor $\bS:\CPG\to\RSS$ in the reverse direction (Theorem \ref{thm:Sfunctor}).
We then prove our main result (Theorem \ref{thm:iso}) in Section~\ref{sect:iso}, by showing that the functors $\bG$ and $\bS$ are mutual inverses.
We return to the ESN Theorem in Section \ref{sect:I}, where we give a new proof as an application of our main result.
We conclude by discussing future directions and open problems in Section \ref{sect:conclusion}.
Several examples are considered throughout the text, in order to provide concrete cases of our abstract constructions, and to highlight some of the subtleties that arise.  %We also pose a number of problems that we believe are worthy of further investigation.
% (see Sections~\ref{subsect:Seg},~\ref{subsect:D} and~\ref{subsect:Kinyon})

\medskip\noindent\textbf{Acknowledgement.}
We are very grateful to an anonymous referee for a number of insightful questions and comments that led to substantial improvements.

%\newpage

\section{Preliminaries}\label{sect:prelim}

We now gather the preliminary definitions and basic results we require in the rest of the paper.  We cover semigroups in Section \ref{subsect:S}, and (ordered) categories and groupoids in Section \ref{subsect:C}.  
Much of what we present here can be found in standard texts, but we include it to establish notation and conventions, and to keep the paper as self-contained as possible.
For further general background, we refer the reader to texts such as \cite{Howie1995,CPbook} for semigroups,~\cite{Lawson1998} for inverse semigroups, \cite{MacLane1998} for categories, and \cite{BS1981} for universal algebra, though the latter is only used tangentially.

\subsection{Semigroups}\label{subsect:S}

A \emph{semigroup} is a set with an associative binary operation, which will typically be denoted by juxtaposition.  A \emph{monoid} is a semigroup with an identity element.  As usual we write $S^1$ for the \emph{monoid completion} of the semigroup $S$.  So $S^1=S$ if $S$ happens to be a monoid; otherwise, $S^1=S\cup\{1\}$, where $1$ is a symbol not belonging to $S$, acting as an adjoined identity element.

\emph{Green's relations} are five equivalence relations, defined on an arbitrary semigroup $S$ as follows~\cite{Green1951}.  First, for $a,b\in S$ we have
\[
a\R b \iff aS^1=bS^1 \COMMA a\L b \iff S^1a=S^1b \COMMA a\J b \iff S^1aS^1=S^1bS^1. 
\]
Note that $a\R b$ precisely when either $a=b$ or else $a=bx$ and $b=ay$ for some $x,y\in S$; similar comments apply to $\L$ and $\J$.  The remaining two of Green's relations are defined by
\[
{\H} = {\R}\cap {\L} \AND {\D} = {\R}\vee{\L},
\]
where the latter is the join in the lattice of equivalences on $S$, i.e.~the least equivalence containing~${\R}\cup{\L}$.  It is well known that ${\D}={\R}\circ{\L}={\L}\circ{\R}$, and that ${\D}={\J}$ if $S$ is finite.  We denote the $\R$-class of $a\in S$ by
\[
R_a = \set{b\in S}{a\R b},
\]
and similarly for $\L$-classes $L_a$, and so on.

An element $a$ of a semigroup $S$ is \emph{regular} if $a=axa$ for some $x\in S$.  The element $y=xax$ then has the property that $a=aya$ and $y=yay$; such an element $y$ is called an \emph{inverse} of $x$.  We say $S$ is \emph{regular} if every element of $S$ is regular.  The set of \emph{idempotents} of $S$ is denoted
\[
E(S) = \set{e\in S}{e=e^2}.
\]
Idempotents are of course regular.  An \emph{inverse semigroup} is a semigroup in which every element has a unique inverse.  It is well known that $S$ is inverse if and only if it is regular and its idempotents commute; in this case, $E(S)$ is a \emph{semilattice}, i.e.~a semigroup of commuting idempotents.

In this paper we are concerned with a class of semigroups contained strictly between regular and inverse semigroups, the so-called \emph{regular $*$-semigroups} of Nordahl and Scheiblich \cite{NS1978}.  We saw their definition in Section \ref{sect:intro}, and we will return to them in Section \ref{sect:RSS}.

%Such a semigroup $S$ has an additional unary operation $S\to S:a\mt a^*$, with respect to which we have
%\[
%(a^*)^* = a = aa^*a \AND (ab)^* = b^*a^* \qquad\text{for all $a,b\in S$.}
%\]
%We will return to these in Section~\ref{sect:RSS}.

%A \emph{congruence} on a semigroup $S$ is an equivalence relation $\si$ that is compatible with the product of $S$, meaning that 
%\begin{align*}
%a\mr\si b &\IMPLIES ax\mr\si bx \text{ and } xa\mr\si xb &&\text{for all $a,b,x\in S$,}
%\intertext{or equivalently that}
%a\mr\si b \text{ and } x\mr\si y &\IMPLIES ax\mr\si by &&\text{for all $a,b,x,y\in S$.}
%\end{align*}
%Given a congruence $\si$, the \emph{quotient semigroup} $S/\si$ consists of all $\si$-classes under the induced product, $[a][b]=[ab]$.  The \emph{kernel} of a semigroup homomorphism $\phi:S\to T$ is the congruence
%\[
%\ker(\phi) = \set{(a,b)\in S\times S}{a\phi=b\phi}.
%\]
%The fundamental homomorphism theorem for semigroups says that the quotient $S/\ker(\phi)$ is isomorphic to $\im(\phi)$, the image of $S$ under $\phi$.

\subsection{Categories}\label{subsect:C}

Broadly speaking, we consider two types of categories in the current paper:
\bit
\item large categories whose objects are algebraic structures, and whose morphisms are structure-preserving mappings (e.g., the category of all semigroups with semigroup homomorphisms),
\item small categories that are thought of as algebraic structures in their own right.
\eit
Categories of the first kind are treated in a completely standard way (see for example \cite{MacLane1998}).  In the current section we explain how we will work with categories of the second kind. %, and gather some preliminary results.  %All categories discussed in the current section are small.

We typically identify a small category $\CC$ with its set of morphisms.  We identify the objects of~$\CC$ with the identities, the set of which is denoted $v\CC$.  We denote the domain and codomain (a.k.a.~range) of $a\in\CC$ by $\bd(a)$ and $\br(a)$, respectively.  We compose morphisms left to right, so $a\circ b$ is defined if and only if $\br(a)=\bd(b)$, in which case $\bd(a\circ b)=\bd(a)$ and $\br(a\circ b)=\br(b)$.  (Often functors will be written to the left of their arguments, and so composed right to left; it should always be clear which convention is being used.)  For $p,q\in v\CC$, we write
\[
\CC(p,q) = \set{a\in\CC}{\bd(a)=p,\ \br(a)=q}
\]
for the set of all morphisms $p\to q$.

All the small categories we study will have an involution and an order, as made precise in the next two standard definitions.

\begin{defn}\label{defn:*cat}
A \emph{$*$-category}\footnote{These are often called \emph{dagger-} or $\dagger$-\emph{categories} in the literature, and the involution denoted $a\mt a^\dagger$.  See for example \cite{HK2019} or \cite[Section 9.7]{HoTT2013}.} is a small category $\CC$ with an involution, i.e.~a map ${\CC\to\CC:a\mt a^*}$ satisfying the following, for all $a,b\in\CC$:
\begin{enumerate}[label=\textup{\textsf{(I\arabic*)}},leftmargin=10mm]
\item \label{I1} $\bd(a^*)=\br(a)$ and $\br(a^*)=\bd(a)$.
\item \label{I2} $(a^*)^*=a$.
\item \label{I3} If $\br(a)=\bd(b)$, then $(a\circ b)^*=b^*\circ a^*$.
\end{enumerate}
A \emph{groupoid} is a $*$-category for which we additionally have:
\begin{enumerate}[label=\textup{\textsf{(I\arabic*)}},leftmargin=10mm]\addtocounter{enumi}{3}
\item \label{I4} $a\circ a^*=\bd(a)$ (and hence also $a^*\circ a=\br(a)$) for all $a\in\CC$.
\end{enumerate}
In a groupoid, we typically write $a^*=a^{-1}$ for $a\in\CC$.  
\end{defn}

It is easy to show that $p^*=p$ for all $p\in v\CC$, when $\CC$ is a $*$-category.  It is also clear that \ref{I1} is a consequence of \ref{I4}, so a groupoid is a category with a map $a\mt a^*$ satisfying~\ref{I2}--\ref{I4}.

\begin{defn}\label{defn:OC}
An \emph{ordered $*$-category} is a $*$-category $\CC$ equipped with a partial order $\leq$ satisfying the following, for all $a,b,c,d\in\CC$ and $p\in v\CC$:
\begin{enumerate}[label=\textup{\textsf{(O\arabic*)}},leftmargin=10mm] 
\item \label{O1} If $a\leq b$, then $\bd(a)\leq\bd(b)$ and $\br(a)\leq\br(b)$.
\item \label{O2} If $a\leq b$, then $a^*\leq b^*$.
\item \label{O3} If $a\leq b$ and $c\leq d$, and if $\br(a)=\bd(c)$ and $\br(b)=\bd(d)$, then $a\circ c\leq b\circ d$.
\item \label{O4} For all $p\leq\bd(a)$, there exists a unique $u\leq a$ with $\bd(u)=p$.
\end{enumerate}
It is easy to see that \ref{O1}--\ref{O4} imply the following, which is a dual of \ref{O4}:
\begin{enumerate}[label=\textup{\textsf{(O\arabic*)$^*$}},leftmargin=10mm] \addtocounter{enumi}{3}
\item \label{O4*} For all $q\leq\br(a)$, there exists a unique $v\leq a$ with $\br(v)=q$.
\end{enumerate}
(Here we have $q\leq\bd(a^*)$, and if $w$ is the unique element with $w\leq a^*$ with $\bd(w)=q$, then we take $v=w^*$.)

The elements $u$ and $v$ in \ref{O4} and \ref{O4*} are denoted $u={}_p\corest a$ and $v=a\rest_q$, respectively, and called the \emph{left restriction of $a$ to $p$} and the \emph{right restriction of $a$ to $q$}.  Some authors call ${}_p\corest a$ a \emph{restriction}, and $a\rest_q$ a \emph{co-restriction}; we prefer the left/right terminology, however, as it does not `prioritise' one over the other.

An \emph{ordered groupoid} is a groupoid with a partial order satisfying \ref{O1}--\ref{O4}.  In fact, when~$\CC$ is a groupoid, \ref{O2}, \ref{O3} and \ref{I4} together imply \ref{O1}.
\end{defn}

It is easy to see that the object set $v\CC$ is an \emph{order ideal} in any ordered $*$-category $\CC$, meaning that the following holds:
\bit
\item For all $a\in\CC$ and $p\in v\CC$, $a\leq p \implies a\in v\CC$.
\eit
Indeed, if $a\leq p$, and if we write $q=\bd(a)\in v\CC$, then by \ref{O1} we have $q\leq\bd(p)=p$.  But then $a,q\leq p$ and $\bd(a)=q=\bd(q)$, so by uniqueness in \ref{O4} we have $a=q\in v\CC$.  It also follows from this that ${}_q\corest p=q$ for any $q\leq p$.  We typically use facts such as these without explicit reference.

In what follows, it is typically more convenient to construct an ordered $*$-category $\CC$ by:
\bit
\item defining an order $\leq$ on the object set $v\CC$,
\item defining left restrictions ${}_p\corest a$, for $a\in\CC$ and $p\leq\bd(a)$, 
\item specifying that $a\leq b$ (for morphisms $a,b\in\CC$) when $a$ is a restriction of $b$.
\eit
The next lemma axiomatises the conditions required to ensure that we do indeed obtain an ordered $*$-category in this way.

\begin{lemma}\label{lem:C}
Suppose $\CC$ is a $*$-category for which the following two conditions hold:
\ben
\item \label{*C1} There is a partial order $\leq$ on the object set $v\CC$.
\item \label{*C2} For all $a\in\CC$, and for all $p\leq\bd(a)$, there exists a morphism ${}_p\corest a\in\CC$, such that the following hold, for all $a,b\in\CC$ and $p,q\in v\CC$:
\begin{enumerate}[label=\textup{\textsf{(O\arabic*)$'$}},leftmargin=10mm]
\item \label{O1'} If $p\leq\bd(a)$, then $\bd({}_p\corest a)=p$ and $\br({}_p\corest a)\leq\br(a)$.
\item \label{O2'} If $p\leq\bd(a)$, and if $q=\br({}_p\corest a)$, then $({}_p\corest a)^*={}_q\corest a^*$.
\item \label{O3'} ${}_{\bd(a)}\corest a = a$.
\item \label{O4'} For all $p\leq q\leq\bd(a)$, we have ${}_p\corest{}_q\corest a = {}_p\corest a$.
\item \label{O5'} If $p\leq\bd(a)$ and $\br(a)=\bd(b)$, and if $q=\br({}_p\corest a)$, then ${}_p\corest(a\circ b) = {}_p\corest a\circ{}_q\corest b$.
\end{enumerate}
\een
Then $\CC$ is an ordered $*$-category with order given by
\begin{equation}\label{eq:aleqb}
a\leq b \IFF a={}_p\corest b \qquad\text{for some $p\leq\bd(b)$.}
\end{equation}
Moreover, any ordered $*$-category has the above form.
\end{lemma}

\pf
Beginning with the final assertion, suppose $\CC$ is an ordered $*$-category.  Then $v\CC$ is a sub-poset of $\CC$, and hence~\ref{*C1} holds.  For \ref{*C2}, we take ${}_p\corest a$ to be the morphism $u\leq a$ from~\ref{O4}, and~\ref{O1'}--\ref{O5'} are all easily checked.  For example, to verify \ref{O2'}, suppose $a\in\CC$ and ${p\leq\bd(a)}$, and let $q=\br({}_p\corest a)$.  Then since ${}_p\corest a\leq a$, it follows from \ref{O2} that $({}_p\corest a)^*\leq a^*$, and we have $\bd(({}_p\corest a)^*)=\br({}_p\corest a)=q$.  But ${}_q\corest a^*$ is the unique element below $a^*$ with domain $q$, so in fact~${({}_p\corest a)^*={}_q\corest a^*}$.

Conversely, suppose conditions \ref{*C1} and \ref{*C2} both hold.  
We first check that the relation $\leq$ in~\eqref{eq:aleqb} is a partial order.  Reflexivity follows immediately from \ref{O3'}, and transitivity from~\ref{O4'}.  For anti-symmetry, suppose $a\leq b$ and $b\leq a$, so that $a={}_p\corest b$ and $b={}_q\corest a$ for some $p\leq\bd(b)$ and $q\leq\bd(a)$.  Then
\[
a = {}_p\corest b = {}_p\corest {}_q\corest a = {}_p\corest a \Implies \bd(a) = \bd({}_p\corest a) = p \leq \bd(b) \ANDSIM \bd(b) = q \leq \bd(a).
\]
It follows that $p=q=\bd(a)=\bd(b)$, and so $a = {}_p\corest b = {}_{\bd(b)}\corest b = b$.  Now that we know $\leq$ is a partial order, we verify conditions \ref{O1}--\ref{O4}.

\pfitem{\ref{O1} and \ref{O2}}  Suppose $a\leq b$, so that $a={}_p\corest b$ for some $p\leq\bd(b)$.  Then \ref{O1'} gives
\[
\bd(a)=\bd({}_p\corest b) = p\leq\bd(b) \AND \br(a)=\br({}_p\corest b)\leq\br(b).
\]
We also have $a^* = ({}_p\corest b)^* = {}_q\corest b^*$ by \ref{O2'}, where $q=\br({}_p\corest b)$, so that $a^*\leq b^*$.

\pfitem{\ref{O3}}  Suppose $a\leq b$ and $c\leq d$ are such that $\br(a)=\bd(c)$ and $\br(b)=\bd(d)$.  So $a={}_p\corest b$ and $c={}_q\corest d$ for some $p\leq\bd(b)$ and $q\leq\bd(d)$.  Since
\[
q = \bd({}_q\corest d) = \bd(c) = \br(a) = \br({}_p\corest b),
\]
it follows from \ref{O5'} that $a\circ c = {}_p\corest b\circ {}_q\corest d = {}_p\corest(b\circ d)$, so that $a\circ c \leq b\circ d$.

\pfitem{\ref{O4}}  Given $p\leq\bd(a)$, we certainly have $u\leq a$ and $\bd(u)=p$, where $u={}_p\corest a$.  For uniqueness, suppose also that $x\leq a$ for some $x\in\CC$ with $\bd(x)=p$.  Since $x\leq a$, we have $x={}_q\corest a$ for some $q\leq\bd(a)$.  But then $q=\bd({}_q\corest a) = \bd(x) = p$, so in fact $x={}_q\corest a={}_p\corest a=u$.
\epf

\begin{rem}\label{rem:dual}
The previous result referred only to left restrictions.  Right restrictions can be defined in terms of left restrictions and the involution:
\[
a\rest_q = ({}_q\corest a^*)^* \qquad\text{for $a\in\CC$ and $q\leq\br(a)$.}
\]
As explained in Definition \ref{defn:OC}, $a\rest_q$ is the unique element below $a$ with codomain $q$.  In particular, the following holds in any ordered $*$-category:
\begin{enumerate}[label=\textup{\textsf{(O\arabic*)$'$}},leftmargin=10mm]\addtocounter{enumi}{5}
\item \label{O6'} If $a\in\CC$ and $p\leq\bd(a)$, and if $q=\br({}_p\corest a)$, then ${}_p\corest a = a\rest_q$.
\end{enumerate}
Each of \ref{O1'}--\ref{O6'} of course have duals.  For example, the dual of \ref{O5'} says:
\bit
\item If $q\leq\br(b)$ and $\br(a)=\bd(b)$, and if $p=\bd(b\rest_q)$, then $(a\circ b)\rest_q=a\rest_p\circ b\rest_q$.
\eit
It is also worth noting that for $a,b\in\CC$ we have
\begin{align*}
a\leq b &\iff a={}_p\corest b &&\text{for some $p\leq\bd(b)$}\\
&\iff a=b\rest_q &&\text{for some $q\leq\br(b)$}\\
&\iff a={}_p\corest b=b\rest_q &&\text{for some $p\leq\bd(b)$ and $q\leq\br(b)$.}
\end{align*}
The $p,q\in P$ here are of course $p=\bd(a)$ and $q=\br(a)$.
\end{rem}

\begin{defn}\label{defn:cong}
A \emph{$v$-congruence} on a small category $\CC$ is an equivalence relation $\approx$ on $\CC$ satisfying the following, for all $a,b,u,v\in\CC$:
\begin{enumerate}[label=\textup{\textsf{(C\arabic*)}},leftmargin=10mm]
\item \label{C1} $a\approx b \implies [\bd(a)=\bd(b)$ and $\br(a)=\br(b)]$,
\item \label{C2} $a\approx b \implies u\circ a\approx u\circ b$, whenever the stated compositions are defined,
\item \label{C3} $a\approx b \implies a\circ v\approx b\circ v$, whenever the stated compositions are defined.
\end{enumerate}
If $\CC$ is a $*$-category, we say that $\approx$ is a \emph{$*$-congruence} if it satisfies \ref{C1}--\ref{C3} and the following, for all $a,b\in\CC$:
\begin{enumerate}[label=\textup{\textsf{(C\arabic*)}},leftmargin=10mm]\addtocounter{enumi}{3}
\item \label{C4} $a\approx b \implies a^*\approx b^*$.
\end{enumerate}
If $\CC$ is an ordered $*$-category, we say that $\approx$ is an \emph{ordered $*$-congruence} if it satisfies \ref{C1}--\ref{C4} and the following, for all $a,b\in\CC$:
\begin{enumerate}[label=\textup{\textsf{(C\arabic*)}},leftmargin=10mm]\addtocounter{enumi}{4}
\item \label{C5} $a\approx b \implies {}_p\corest a \approx {}_p\corest b$ for all $p\leq\bd(a)$.
\end{enumerate}
\end{defn}

Given a $v$-congruence $\approx$ on a small category $\CC$, the quotient category $\CC/{\approx}$ consists of all $\approx$-classes, under the induced composition.  We typically write $[a]$ for the $\approx$-class of $a\in\CC$.  It follows immediately from \ref{C1} that $p\approx q\implies p=q$ for objects $p,q\in v\CC$, so we can identify the object sets of $\CC$ and $\CC/{\approx}$, \emph{viz.}~$p\equiv[p]$.  In this way, for $a\in\CC$ we have $\bd[a]=\bd(a)$ and $\br[a]=\br(a)$ for all $a\in\CC$, and
\[
[a]\circ[b]=[a\circ b] \qquad\text{whenever $\br[a]=\bd[b]$.}
\]

\begin{lemma}\label{lem:approx}
\ben
\item \label{approx1} If $\approx$ is a $*$-congruence on a $*$-category $\CC$, then $\CC/{\approx}$ is a $*$-category, with involution given by
\begin{equation}
\label{eq:a*}
[a]^* = [a^*] \qquad\text{for all $a\in\CC$.}
\end{equation}
If also $a\circ a^*\approx\bd(a)$ for all $a\in\CC$, then $\CC/{\approx}$ is a groupoid.  
\item \label{approx2} If $\approx$ is an ordered $*$-congruence on an ordered $*$-category $\CC$, then $\CC/{\approx}$ is an ordered $*$-category, with involution \eqref{eq:a*}, and order given by
\begin{equation}\label{eq:alleqbe}
\al\leq\be \Iff  a\leq b \qquad\text{for some $a\in\al$ and $b\in\be$.}
\end{equation}
\een
\end{lemma}

\pf
Part \ref{approx1} is routine, and \ref{approx2} only a little more difficult, so we just give a sketch for the latter.  For this, one begins by showing that $\al\leq\be$ is equivalent to the ostensibly stronger condition:
\bit
\item For all $b\in\be$, there exists $a\in\al$ such that $a\leq b$.
\eit
It is then easy to check that $\leq$ is indeed a partial order on $\CC/{\approx}$.  Conditions \ref{O1}--\ref{O4} for $\CC/{\approx}$ follow quickly from the corresponding conditions for $\CC$.
\epf

We will shortly return to congruences in Lemma \ref{lem:Om} below, where we give simpler criteria for checking that a $v$-congruence satisfies \ref{C4} or \ref{C5}.  The proof of the lemma will require certain maps defined on any ordered $*$-category.  These maps will also be used extensively in the rest of the paper.  To define them, fix some such ordered $*$-category $\CC$ with object set $P=v\CC$.  For $p\in P$, we write
\[
p^\da = \set{q\in P}{q\leq p} 
\]
for the down-set of $p$ in the poset $(P,\leq)$.  Consider a morphism $a\in\CC$.  A restriction ${}_p\corest a$ is defined precisely when $p\leq\bd(a)$, and we write $p\vt_a = \br({}_p\corest a)$ for the codomain of ${}_p\corest a$; by \ref{O1'}, we have $p\vt_a\leq\br(a)$.  In other words, we have a map
\begin{equation}\label{eq:vta}
\vt_a : \bd(a)^\da\to\br(a)^\da \GIVENBY p\vt_a = \br({}_p\corest a).
\end{equation}
Note that for $q\leq\br(a)$, we have
\begin{equation}\label{eq:vta*}
\bd(a\rest_q) = \bd(({}_q\corest a^*)^*) = \br({}_q\corest a^*) = q\vt_{a^*}.
\end{equation}

\begin{lemma}\label{lem:vtavta*}
If $\CC$ is an ordered $*$-category, then for any $a\in\CC$,
\[
\vt_a:\bd(a)^\da\to\br(a)^\da \AND \vt_{a^*}:\br(a)^\da\to\bd(a)^\da
\]
are mutually inverse order-isomorphisms.
\end{lemma}

\pf
By symmetry, it is enough to show that 
\ben
\item \label{vtavta*1} $\vt_a$ is order-preserving, in the sense that $p\leq q \implies p\vt_a\leq q\vt_a$ for all $p,q\in\bd(a)^\da$, and
\item \label{vtavta*2} $\vt_a\vt_{a^*}=\id_{\bd(a)^\da}$.  
\een

\pfitem{\ref{vtavta*1}}  If $p\leq q\leq\bd(a)$, then ${}_p\corest a={}_p\corest{}_q\corest a \leq {}_q\corest a$, and so $\br({}_p\corest a)\leq\br({}_q\corest a)$, i.e.~$p\vt_a\leq q\vt_a$.

\pfitem{\ref{vtavta*2}}  Let $p\leq\bd(a)$, and write $q=p\vt_a=\br({}_p\corest a)$.  By \ref{O6'} we have ${}_p\corest a=a\rest_q$.  Combining this with \ref{O1'} and \eqref{eq:vta*}, we deduce that $p = \bd({}_p\corest a) = \bd(a\rest_q) = q\vt_{a^*} = p\vt_a\vt_{a^*}$, as required.
\epf

Note that property \ref{O5'} and its dual (cf.~Remark \ref{rem:dual}) can be rephrased in terms of the $\vt$ maps.  Specifically, if $a,b\in\CC$ (an ordered $*$-category) are such that $\br(a)=\bd(b)$, then (keeping~\eqref{eq:vta*} and Lemma \ref{lem:vtavta*} in mind) we have
\begin{equation}\label{eq:pab}
{}_p\corest(a\circ b) = {}_p\corest a \circ {}_{p\vt_a}\corest b \AND (a\circ b)\rest_q = a\rest _{q\vt_b^{-1}} \circ b\rest_q \qquad\text{for all $p\leq\bd(a)$ and $q\leq\br(b)$.}
\end{equation}
This can be extended to restrictions of compositions with an arbitrary number of terms.  In a sense, the first part of the next result is a formalisation of the previous observation.
Throughout the paper we write $f|_A$ for the (set-theoretic) restriction of a function $f$ to a subset $A$ of its domain.

\begin{lemma}\label{lem:vt}
Let $\CC$ be an ordered $*$-category, and let $a,b\in\CC$ and $p\in v\CC$.
\ben
\item \label{vt1} If $\br(a)=\bd(b)$, then $\vt_{a\circ b} = \vt_a\vt_b$.
\item \label{vt2} $\vt_p=\id_{p^\da}$.
\item \label{vt3} If $p\leq\bd(a)$, then $\vt_{{}_p\corest a} = \vt_a|_{p^\da}$.  
\een
\end{lemma}

\pf
\firstpfitem{\ref{vt1}}  Write $p=\bd(a)$, $q=\br(a)=\bd(b)$ and $r=\br(b)$.  Note that
\[
\vt_a:p^\da\to q^\da \COMMA \vt_b:q^\da\to r^\da \AND \vt_{a\circ b}:p^\da\to r^\da,
\]
so that $\vt_{a\circ b}$ and $\vt_a\vt_b$ are defined on the same domain, i.e.~on $p^\da$.  Now let $s\leq p$, and put $t=s\vt_a=\br({}_s\corest a)$.  Then by \ref{O5'} we have
\[
s\vt_{a\circ b} = \br({}_s\corest(a\circ b)) = \br({}_s\corest a\circ{}_t\corest b) = \br({}_t\corest b) = t\vt_b = s\vt_a\vt_b.  
\]

\pfitem{\ref{vt2}}  By part \ref{vt1}, we have $\vt_p = \vt_{p\circ p} = \vt_p\vt_p$, and the result follows since the only idempotent bijection $p^\da\to p^\da$ is the identity.

\pfitem{\ref{vt3}}  Both maps have domain $p^\da$, so it suffices to show that $t\vt_{{}_p\corest a} = t\vt_a$ for all $t\leq p$.  For this we use \ref{O4'} to calculate
\[
t\vt_{{}_p\corest a} = \br({}_t\corest{}_p\corest a) = \br({}_t\corest a) = t\vt_a.  \qedhere
\]
\epf

%We will also need the following simple property of the $\vt$ maps:
%
%\begin{lemma}\label{lem:vtaOP}
%If $\CC$ is an ordered $*$-category, and if $a\in\CC$, then $\vt_a$ is order-preserving, in the sense that
%\[
%p\leq q \Implies p\vt_a\leq q\vt_a \qquad\text{for all $p,q\in\bd(a)^\da$.}
%\]
%\end{lemma}
%
%\pf
%If $p\leq q\leq\bd(a)$, then ${}_p\corest a={}_p\corest{}_q\corest a \leq {}_q\corest a$, and so $\br({}_p\corest a)\leq\br({}_q\corest a)$, i.e.~$p\vt_a\leq q\vt_a$.
%\epf

In later sections we will define $v$-congruences by specifying sets of generating pairs:

\begin{defn}\label{defn:Om}
Suppose $\CC$ is a small category, and $\Om\sub\CC\times\CC$ a set of pairs satisfying
\bit
\item $(u,v)\in\Om \implies [\bd(u)=\bd(v)$ and $\br(u)=\br(v)]$.
\eit
The \emph{($v$-)congruence generated by $\Om$}, denoted $\Om^\sharp$, is the least congruence on $\CC$ containing $\Om$.  Specifically, we have $(a,b)\in\Om^\sharp$ if and only if there is a sequence
\[
a = a_1 \to \cdots \to a_k = b
\]
such that for each $1\leq i<k$ we have
\[
a_i = b_i\circ u_i\circ c_i \AND a_{i+1} = b_i\circ v_i\circ c_i \qquad\text{for some $b_i,c_i\in\CC$ and $(u_i,v_i)\in\Om\cup\Om^{-1}$.}
\]
(Here $\Om^{-1} = \set{(v,u)}{(u,v)\in\Om}$.)
\end{defn}

The next result shows how conditions \ref{C4} and \ref{C5} regarding a $v$-congruence ${\approx}=\Om^\sharp$ can be deduced from the corresponding properties of its generating set $\Om$.

\begin{lemma}\label{lem:Om}
Suppose ${\approx}=\Om^\sharp$ is a $v$-congruence on an ordered $*$-category $\CC$.
\ben
\item If $\Om$ satisfies the condition
\begin{equation}\label{eq:Om*}
(a,b)\in\Om \Implies a^* \approx b^* ,
\end{equation}
then $\approx$ satisfies condition \ref{C4}.
\item If $\Om$ satisfies the condition
\begin{equation}\label{eq:Om}
%(a,b)\in\Om \Implies \vt_a=\vt_b \ANd {}_p\corest a\approx{}_p\corest b \text{ for all $p\leq\bd(a)$,}
(a,b)\in\Om \Implies  {}_p\corest a\approx{}_p\corest b \text{ for all $p\leq\bd(a)$,}
\end{equation}
then $\approx$ satisfies condition \ref{C5}.
\een
\end{lemma}

\pf
We just prove the second part, as the first is easier.  We first claim that
\begin{equation}\label{eq:Om'}
\vt_a=\vt_b \qquad\text{for all $(a,b)\in\Om$.}
\end{equation}
Indeed, for all $(a,b)\in\Om$ we have ${}_p\corest a\approx{}_p\corest b$ for all $p\leq\bd(a)$, by \eqref{eq:Om}.  Since $\approx$ is a $v$-congruence, it follows that $\br({}_p\corest a) = \br({}_p\corest b)$, i.e.~$p\vt_a=p\vt_b$, for all such $p$, and this says precisely that $\vt_a=\vt_b$.

To show that \ref{C5} holds, suppose $a\approx b$, and let the $a_i,b_i,c_i,u_i,v_i$ be as in Definition \ref{defn:Om}.  Of course it suffices to show that ${}_p\corest a_i \approx {}_p\corest a_{i+1}$ for each $i$.  For this we define
\[
q = p\vt_{b_i} \AND r = q\vt_{u_i}.
\]
Since $\vt_{u_i}=\vt_{v_i}$ by \eqref{eq:Om'}, we also have $r = q\vt_{v_i}$.  We then use \ref{O5'} and \eqref{eq:Om} to calculate
\[
{}_p\corest a_i = {}_p\corest(b_i\circ u_i\circ c_i) = {}_p\corest b_i \circ {}_q\corest u_i \circ {}_r\corest c_i \approx {}_p\corest b_i \circ {}_q\corest v_i \circ {}_r\corest c_i = {}_p\corest(b_i\circ v_i\circ c_i) = {}_p\corest a_{i+1}.  \qedhere
\]
\epf

\section{\boldmath Regular $*$-semigroups}\label{sect:RSS}

We now arrive at the main subject of our study, the class of regular $*$-semigroups.  In Section~\ref{subsect:Sprelim} we recall the relevant definitions, and give some basic results.  In Section \ref{subsect:StoG} we construct a functor into the category of ordered groupoids; see Proposition \ref{prop:calGfunctor}.  We then show how the abstract theory applies in a number of concrete examples: adjacency semigroups and Rees 0-matrix semigroups in Section \ref{subsect:Seg}, and diagram monoids in Section \ref{subsect:D}.  We will return to these examples often throughout the paper.

\subsection{Definitions and basic properties}\label{subsect:Sprelim}

Here is the main definition, given by Nordahl and Scheiblich in \cite{NS1978}:

\begin{defn}\label{defn:RSS}
A \emph{regular $*$-semigroup} is a semigroup $S$ with a unary operation ${{}^*:S\to S:a\mt a^*}$ satisfying
\[
(a^*)^* = a = aa^*a \AND (ab)^* = b^*a^* \qquad\text{for all $a,b\in S$.}
\]
We denote by $\RSS$ the category of all regular $*$-semigroups.  A morphism in $\RSS$ is a \emph{$*$-semigroup homomorphism} (a.k.a.~\emph{$*$-morphism}), i.e.~a map $\phi:S\to S'$ (for regular $*$-semigroups $S$ and $S'$) satisfying
\[
(ab)\phi = (a\phi)(b\phi) \AND (a^*)\phi = (a\phi)^* \qquad\text{for all $a,b\in S$,}
\]
where for simplicity we have written ${}^*$ for the unary operation on both $S$ and $S'$.
\end{defn}

Note that a regular $*$-semigroup $S$ is considered to be an algebra of type $(2,1)$, i.e.~${}^*$ is a basic operation of $S$.  It is also worth noting that a semigroup homomorphism between regular $*$-semigroups need not be a $*$-morphism.  For example, if $e\in S$ satisfies $e^2=e\not=e^*$, then the constant map with image $e$ is a semigroup homomorphism $S\to S$, but not a $*$-morphism.  This means that $\RSS$ is not a \emph{full} subcategory of the category of all semigroups.

Any inverse semigroup is a regular $*$-semigroup, with the involution being inversion, ${a^*=a^{-1}}$.  In fact, it has long been known \cite{Schein1963} that inverse semigroups are precisely the regular $*$-semigroups satisfying the additional identity
\begin{equation}\label{eq:aa*bb*}
aa^*bb^* = bb^*aa^* \qquad\text{for all $a,b\in S$.}
\end{equation}
As we will see, inverse semigroups represent a somewhat `degenerate' case in the theory we develop; see Section \ref{sect:I}.

For a regular $*$-semigroup $S$, we write
\[
E = E(S) = \set{e\in S}{e^2=e} \AND P = P(S) = \set{p\in S}{p^2=p=p^*}.
\]
So $E$ consists of all idempotents of $S$.  The elements of $P$ are called \emph{projections}, and of course we have $P\sub E$. Projections play a very important role in virtually every study of regular $*$-semigroups in the literature, and this is also true in the current work.  
The following result gathers some of the basic properties of idempotents and projections.  Proofs can be found in many places (e.g., \cite{Imaoka1980}), but we give some simple arguments here to keep the paper self-contained.

\begin{lemma}\label{lem:PS1}
If $S$ is a regular $*$-semigroup, with sets of idempotents $E=E(S)$ and projections $P=P(S)$, then
\ben
\item \label{PS11} $P=\set{aa^*}{a\in S}=\set{a^*a}{a\in S}$,
\item \label{PS12} $E=P^2=\set{pq}{p,q\in P}$, and consequently $\la E\ra = \la P\ra$,
\item \label{PS13} $a^*Pa\sub P$ for all $a\in S$.  
\een
\end{lemma}

\pf
\firstpfitem{\ref{PS11}}  It is enough to show that $P=\set{aa^*}{a\in S}$.  If $p\in P$, then $p=pp=pp^*$.  Conversely, if $a\in S$, then $(aa^*)^2=aa^*aa^*=aa^*$, and $(aa^*)^*=(a^*)^*a^*=aa^*$.

\pfitem{\ref{PS12}}  If $e\in E$, then $e=ee^*e=e(ee)^*e=ee^*\cdot e^*e$, with $ee^*,e^*e\in P$, by part \ref{PS11}.  Conversely, for any $p,q\in P$, we have $pq=pq(pq)^*pq=pqq^*p^*pq=pqqppq=pqpq=(pq)^2$.

\pfitem{\ref{PS13}}  If $p\in P$, then $a^*pa=a^*p^*pa=(pa)^*pa\in P$ by part \ref{PS11}.
\epf

Comparing \eqref{eq:aa*bb*} with Lemma \ref{lem:PS1}\ref{PS11}, it follows that inverse semigroups are precisely the regular $*$-semigroups with commuting projections.  In fact, it is easy to see that any idempotent of an inverse semigroup is a projection, so another equivalent condition for a regular $*$-semigroup~$S$ to be inverse is that $E(S)=P(S)$.  Some other conditions will be listed below in Lemma \ref{lem:inv}.

Consider a regular $*$-semigroup $S$, and let $P=P(S)$.  If $p,q\in P$, then it follows from Lemma~\ref{lem:PS1}\ref{PS13} that $pqp = p^*qp\in P$.  Thus, for all $p\in P$, we have a well-defined map
%\begin{equation}\label{eq:thpS}
\[
\th_p:P\to P \GIVENBY q\th_p = pqp \qquad\text{for $q\in P$.}
\]
%\end{equation}
The next result summarises some of the key properties of these maps.

\begin{lemma}\label{lem:PS2}
If $S$ is a regular $*$-semigroup, then for any $p,q\in P = P(S)$ we have
\[
p\theta_p = p \COMMA 
\theta_p\theta_p=\theta_p \COMMA 
p\theta_q\theta_p=q\theta_p \COMMA 
\theta_p\theta_q\theta_p=\theta_{q\theta_p} \AND
\theta_p\theta_q\theta_p\theta_q=\theta_p\theta_q.
\]
\end{lemma}

\pf
The first two follow from the fact that projections are idempotents.
The third and fifth follow from the fact that the product of two projections is an idempotent (cf.~Lemma~\ref{lem:PS1}\ref{PS12}).  
For the fourth, given any $t\in P$ we have
\[
t\th_{q\th_p} = t\th_{pqp} = pqp\cdot t\cdot pqp = t\th_p\th_q\th_p. \qedhere
\]
\epf

\begin{rem}\label{rem:Bana1}
In \cite{Banaschewski1979}, Banaschewski studied the special class of regular $*$-semigroups satisfying the additional identity $abb^*a^* = aa^*$.  It is not hard to see that these are precisely the regular $*$-semigroups for which each $\th_p$ is a constant map (with image $p$).  Indeed, for $p,q\in P$, we take $a=p$ and $b=q$ in the above identity to see that $q\th_p = pqp = pqq^*p^* = pp^* = p$.  Conversely, given any $a,b\in S$, if $\th_{a^*a}$ is constant, then
\[
abb^*a^* = aa^*abb^*a^*aa^* = a\cdot (bb^*)\th_{a^*a}\cdot a^* = a\cdot a^*a\cdot a^* = aa^*.
\]
We will say more about these semigroups in Remark \ref{rem:Bana2}.
\end{rem}

Next we define a relation $\leq$ on $P=P(S)$  by
\begin{equation}\label{eq:leq}
p\leq q \Iff p=p\th_q = qpq \qquad\text{for $p,q\in P$}.
\end{equation}
Note that $p\leq q$ precisely when $q$ is a left and right identity for $p$.  Using this it is easy to see that~$\leq$ is a partial order.  It is also possible to deduce this fact from the identities listed in Lemma~\ref{lem:PS2}, as we do in Lemma \ref{lem:leqP} below.  In fact, since $p=qpq \iff p=pq=qp$, and since
%\begin{equation}\label{eq:p=pq}
\[
p=pq \implies p=p^*=(pq)^*=q^*p^*=qp,
\]
%\end{equation}
and conversely by symmetry, it follows that
\begin{equation}\label{eq:leq2}
p\leq q \Iff p=pq \Iff p=qp.
\end{equation}
The partial order on projections can be used to verify the following neat equational (rather than existential) characterisation of Green's relations on a regular $*$-semigroup.  It is again well known, but we provide a simple proof for completeness.

\begin{lemma}\label{lem:Green}
If $S$ is a regular $*$-semigroup, and if $a,b\in S$, then
\[
a\R b \iff aa^*=bb^* \AND a\L b \iff a^*a=b^*b.  
\]
\end{lemma}

\pf
We give the argument for $\R$, with $\L$ being dual.
Suppose first that $a\R b$, so that $a=bx$ for some $x\in S$.  Using \eqref{eq:leq2}, it then follows that 
\[
aa^*=bxa^*=bb^*bxa^*=bb^*aa^* \Implies aa^*\leq bb^*.
\]
By symmetry, we also have $bb^*\leq aa^*$, so that $aa^*=bb^*$.

Conversely, if $aa^*=bb^*$, then $a=aa^*a=b(b^*a)$, and similarly $b=a(a^*b)$, so that $a\R b$.
\epf

The above properties of projections, and their associated $\th$ mappings, will be of fundamental importance in all that follows.  In particular, even though $P=P(S)$ is generally not a subsemigroup of the regular $*$-semigroup $S$, we can regard it as a \emph{unary algebra}, whose basic (unary) operations are the $\th_p$ ($p\in P$).  In Section \ref{sect:P}, we take the properties from Lemma \ref{lem:PS2} as the axioms for what we will call a \emph{projection algebra} (see Definition \ref{defn:P}).  It has been known for some time that (abstract) projection algebras are precisely the unary algebras of projections of regular $*$-semigroups (see for example \cite{Imaoka1983}), but we will also verify this with a simple example in Section~\ref{subsect:PtoS}.  When $S$ is inverse, $P$ is simply a subsemigroup of~$S$, and there is no real advantage in considering $P$ as a unary algebra.  In a sense, a projection algebra is the `$*$-analogue' of the semilattice of idempotents of an inverse semigroup.
For this reason, it is convenient to record the following simple result characterising inverse semigroups in terms of properties of their projections.  %It is a straightforward consequence of commutativity of idempotents, which implies that $q\th_p=pq$ for projections $p,q$.

\begin{lemma}\label{lem:inv}
For a regular $*$-semigroup $S$, the following are equivalent:
\ben\bmc2
\item \label{inv1} $S$ is inverse,
\item \label{inv2} $E(S)=P(S)$,
\item \label{inv3} $P(S)$ is a subsemigroup of $S$,
\item \label{inv4} $pq=qp$ for all $p,q\in P(S)$,
\item \label{inv5} $p\th_q=q\th_p$ for all $p,q\in P(S)$.  
\item[] ~ %\epfres
\emc\een
\end{lemma}

\pf
This is straightforward, and mostly well known.  We just give a proof that \ref{inv5}$\implies$\ref{inv3}.  For this, we first note that \ref{inv5} says $pqp=qpq$ for any $p,q\in P(S)$.  Combining this with Lemma~\ref{lem:PS1}\ref{PS12}, it then follows that
\[
pq = (pq)^2 = p(qpq) 
%= p(p\th_q) = p(q\th_p) 
= p(pqp) = pqp \in P.
%\IMPLIES qp = (pq)^* = (pqp)^* = pqp = pq.  
\qedhere
\]
\epf

%\pf
%\firstpfitem{\ref{inv1}$\implies$\ref{inv2}}  
%
%\pfitem{\ref{inv2}$\implies$\ref{inv3}}  
%
%\pfitem{\ref{inv3}$\implies$\ref{inv4}}  
%
%\pfitem{\ref{inv4}$\implies$\ref{inv5}}  
%
%\pfitem{\ref{inv5}$\implies$\ref{inv1}}  
%
%\epf

%\begin{rem}
%When $S$ is inverse, commutativity of projections implies that $q\th_p=p\th_q=pq=qp$ for all~${p,q\in P}$.  This then leads to a significant simplification in the theory developed below, as we will explore in detail in Section \ref{sect:I}.  \JE{Repetitive?}
%\end{rem}

\begin{rem}\label{rem:inv}
One might wonder why the following condition was not included in Lemma \ref{lem:inv}:
\bit
\item the poset $(P(S),{\leq})$ is a $\wedge$-semilattice.
\eit
While this is a necessary condition for $S$ to be inverse, it is not sufficient, as we will see in Example \ref{eg:AG} below.
\end{rem}

As we have already mentioned, projections have played a crucial role in virtually all studies of regular $*$-semigroups.  Of particular significance in papers such as \cite{DEG2017,EG2017,JEgrpm} are pairs of projections ${p,q\in P=P(S)}$ that are mutual inverses, i.e.~pairs that satisfy $p=pqp$ and $q=qpq$.  This is very much the case in the current work: so much so, in fact, that we define relations $\leqF$ and $\F$ on~$P$ by
\begin{equation}\label{eq:leqFF}
p\leqF q \Iff p=pqp=q\th_p \AND p\F q \Iff[p\leqF q\text{ and }q\leqF p].
\end{equation}
We think of pairs of $\F$-related projections as \emph{friends} (hence the symbol $\F$).  It is clear that~$\leqF$ and $\F$ are both reflexive, and that $\F$ is symmetric.  But neither $\leqF$ nor $\F$ is transitive in general (and neither is friendship in `real life'); for some examples of this, see Sections~\ref{subsect:Seg} and~\ref{subsect:D}.  It is easy to see that $p\leq q \implies p\leqF q$, where $\leq$ is the partial order in~\eqref{eq:leq}.  It follows quickly from Lemma~\ref{lem:Green} that 
\begin{equation}\label{eq:pRpqLq}
p\F q \Iff p\R pq \L q \Iff p\L qp\R q  .
\end{equation}
It of course follows from this that
\begin{equation}\label{eq:FD}
p\F q \Implies p\D q,
\end{equation}
though the converse does not hold in general; again see Sections~\ref{subsect:Seg} and~\ref{subsect:D}.

\begin{rem}\label{rem:Bana2}
Following on from Remark \ref{rem:Bana1}, when the $\th$ maps in a regular $*$-semigroup $S$ are all constant, it follows that all projections are $\F$-related (cf.~\eqref{eq:leqFF}), and hence $\D$-related (cf.~\eqref{eq:FD}).  Consequently, such an $S$ is $\D$-simple.  This gives a new (simple) proof of \cite[Theorem~4.1]{Petrich1985}.
\end{rem}

Combining \eqref{eq:pRpqLq} with \cite[Theorem 3]{MC1956}, it follows that
\[
p\F q \Iff L_p\cap R_q \text{ contains an idempotent} \Iff R_p\cap L_q \text{ contains an idempotent} ,
\]
and of course these idempotents are
\[
qp \in L_p\cap R_q \AND pq \in R_p\cap L_q.
\]
This all means that a pair of $\F$-related projections $p,q\in P(S)$ generates a $2\times2$ rectangular band in $S$:
\[
\begin{tikzpicture}[scale=1]
\draw(0,0)--(2,0)--(2,2)--(0,2)--(0,0) (0,1)--(2,1) (1,0)--(1,2);
\node()at(.5,1.45){$p$};
\node()at(1.5,1.45){$pq$};
\node()at(.5,.45){$qp$};
\node()at(1.5,.45){$q$};
\end{tikzpicture}
\]
One can define a graph $\Ga(S)$ with vertex set $P=P(S)$, and with (undirected) edges $\{p,q\}$ whenever $p\not=q$ and $p\F q$.  Connectivity properties of (certain subgraphs of) these graphs played a vital role in \cite{EG2017}, in studies of idempotent-generated regular $*$-semigroups and their ideals.  We will say more about this in Section \ref{subsect:D}.

The (standard) proof of Lemma \ref{lem:PS1}\ref{PS12} given above involved showing that any idempotent $e\in E=E(S)$ of a regular $*$-semigroup $S$ satisfies $e = ee^*\cdot e^*e$, with $ee^*,e^*e\in P=P(S)$.  It quickly follows that the idempotents $e,e^*\in E(S)$ also generate a $2\times2$ rectangular band in $S$:
\[
\begin{tikzpicture}[scale=1]
\draw(0,0)--(2,0)--(2,2)--(0,2)--(0,0) (0,1)--(2,1) (1,0)--(1,2);
\node()at(.5,1.45){$\phantom{{}_*}e\phantom{{}^*}$};
\node()at(1.5,1.45){$ee^*_{\phantom*}$};
\node()at(.5,.45){$e^*_{\phantom*}e$};
\node()at(1.5,.45){$e^*_{\phantom*}$};
\end{tikzpicture}
\]
This is in fact a characterisation of the biordered sets arising from regular $*$-semigroups, in a sense made precise in \cite[Corollary 2.7]{NP1985}.  In any case, it quickly follows that $ee^*\F e^*e$ for any $e\in E$, so this shows that any idempotent is a product of two $\F$-related projections, and it follows that not only do we have $E=P^2$ (cf.~Lemma \ref{lem:PS1}\ref{PS12}), but in fact
\[
E = \set{pq}{(p,q)\in{\F}}.
\]
Actually, such factorisations for idempotents are \emph{unique}.  Indeed, if $e=pq$ with $(p,q)\in{\F}$, then
\begin{equation}\label{eq:ee*e*e}
ee^* = pq(pq)^* = pqq^*p^* = pqp = p \ANDSIM e^*e = q.
\end{equation}
Consequently, we have $|E|=|{\F}|$.  A number of studies have enumerated the idempotents in various classes of regular $*$-semigroups \cite{DEEFHHLM2019,DEEFHHL2015,Larsson2006}.  The identity $|E|=|{\F}|$ played an implicit role in~\cite{DEEFHHLM2019}.

So idempotents are products of $\F$-related projections.  As an application of the theory developed in this paper, we will be able to prove the following result, which extends this to arbitrary products of idempotents, though uniqueness is lost, in general (see for example~\eqref{eq:triangle}); see Proposition~\ref{prop:ES}\ref{ES3} and Remark~\ref{rem:ES}.   The result might be known, but we are not aware of its existence in the literature; it does bear some resemblance, however, to a classical result of FitzGerald \cite{FitzGerald1972} on products of idempotents in regular semigroups.  Moreover, we believe it could be useful in studies of idempotent-generated regular $*$-semigroups, and might have even led to some simpler arguments in existing studies such as~\cite{EF2012,EG2017,DEG2017}.

\begin{prop}\label{prop:ERSS}
Let $S$ be a regular $*$-semigroup, with $E=E(S)$ and $P=P(S)$.  Then for any $a\in\la E\ra=\la P\ra$ we have
\[
a = p_1\cdots p_k \qquad\text{for some $p_1,\ldots,p_k\in P$ with $p_1\F\cdots\F p_k$.}
\]
\end{prop}

\subsection[From regular $*$-semigroups to ordered groupoids]{\boldmath From regular $*$-semigroups to ordered groupoids}\label{subsect:StoG}

We now show how to use the projections of a regular $*$-semigroup $S$ to construct an ordered groupoid $\G=\G(S)$.  The groupoid~$\G$ has the same underlying set as $S$, and the (partially defined) composition in $\G$ is a restriction of the (totally defined) product in $S$.  Roughly speaking,~$\G$ remembers only the `nice' or `easy' products; see Sections \ref{subsect:Seg} and \ref{subsect:D} for some examples justifying our use of these words.  More importantly, the construction leads to a functor from regular $*$-semigroups to ordered groupoids; see Proposition \ref{prop:calGfunctor}.

As it happens, the construction of $\G=\G(S)$ is \emph{exactly} the same as in the inverse case; see \cite[Chapter~4]{Lawson1998}.  However, we will see that the groupoid $\G(S)$ is not a \emph{total invariant} of the regular $*$-semigroup $S$, unlike the case for inverse semigroups; see Examples \ref{eg:AG} and \ref{eg:Rees}.

For a regular $*$-semigroup $S$, and for $a\in S$, we write
\[
\bd(a) = aa^* \AND \br(a) = a^*a.
\]
Both of these elements are projections (cf.~Lemma \ref{lem:PS1}\ref{PS11}), and the identity $a=aa^*a$ gives
\[
a = \bd(a)\cdot a = a\cdot \br(a).
\]
Moreover, Lemma \ref{lem:Green} says that
\[
a\R b \iff \bd(a)=\bd(b) \AND a\L b \iff \br(a)=\br(b).
\]
Since $p=pp^*=p^*p$ for any projection $p$, we have $\bd(p)=\br(p)=p$ for such $p$.  It then follows from Lemma \ref{lem:Green} that the $\R$- and $\L$-classes of $p$ are given by
\begin{equation}\label{eq:RpLp}
R_p = \set{a\in S}{aa^*=pp^*} = \set{a\in S}{\bd(a)=p} \ANDSIM L_p = \set{a\in S}{\br(a)=p}.
\end{equation}

\begin{lemma}\label{lem:dab}
If $S$ is a regular $*$-semigroup, and if $a,b\in S$ satisfy $\br(a)=\bd(b)$, then
\[
\bd(ab) = \bd(a) \AND \br(ab) = \br(b).
\]
\end{lemma}

\pf
We just prove the first claim, as the second is symmetrical.  If $\br(a)=\bd(b)$, then $a^*a = bb^*$, and then
\[
\bd(ab) = ab(ab)^* = abb^*a = aa^*aa^* = aa^* = \bd(a).  \qedhere
\]
\epf

This allows us to make the following definition.

\newpage

\begin{defn}\label{defn:GS}
Given a regular $*$-semigroup $S$, we define the (small) category $\G=\G(S)$ as follows.
\bit
\item The underlying (morphism) set of $\G$ is $S$.
\item The object set of $\G$ is $v\G=P=P(S)$.
\item For $a\in\G$ we have $\bd(a) = aa^*$ and $\br(a) = a^*a$.
\item For $a,b\in\G$ with $\br(a)=\bd(b)$, we have $a\circ b = ab$.
\eit
Since the underlying set of $\G$ is $S$, the unary operation ${}^*:S\to S$ can be thought of as a map ${}^*:\G\to\G$.  It is a routine matter to verify conditions \ref{I1}--\ref{I4} from Definition \ref{defn:*cat}, so that~$\G$ is a groupoid, with inversion given by ${}^{-1} = {}^*$:
\bit
\item For $a\in\G$ we have $a^{-1}=a^*$.
\eit
We will soon show that $\G$ is in fact an \emph{ordered} groupoid, for which we need to define suitable restrictions.  Let $\leq$ be the ordering on $P=v\G$ given in \eqref{eq:leq}.
\bit
\item For $a\in\G$ and $p\leq\bd(a)$ we define ${}_p\corest a = pa$.
\eit
\end{defn}

\begin{rem}\label{rem:GS0}
For projections $p,q\in P$, it follows from \eqref{eq:RpLp} that
\[
\G(p,q) = \set{a\in S}{\bd(a)=p,\ \br(a)=q} = R_p \cap L_q.
\]
Since ${\D}={\R}\circ{\L}$, such a morphism set is non-empty precisely when $p\D q$, and is then in fact an $\H$-class.  It follows that the connected components of $\G=\G(S)$ are in one-one correspondence with the $\D$-classes of $S$.
%
%It is also worth noting that if $p,q\in P$ are such that $p\leq q(=\bd(q))$, then ${}_p\corest q=p$, as quickly follows from \eqref{eq:leq2}.  In other words, $P$ is an order ideal of the poset $(\G,{\leq})$.
\end{rem}

%Since the underlying set of $\G=\G(S)$ is $S$, the unary operation ${}^*:S\to S$ can be thought of as a map ${}^*:\G\to\G$.  It is a routine matter to verify conditions \ref{I1}--\ref{I4} from Definition \ref{defn:*cat}, so that $\G$ is a groupoid, with inversion given by ${}^{-1} = {}^*$.  In fact, we will soon see that $\G$ is an \emph{ordered} groupoid.  To set this up, we make the following definition.
%
%
%\begin{defn}\label{defn:GS2}
%Given a regular $*$-semigroup $S$, let $\G=\G(S)$ be the groupoid from Definition~\ref{defn:GS}, and let $P=P(S)=v\G$.  Let $\leq$ be the ordering on $P=v\G$ given in \eqref{eq:leq}.  %Left restrictions in $\G$ are defined as follows.
%\bit
%\item For $a\in\G$ and $p\leq\bd(a)$ we define ${}_p\corest a = pa$.
%\eit
%\end{defn}

\begin{prop}\label{prop:GS}
If $S$ is a regular $*$-semigroup, then $\G=\G(S)$ is an ordered groupoid.
\end{prop}

\pf
We need to verify conditions \ref{O1'}--\ref{O5'} from Lemma \ref{lem:C}, with respect to the ordering~\eqref{eq:leq}, and the restrictions given in Definition \ref{defn:GS}.

\pfitem{\ref{O1'}}  Consider a morphism $a\in\G$, and let $p\leq \bd(a)=aa^*$.  It follows from $p\leq aa^*$ that $p=paa^*$, and so
\[
\bd({}_p\corest a) = \bd(pa) = pa(pa)^* = paa^*p^* = pp = p .
\]
We also have
\begin{equation}\label{eq:a*pa}
\br({}_p\corest a) = \br(pa) = (pa)^*pa = a^*p^*pa = a^*pa .
\end{equation}
It then follows that
\[
\br({}_p\corest a) = a^*pa = a^*a \cdot a^*pa\cdot a^*a = \br(a) \cdot \br({}_p\corest a) \cdot \br(a).
\]
By \eqref{eq:leq}, the previous conclusion says precisely that $\br({}_p\corest a) \leq \br(a)$.  

\pfitem{\ref{O2'}}  By \eqref{eq:a*pa} we have $q=a^*pa$, and again $p=paa^*$ follows from $p\leq\bd(a)$.  It follows that
\[
({}_p\corest a)^* = (pa)^* = a^*p^* = a^*p = a^*paa^* = qa^* = {}_q\corest a^*.
\]

\pfitem{\ref{O3'}}  Here ${}_{\bd(a)}\corest a = \bd(a)a = aa^*a = a$.

\pfitem{\ref{O4'}}  It follows from $p\leq q$ that $p=pq$, and so ${}_p\corest{}_q\corest a = p(qa) = (pq)a = pa = {}_p\corest a$.

\pfitem{\ref{O5'}}  Here we again have $q=a^*pa$.  Since $p\leq\bd(a)$ we have $p=aa^*p$, so $pa=aa^*pa=aq$.  It follows that
\[
{}_p\corest(a\circ b) = p(ab) = pp(ab) = paqb = {}_p\corest a \circ {}_q\corest b.  \qedhere
\]
\epf

As in Remark \ref{rem:dual}, right restrictions in $\G=\G(S)$ are given by $a\rest_q = ({}_q\corest a^{-1})^{-1}$, for $q\leq\br(a)$, and in fact we have
\[
a\rest_q = (qa^*)^* = aq
\]
for such $q$.  For $a,b\in\G$, we have
\begin{align}
\nonumber a\leq b &\IFF a={}_p\corest b &&\hspace{-1.5cm}\text{for some $p\leq\bd(b)$}\\
\nonumber &\IFF a=b\rest_q &&\hspace{-1.5cm}\text{for some $q\leq\br(b)$}\\
\label{eq:aleqb2} &\IFF a={}_p\corest b=b\rest_q &&\hspace{-1.5cm}\text{for some $p\leq\bd(b)$ and $q\leq\br(b)$,}
\end{align}
and then of course $p=\bd(a)$ and $q=\br(a)$.

\begin{rem}
In particular, the relation $\leq$ on the regular $*$-semigroup $S(=\G)$ given in \eqref{eq:aleqb2} is a partial order.  This order was in fact used by Imaoka in \cite{Imaoka1981}, but in a different form.  Imaoka's order, which for clarity we will denote by $\leq'$, was defined by
\begin{equation}\label{eq:leq'}
a\leq' b \Iff a\in Pb\cap bP \qquad\text{where $P=P(S)$.}
\end{equation}
In light of the rules ${}_p\corest b=pb$ and $b\rest_q=bq$, it is clear that $a\leq b\implies a\leq' b$.  For the converse, suppose $a\leq'b$, so that
\[
a = pb = bq \qquad\text{for some $p,q\in P$.}
\]
In the following calculations, we make extensive use of \eqref{eq:leq2}.  We first note that
\[
a = bq \Implies a = bb^*a \Implies aa^* = bb^*aa^* \Implies aa^* \leq bb^* .
\]
We also have
\[
a = pb \Implies a = pa \Implies aa^* = paa^* \Implies aa^* \leq p \Implies aa^* = aa^*p.
\]
We then calculate
\[
a = aa^*a = aa^*(pb) = (aa^*p)b = aa^*b.
\]
Writing $r=aa^*\in P$, we have already seen that $r = aa^* \leq bb^* = \bd(b)$, so that in fact $a = aa^*b = rb = {}_r\corest b$, and hence $a\leq b$.
\end{rem}

\begin{rem}
Recall that any regular semigroup $S$ has a so-called \emph{natural partial order}.  This order, which for clarity we will denote by $\preceq$, has \emph{many} different formulations; see for example \cite{Hartwig1980,Nambooripad1980,Hickey1983}, and especially \cite{Mitsch1986} for a discussion of the variations, and an extension to arbitrary semigroups.  The most straightforward definition of the natural partial order on the regular semigroup $S$ is:
\[
a\preceq b \Iff a\in Eb\cap bE \qquad\text{where $E=E(S)$.}
\]
It is natural to ask how the orders $\leq$ and $\preceq$ are related when $S$ is a regular $*$-semigroup, where~$\leq$ is the order from \eqref{eq:aleqb2}.  Since $P\sub E$ (and since ${}_p\corest b\in Pb$ and $b\rest_q\in bP$ for $p\leq\bd(b)$ and $q\leq\br(b)$), the order~$\leq$ is of course contained in $\preceq$, meaning that $a\leq b\implies a\preceq b$.  But the converse does not hold in general.  For example, if $S$ is a regular $*$-\emph{monoid} with identity~$1$, and if $e\in E\setminus P$ is any non-projection idempotent, then $e\preceq1$ but $e\not\leq1$.
\end{rem}

The construction of $\G(S)$ from $S$ can be thought of as an object map from the category~$\RSS$ of regular $*$-semigroups (with $*$-morphisms) to the category $\OG$ of ordered groupoids (with ordered groupoid functors).  Since the underlying sets of $S$ and $\G(S)$ are the same, any $*$-morphism $\phi:S\to S'$ can be thought of as a map $\G(\phi)=\phi:\G(S)\to\G(S')$.  The main result of this section shows that $\G$ (interpreted in this way) is a functor:

\begin{prop}\label{prop:calGfunctor}
$\G$ is a functor $\RSS\to\OG$.
\end{prop}

\pf
%By Proposition \ref{prop:GS}, $\G$ maps $v\RSS$ to $v\OG$.  Now consider regular $*$-semigroups $S$ and~$S'$, and write $\G=\G(S)$ and $\G'=\G(S')$.  We need to show that 
It remains to show that any $*$-morphism $\phi:S\to S'$ in $\RSS$ is an ordered groupoid morphism $\G\to\G'$, where $\G=\G(S)$ and $\G'=\G(S')$, i.e.~that
\ben
\item \label{RSSOG1} $(a\circ b)\phi = (a\phi)\circ (b\phi)$ for $a,b\in\G$ with $\br(a)=\bd(b)$, and
\item \label{RSSOG2} $({}_p\corest a)\phi = {}_{p\phi}\corest a\phi$ for $a\in\G$ and $p\leq\bd(a)$.%
\footnote{Strictly speaking, we must also check that $\G(\phi\circ\phi')=\G(\phi)\circ\G(\phi')$ for composable morphisms $\phi,\phi'$ in $\RSS$, and that $\G(\id_S)=\id_{\G(S)}$ for all $S$.  However, these are trivial.  This is the case for all other functors constructed in the paper, so we will pass over this step in all such proofs.}
\een
Before we do this, we first note that for any $a\in\G$ we have
\[
\bd(a\phi) = (a\phi)(a\phi)^* = (aa^*)\phi = \bd(a)\phi \ANDSIM \br(a\phi) = \br(a)\phi.
\]
We make repeated use of this in the following, as well as \eqref{eq:leq2} and Definition \ref{defn:GS}.

\pfitem{\ref{RSSOG1}}  If $\br(a)=\bd(b)$, then $\br(a\phi) = \br(a)\phi = \bd(b)\phi = \bd(b\phi)$, and so 
\[
(a\circ b)\phi = (ab)\phi = (a\phi)(b\phi) = (a\phi) \circ (b\phi).
\]
\pfitem{\ref{RSSOG2}}  If $p\leq\bd(a)$, then $p=p\cdot\bd(a)$ and so $p\phi = (p\phi)(\bd(a)\phi) = (p\phi)\cdot\bd(a\phi)$, so that $p\phi\leq\bd(a\phi)$, and then
\[
({}_p\corest a)\phi = (pa)\phi = (p\phi)(a\phi) = {}_{p\phi}\corest a\phi.  \qedhere
\]
\epf

\begin{rem}
The functor $\G$ is not surjective (on objects).  To see this, we first note that any poset can be thought of as an ordered groupoid consisting entirely of (pairwise-noncomposable) identities.  Consider the case in which $\G=\{p,q\}$ is a two-element anti-chain.  Aiming for a contradiction, suppose $\G=\G(S)$ for some regular $*$-semigroup $S$.  We then have $P(S)=\{p,q\}$, and we claim that $p\F q$.  Indeed, if this was not the case then either $p=p\th_q\leq q$ or $q=q\th_p\leq p$, neither of which hold.  Since~$p\F q$, it follows that $e=pq\in S(=\G)$ satisfies $\bd(e)=p$ and $\br(e)=q$ (cf.~\eqref{eq:ee*e*e}), i.e.~$e\in\G(p,q)$, contradicting $\G(p,q)=\es$.

%since $p\not\leq q$, we cannot have $p\th_q=p$ (cf.~\eqref{eq:leq}) so in fact $p\th_q=q$, and similarly $q\th_p=p$.  
%
%, as an ordered groupoid need not have a (so-called) projection algebra for its object set.  

We will see in Section \ref{subsect:Seg} that $\G$ is not injective (on objects) either.  That is, it is possible to have $\G(S)=\G(S')$ for distinct regular $*$-semigroups $S$ and $S'$.  In a sense, the rest of the paper is devoted to the task of finding a suitable enrichment of the groupoids $\G(S)$ to obtain a faithful representation of regular $*$-semigroups by ordered groupoids.

On the other hand, the ESN Theorem tells us that the restriction of the functor $\G$ to the category of inverse semigroups $\IS\sub\RSS$ \emph{is} injective, with image $\IG$, the category of inductive groupoids.
\end{rem}

%\subsection{\JErev{Going backwards}}

\begin{rem}\label{rem:GS}
The (partial) composition of the groupoid $\G=\G(S)$ is a restriction of the (total) product of the regular $*$-semigroup~$S$.  Nevertheless, the order on $\G$ can be used to reduce the calculation of an arbitrary product $ab$ in $S$ to an associated composition in $\G$:
\begin{equation}\label{eq:ab}
ab = a' \circ e \circ b',
\end{equation}
where $a'\leq a$ and $b'\leq b$, and where $e\in\G$ is a special morphism $\br(a')\to\bd(b')$.  Specifically, let us write
\[
p=\br(a)=a^*a \AND q=\bd(b)=bb^*,
\]
noting that $a=ap$ and $b=qb$.  We also define the additional projections
\[
p' = q\th_p = pqp \leq p \AND q' = p\th_q = qpq \leq q.
\]
Then since $pq$ is an idempotent (cf.~Lemma \ref{lem:PS1}\ref{PS12}), we have $p'q'=(pq)^3=pq$, and so
\[
ab = ap \cdot qb = a\cdot p'q'\cdot b = ap' \cdot p'q' \cdot q'b.
\]
Then with $a'=a\rest_{p'}=ap'$, $e=p'q'(=pq)$ and $b'={}_q\corest b = q'b$, one can check that $p'\F q'$, from which it follows (as in \eqref{eq:ee*e*e}) that
\[
\bd(e)=p' = \br(a') \AND \br(e)=q'=\bd(b').
\]
%For example, $\bd(e) = ee^* = p'q'(p'q')^* = p'q'p' = (pqp)(qpq)(pqp) = pqp = p'$.
%\[
%\br(a') = (ap')^*ap' = (apqp)^*apqp = pqpa^*apqp = pqp\cdot p\cdot pqp = pqpqp = pqp = p'.
%\]
Thus, continuing from above, it follows that indeed
\[
ab = ap' \cdot p'q' \cdot q'b = a' \cdot e \cdot b' = a' \circ e \circ b'.
\]
%In what follows, we will think of $a'=ap'$ and $b'=q'b$ as `restrictions' $a'=a\rest_{p'}$ and $b'={}_{q'}\corest b$ in the groupoid $\G$.
\end{rem}

\begin{rem}\label{rem:e}
There is a rather subtle point regarding Remark \ref{rem:GS} that is far from obvious at this stage, but ought to be mentioned now.  It may seem as if we are implying that the structure of the regular $*$-semigroup $S$ is completely determined by that of the ordered groupoid $\G=\G(S)$.  That is, if we are given a complete description of the groupoid $\G$, including its composition, inversion and ordering, 
%and the (so-called) projection algebra structure of its object set $P=v\G$, 
then it might seem that we ought to be able to construct the entire multiplication table of $S$.  
However, this is far from the truth.  Indeed, in Section \ref{subsect:Seg}, we will see that it is possible for distinct (even non-isomorphic) regular $*$-semigroups $S$ and~$S'$ to have exactly the same groupoids $\G(S)=\G(S')$; see Examples \ref{eg:AG} and \ref{eg:Rees}.  %(We postpone the definition of these semigroups, as we have already gone on a somewhat lengthy tangent, and also since we will only be able to fully appreciate their properties once we have developed more theory.)

If the reader is worried that this contradicts the previous remark,  the source of the subtlety is in the precise identity of the element~$e$ from \eqref{eq:ab} that allowed us to reduce a product $ab$ in~$S$ to a composition in $\G$:
\[
ab = a\rest_{p'}\circ e \circ {}_{q'}\corest b.
\]
(On the other hand, the projections $p',q'$ can be found directly using the so-called \emph{projection algebra} structure of $P=v\G$.)  We were able to locate this element $e$ in Remark \ref{rem:GS} because we began with full knowledge of the semigroup $S$; we simply took $e=pq=p'q'$.  However, if we begin with an ordered groupoid $\G$, it is not so obvious what this element $e$ should be, even if we know that $\G=\G(S)$ for some (unknown) regular $*$-semigroup $S$.  Certainly $e$ should be a morphism $p'\to q'$, but when $p'\not=q'$, it is not so easy to distinguish any such morphism.  Getting around this problem is one of the main sources of difficulty encountered in the current work.  Our eventual solution utilises what we will call the `chain groupoid' $\C=\C(P)$ associated to an (abstract) projection algebra $P$, and a certain `evaluation map' $\ve:\C\to\G$.  We will find our element $e$ in the image of this map.
\end{rem}

\begin{rem}\label{rem:trace}
Before moving on, it is worth commenting on another related notion, namely the so-called \emph{trace} of a semigroup.  This goes back to a classical result of Miller and Clifford.  Among other things, \cite[Theorem 3]{MC1956} says that for elements $a,b$ of a semigroup $S$, we have
\[
ab \in R_a \cap L_b \Iff L_a \cap R_b \text{ contains an idempotent}.
\]
Such products $ab \in R_a \cap L_b$ are often called \emph{trace products} in the literature.  We can then define a partial binary operation~$\bc$ on $S$ by
\[
a\bc b = \begin{cases}
ab &\text{if $L_a \cap R_b$ contains an idempotent}\\
\text{undefined} &\text{otherwise,}
\end{cases}
\]
and the resulting partial algebra $(S,{\bc})$ is often called the \emph{trace of $S$}.  Following Definition \ref{defn:GS}, it is not hard to see that in a regular $*$-semigroup $S$, a composition $a\circ b$ exists precisely when $L_a\cap R_b$ contains a projection, which must of course be $a^*a=bb^*$, i.e.~$\br(a)=\bd(b)$.  Since projections are idempotents, it follows that $a\circ b$ being defined forces $a\bc b$ to be defined, and then of course $a\circ b=a\bc b=ab$, meaning that ${\circ}\sub{\bc}$, i.e.~that $\bc$ is an extension of $\circ$.  We have already noted that $S$ is inverse if and only if $P(S)=E(S)$, and it now quickly follows that this is also equivalent to having ${\circ}={\bc}$.  In this inverse case, the trace $(S,{\bc})$ is then precisely the groupoid $\G(S)$ from Definition \ref{defn:GS}, but this is not true of non-inverse regular $*$-semigroups.
\end{rem}

\subsection{Adjacency semigroups and Rees 0-matrix semigroups}\label{subsect:Seg}

The rest of Section \ref{sect:RSS} is devoted to a discussion of some natural examples of regular $*$-semigroups.  We will return to these over the course of the theoretical development below, as they will provide concrete instances of some of the subtleties that arise; see for example Remarks \ref{rem:Rees1}, \ref{rem:Rees2} and~\ref{rem:G3'}.  In the current section we treat adjacency semigroups and Rees 0-matrix semigroups.  We look at diagram semigroups in Section \ref{subsect:D}.

%\begin{eg}\label{eg:PxP}
%%\JE{Not sure if we want to include this, as it was maybe just to illustrate a point in the previous version, maybe comparing to a Brandt semigroup?}
%Let $P$ be an arbitrary set, and let $S=P\times P=\set{(p,q)}{p,q\in P}$ be the cartesian product of two copies of $P$.  We define binary and unary operations on $S$ by:
%\[
%(p,q)(r,s)=(p,s) \AND (p,q)^*=(q,p).
%\]
%It is then routine to check that $S$ is a regular $*$-semigroup.  This semigroup is known as the \emph{square band} over $P$.  (More generally, a \emph{rectangular band} has the form $P\times Q$, for possibly different sets~$P$ and~$Q$, with the same product as above; this is still a semigroup, but need not have an involution.)  Identifying $p\equiv(p,p)$ for each $p$, we see that $P(S) \equiv P$, and moreover each~$\th_p$ operation is the constant map with image $p$.  It follows from this that the ordering $\leq$ on~$P$ from \eqref{eq:leq} is the equality relation.
%\end{eg}

\begin{eg}\label{eg:AG}
Let $\Ga=(P,E)$ be a digraph, so that $P$ is a vertex set, and $E\sub P\times P$ an edge set.  The \emph{adjacency semigroup} of $\Ga$, denoted $A(\Ga)$, has underlying set $(P\times P)\cup\{0\}$, and product defined by
\begin{equation}\label{eq:AG.}
a0=0a=0 \quad\text{for all $a\in A(\Ga)$} \AND (p,q)(r,s) = \begin{cases}
(p,s) &\text{if $(q,r)\in E$}\\
0 &\text{otherwise.}
\end{cases}
\end{equation}
We also have a unary operation ${}^*$ defined by
\begin{equation}\label{eq:AG*}
0^* = 0 \AND (p,q)^*=(q,p).
\end{equation}
Clearly $A(\Ga)$ satisfies the identity $(a^*)^* = a$.  Adjacency semigroups were studied by Jackson and Volkov in \cite{JV2010}, who observed that:
\bit
\item the identity $(ab)^*=b^*a^*$ is satisfied if and only if the adjacency relation is symmetric, i.e.~$(p,q)\in E \iff (q,p)\in E$ for all $p,q\in P$, and
\item the identity $a=aa^*a$ is satisfied if and only if the adjacency relation is reflexive, i.e.~$(p,p)\in E$ for all $p\in P$.
\eit
Thus, $S=A(\Ga)$ is a regular $*$-semigroup (under the above operations) precisely when $E$ is reflexive and symmetric, i.e.~when $\Ga$ is an undirected graph with loops at every vertex.  When these conditions hold, and identifying $p\equiv(p,p)$ for $p\in P$, we have $P(S) \equiv P\cup\{0\}$, which for convenience we will denote by $P_0$.

The $\th_0$ operation on $P(S)=P_0$ is of course the constant map with image $0$, while for $p\in P$ and $q\in P_0$ we have 
\begin{equation}\label{eq:AGth}
q\th_p = \begin{cases}
p &\text{if $q\in P$ and $(p,q)\in E$}\\
0 &\text{otherwise.}
\end{cases}
\end{equation}
It follows from this that the poset $(P_0,{\leq})$, as in \eqref{eq:leq}, is the following \emph{fan semilattice},%
\footnote{Note that the ordered set $(P_0,{\leq})$ being a semilattice does not imply that $(P_0,{\cdot})$ is a semilattice in the semigroup-theoretic sense (i.e.~a commutative band).  Consulting Lemma \ref{lem:inv} (or using \eqref{eq:AG.}), the latter is the case precisely when $S$ is inverse, which occurs when $\Ga$ is a disjoint union of loops.}
shown in the case $P=\{p_1,\ldots,p_5\}$:
\begin{equation}\label{eq:AGleq}
\begin{tikzpicture}
\foreach \x in {1,...,5} {\node (p\x) at (\x,2) {$p_\x$};}
\node (0) at (3,0) {$0$};
\foreach \x in {1,...,5} {\draw (p\x)--(0);}
\end{tikzpicture}
\end{equation}
Indeed, we clearly have $0\leq p$ for all $p\in P_0$.  If $p\leq q$ for $p,q\in P$, then $p=p\th_q$; since $\im(\th_q)=\{0,q\}$ (cf.~\eqref{eq:AGth}), it follows that $p=q$.
On the other hand, it follows from \eqref{eq:AGth} that the $\F$ relation on $P_0$ is given by ${\F} = E\cup\{(0,0)\}$, i.e.~that
\begin{equation}\label{eq:AGF}
p\F q \IFF p=q=0 \text{ \ \ or \ \ $(p,q)$ is an edge of $\Ga$.}
\end{equation}
Consequently, $\F$ is transitive if and only if $\Ga$ is a disjoint union of complete graphs.  Combining ${\F}=E\cup\{(0,0)\}$ with Lemma \ref{lem:PS1}\ref{PS12}, it also follows that $E(S)=P_0^2=E\cup\{0\}$.

We now consider the ordered groupoid $\G=\G(S)$ for the adjacency graph $S=A(\Ga)$.  First, we have $v\G=P(S) = P_0$.  It is easy to check that
\[
(p,q) \R (r,s) \iff p=r \AND (p,q) \L (r,s) \iff q=s,
\]
so that the non-empty morphism sets of $\G$ are
\[
\G(0,0) = \{0\} \AND \G(p,q) = \{(p,q)\} \quad\text{for $p,q\in P$.}
\]
Thus, the only compositions in $\G$ are $(p,q)\circ(q,r) = (p,r)$ for all $p,q,r\in P$, and $0\circ0=0$.  Moreover, due to the nature of the order $(P_0,{\leq})$, the only non-trivial restrictions in $\G$ are ${}_0\corest(p,q) = (p,q)\rest_0 = 0$ for $p,q\in P$.  This all shows that the groupoid $\G=\G(S)$ depends only on the vertex set $P$ of the underlying graph $\Ga=(P,E)$, and not on its edge set $E$.  
%
%Consequently, two different graphs $\Ga=(P,E)$ and $\Ga'=(P,E')$ give rise to different adjacency semigroups $S=A(\Ga)$ and $S'=A(\Ga')$, even though $\G(S)=\G(S')$ are identical as ordered groupoids.  As an extreme example, we can take $\Ga$ to be the discrete graph (where $E$ consists only of loops at each vertex), and $\Ga'$ to be the complete graph.  Here $S$ is a $P\times P$ \emph{combinatorial Brandt semigroup} and is hence inverse, while $S'$ is a $P\times P$ square band with a zero adjoined.  Since $P_0$ is a semilattice (cf.~\eqref{eq:AGleq}), the groupoid $\G(S)=\G(S')$ is inductive, even though $S'$ is not inverse.

Consequently, two different graphs $\Ga=(P,E)$ and $\Ga'=(P,E')$ give rise to distinct adjacency semigroups $S=A(\Ga)$ and $S'=A(\Ga')$, even though $\G(S)=\G(S')$ are identical as ordered groupoids.  Moreover, the semigroups $S$ and $S'$ need not even be isomorphic.  As an extreme example, we can take $\Ga$ to be the discrete graph (where $E$ consists only of loops at each vertex), and $\Ga'$ to be the complete graph.  Here $S$ is a $P\times P$ \emph{combinatorial Brandt semigroup} and is hence inverse, while $S'$ is a $P\times P$ square band with a zero adjoined.  Since $(P_0,{\leq})$ is a semilattice (cf.~\eqref{eq:AGleq}), the groupoid $\G(S)=\G(S')$ is inductive, even though $S'$ is not inverse for $|P|\geq2$.
\end{eg}

Regular Rees 0-matrix semigroups are one of the most important families of semigroups, being (up to isomorphism) precisely the completely 0-simple semigroups, as per the celebrated \emph{Rees Theorem}; see \cite{Rees1940}, and also \cite[Chapter~3]{CPbook} or \cite[Section 3.2]{Howie1995}.  It is therefore important to know which of these are regular $*$-semigroups:

\begin{eg}\label{eg:Rees}
Let $P$ and $Q$ be arbitrary non-empty sets, $G$ a group, $0$ a symbol not belonging to $G$, and $M=(m_{qp})$ a $Q\times P$ \emph{regular sandwich matrix} over $G\cup\{0\}$.  Regularity here means that~$M$ has at least one non-zero entry in each row and in each column.  The \emph{Rees 0-matrix semigroup} $S = \M^0(P,Q,G,M)$ has underlying set $(P\times G\times Q)\cup\{0\}$ and product
\begin{equation}\label{eq:M0.}
a0=0a=0 \quad\text{for all $a\in S$} \AND (p,g,q)(r,h,s) = \begin{cases}
(p,gm_{qr}h,s) &\text{if $m_{qr}\not=0$}\\
0 &\text{otherwise.}
\end{cases}
\end{equation}
The non-zero $\R$- and $\L$-classes of $S$ are respectively the sets
\[
R_p = \{p\}\times G\times Q \AND L_q = P\times G\times\{q\} \qquad\text{for $p\in P$ and $q\in Q$,}
\]
the non-zero $\H$-classes are the sets
\[
H_{pq} = R_p\cap L_q = \{p\}\times G\times\{q\} \qquad\text{for $p\in P$ and $q\in Q$,}
\]
and the non-zero idempotents have the form $(p,m_{qp}^{-1},q)$ for each $p\in P$ and $q\in Q$ with $m_{pq}\not=0$.  In particular, $H_{pq}$ contains an idempotent if and only if $m_{qp}\not=0$.  One can \emph{scale} the matrix $M$ by left-multiplying each column by any non-zero group element to obtain an isomorphic Rees 0-matrix semigroup.

Suppose now that $S$ is a regular $*$-semigroup, with respect to some involution ${}^*:S\to S$.  Since~${}^*$~maps $\R$-classes to $\L$-classes, and vice versa, there is a bijection $P\to Q:p\mt p'$ such that $R_p^* = L_{p'}$ for all $p\in P$.  Thus, re-naming the elements of $Q$, we can assume that $Q=P$ and $p'=p$ for all $p\in P$.  Note then that $H_{pq}^*=H_{qp}$ for all $p,q\in P$.  In particular, the only $\H$-classes fixed by ${}^*$ are $\{0\}$ and $H_{pp}$ ($p\in P$).  Since each $\R$-class contains a projection, it follows that the projection in the non-zero $\R$-class $R_p$ is the idempotent of $H_{pp}$.  In particular, we have $m_{pp}\not=0$ for all $p\in P$.  We can scale the sandwich matrix $M$ by multiplying each column $p$ by $m_{pp}^{-1}$, in order to assume that in fact $m_{pp}=1$ for all $p$.  In this way, we can identify 
\[
p \equiv (p,1,p) \qquad\text{for all $p\in P$,}
\]
and again we have $P(S)\equiv P_0 = P\cup\{0\}$.
%so that the projections and idempotents of $S$ are precisely
%\[
%P(S) \equiv P\cup\{0\} \AND E(S) = \set{(p,m_{qp}^{-1},q)}{p,q\in P,\ m_{qp}\not=0} \cup \{0\}.
%\]
Since ${}^*$ maps group $\H$-classes to group $\H$-classes, we must have the condition
\[
m_{qp}=0 \iff m_{pq}=0 \qquad\text{for all $p,q\in P$.}
\]
Now consider a non-zero idempotent $e$, say with $e\in H_{pq}$.  As noted above, we have $e=(p,m_{qp}^{-1},q)$.  But since also $e\in H_{pq}=R_p\cap L_q$, we have $e=pq\equiv(p,1,p)(q,1,q) = (p,m_{pq},q)$, and this shows that $m_{pq}=m_{qp}^{-1}$.  If we write $0^{-1}=0$ as well, then it follows that the sandwich matrix $M$ has the skew-symmetric condition
\begin{align}
\label{eq:ss} m_{pq} &= m_{qp}^{-1} &&\hspace{-2.5cm}\text{for all $p,q\in P$,}
\intertext{and we recall that necessarily}
\label{eq:sss} m_{pp} &= 1 &&\hspace{-2.5cm}\text{for all $p\in P$.}
\end{align}
Now consider some non-zero element $a=(p,g,q)\in S$.  Since $aa^*$ is the unique projection in the $\R$-class $R_a=R_p$, and since $a^*\in H_{pq}^*=H_{qp}$, we have $a^*=(q,h,p)$ for some $h\in G$, and then
\[
(p,1,p) \equiv p = aa^* = (p,g,q)(q,h,p) = (p,gm_{qq}h,p) = (p,gh,p),
\]
so that $h=g^{-1}$.  In other words, we have
\begin{equation}\label{eq:M0*}
(p,g,q)^* = (q,g^{-1},p) \quad\text{for all $p,q\in P$ and $g\in G$,} \qquad\text{and of course}\qquad 0^*=0.
\end{equation}
To summarise, up to isomorphism, the Rees 0-matrix regular $*$-semigroups are precisely those of the form $\M^0(P,P,G,M)$ whose sandwich matrix $M=(m_{qp})$ satisfies~\eqref{eq:ss} and \eqref{eq:sss}, and then the product and involution are given in \eqref{eq:M0.} and \eqref{eq:M0*}.  

Again $\th_0$ is a constant map, while for $p\in P$ and $q\in P_0$ we have
\begin{equation}\label{eq:M0th}
q\th_p = \begin{cases}
p &\text{if $q\in P$ and $m_{pq}\not=0$}\\
0 &\text{otherwise.}
\end{cases}
\end{equation}
The poset $(P_0,{\leq})$ is again the fan semilattice, as in \eqref{eq:AGleq}.  As in \eqref{eq:AGF}, the $\F$ relation on $P_0$ is given by
\begin{equation}\label{eq:M0F}
p\F q \IFF p=q=0 \text{ \ \ or \ \ $m_{pq}\not=0$.}
\end{equation}

The groupoid $\G=\G(S)$ can also be easily described.  We again have $v\G=P(S)=P_0$, but this time the morphism sets are
\[
\G(0,0)=\{0\} \AND \G(p,q) = H_{pq} = \{p\}\times G\times\{q\} \quad\text{for $p,q\in P$.}
\]
Composition is given by $(p,g,q)\circ(q,h,r) = (p,gh,r)$, remembering that $m_{qq}=1$, and inversion by $(p,g,q)^{-1}=(q,g^{-1},r)$.  Again, the only non-trivial restrictions are ${}_0\corest(p,g,q) = (p,g,q)\rest_0 = 0$.  It follows that the groupoid $\G(S)$ is independent of the sandwich matrix $M$, and depends only on the indexing set $P$ and the group $G$.
\end{eg}

\begin{rem}
If all entries of the sandwich matrix $M=(m_{qp})$ are non-zero, then the \emph{Rees matrix semigroup} $\M(P,P,G,M) = P\times G\times P$ is a subsemigroup of $\M^0(P,P,G,M)$.  It is a regular $*$-semigroup under
\[
(p,g,q)(r,h,s) = (p,gm_{qr}h,s) \AND (p,g,q)^* = (q,g^{-1},p).
\]
Square bands are a special case, where $G=\{1\}$ is trivial.  Involutions on Rees matrix semigroups were described in \cite[Section 3]{Petrich1985}, but Rees 0-matrix semigroups were not discussed.
\end{rem}

\begin{rem}
Adjacency semigroups are special cases of Rees 0-matrix semigroups.  Given a graph $\Ga=(P,E)$, we take $G=\{1\}$ to be a trivial group, and we have $A(\Ga) \cong \M^0(P,P,G,M)$, where $M$ is the adjacency matrix of $\Ga$, i.e.~$m_{pq}=1 \iff (p,q)\in E$.
\end{rem}

%\begin{eg}
%Finally, we briefly discuss the \emph{bicyclic monoid} $B$, as defined by the presentation
%\[
%B = \pres{a,b}{ba=1}.
%\]
%The elements of $B$ can be identified with words of the form $a^ib^j$.  The product of such words $a^ib^j$ and $a^kb^l$ is obtained by first forming the word $a^ib^ja^kb^l$, and then using the relation $ba=1$ to simplify $b^ja^k$ to a power of $b$ or $a$.  Identifying $a^ib^j$ with the pair $(i,j)\in\N\times\N$, where $\N=\{0,1,2,\ldots\}$, it is easy to see that $B$ is inverse, with $(i,j)^{-1}=(j,i)$.  Moreover, we have
%\[
%E(B) = P(B) = \set{(i,i)\equiv i}{i\in\N} \equiv \N.
%\]
%The $\th$ operations are given by $j\th_i = \max(i,j)$, and the $\leq$ ordering on $P(B)\equiv\N$ is the reverse of the usual ordering on natural numbers.  In the groupoid $\G=\G(B)$ we have
%\[
%\bd(i,j)=i \COMMA \br(i,j)=j \AND (i,j)\circ(j,k)=(i,k).
%\]
%In particular, $\G(B)$ is precisely the same as $\G(S)$, where $S=\M(\N,\N,G,M)$, where $G$ is the trivial group, and $M$ the all-$1$ matrix.  This is only true at the groupoid level, however, as the orders on $\G(S)$ and $\G(B)$ are different.
%
%\end{eg}

\subsection{Diagram monoids}\label{subsect:D}

%One of the most important and well-studied collections of non-inverse regular $*$-semigroups are the so-called \emph{diagram monoids}, including the Brauer monoids \cite{Brauer1937}, Temperley-Lieb monoids \cite{TL1971}, partition monoids \cite{Jones1994_2,Martin1994}, Motzkin monoids~\cite{BH2014} and several others.  

The theory of regular $*$-semigroups has seen quite a resurgence in recent years, partly due to the importance of so-called \emph{diagram monoids}.  Here we recall the definition of these monoids, and use them to illustrate some of the new ideas and constructions of the current paper, but we also relate this to previous studies.  We will return to these monoids at various points during the rest of the paper (see for example Remarks \ref{rem:notFfunctor} and \ref{rem:G3'}), and in the sequels \cite{Paper2,Paper3}.  Here we focus exclusively on the \emph{partition monoids}~\cite{Martin1994,Jones1994_2}, although similar comments could be made for other diagram monoids, such as (partial) Brauer, Temperley-Lieb and Motzkin monoids \cite{Brauer1937,TL1971,DEG2017,EG2017,BH2014,MM2014}.  

Let $X$ be a set, and let $X'=\set{x'}{x\in X}$ be a disjoint copy of $X$.  The \emph{partition monoid over $X$}, denoted $\PP_X$, is defined as follows.  The elements of $\PP_X$ are the set partitions of $X\cup X'$, and the operation will be defined shortly.  A partition from $\PP_X$ will be identified with any graph on vertex set $X\cup X'$ whose connected components are the blocks of the partition.  The elements of $X$ and $X'$ are called \emph{upper} and \emph{lower} vertices, respectively, and are generally displayed in two parallel rows, as in the various figures below.

To define the product on $\PP_X$, consider two partitions $\al,\be\in\PP_X$.  Let $X''=\set{x''}{x\in X}$ be another disjoint copy of $X$, and define three new graphs:
\bit
\item $\al^\vee$, the graph on vertex set $X\cup X''$ obtained from $\al$ by changing each lower vertex $x'$ to $x''$, 
\item $\be^\wedge$, the graph on vertex set $X''\cup X'$ obtained from $\be$ by changing each upper vertex $x$ to $x''$,
\item $\Pi(\al,\be)$, the graph on vertex set $X\cup X''\cup X'$, whose edge set is the union of the edge sets of~$\al^\vee$ and $\be^\wedge$.
\eit
The graph $\Pi(\al,\be)$ is called the \emph{product graph} of $\al$ and $\be$; it is generally drawn with $X''$ as the middle row of vertices.  The product $\al\be$ is defined to be the partition of $X\cup X'$ with the property that $u,v\in X\cup X'$ belong to the same block of $\al\be$ if and only if $u,v$ are connected by a path in~$\Pi(\al,\be)$.  

If $X=\{1,\ldots,n\}$ for a non-negative integer $n$, we write $\PP_n=\PP_X$.  

\begin{figure}[ht]
\begin{center}
\begin{tikzpicture}[scale=.53]
\begin{scope}[shift={(0,0)}]	
\uvs{1,...,6}
\lvs{1,...,6}
\uarc12
\uarc23
\darc23
\darc45
\darc56
\stline11
\stline44
\draw(0.6,1)node[left]{$\al=$};
\draw[->](7.5,-1)--(9.5,-1);
\end{scope}
\begin{scope}[shift={(0,-4)}]	
\uvs{1,...,6}
\lvs{1,...,6}
\uarc23
\uarc45
\uarc56
\darc12
\darc23
\darc46
\stline14
\stline43
\draw(0.6,1)node[left]{$\be=$};
\end{scope}
\begin{scope}[shift={(10,-1)}]	
\uvs{1,...,6}
\lvs{1,...,6}
\uarc12
\uarc23
\darc23
\darc45
\darc56
\stline11
\stline44
\draw[->](7.5,0)--(9.5,0);
\end{scope}
\begin{scope}[shift={(10,-3)}]	
\uvs{1,...,6}
\lvs{1,...,6}
\uarc23
\uarc45
\uarc56
\darc12
\darc23
\darc46
\stline14
\stline43
\end{scope}
\begin{scope}[shift={(20,-2)}]	
\uvs{1,...,6}
\lvs{1,...,6}
\uarc12
\uarc23
\darc12
\darc23
\darc46
\stline34
\stline43
\draw(6.4,1)node[right]{$=\al\be$};
\end{scope}
\end{tikzpicture}
\caption{Left to right: partitions $\al,\be\in\PP_6$, the product graph $\Pi(\al,\be)$, and the product $\al\be\in\PP_6$.  For more information, see Example \ref{eg:PX1}.}
\label{fig:PX1}
\end{center}
\end{figure}
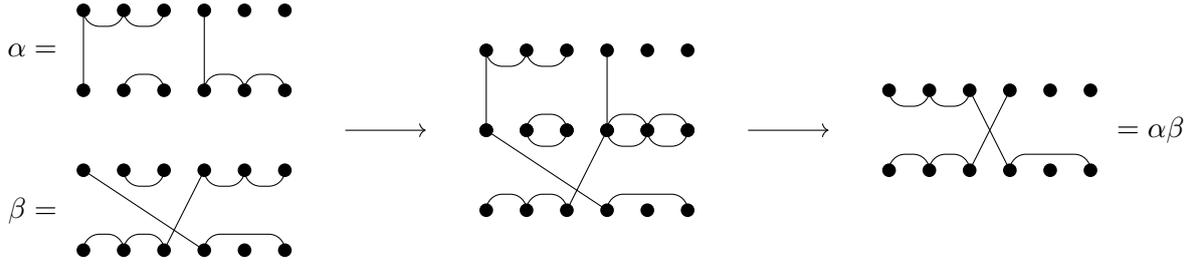

\begin{eg}\label{eg:PX1}
To illustrate the product above, consider the partitions $\al,\be\in\PP_6$ defined by
\begin{align*}
\al &= \big\{ \{1,2,3,1'\}, \{4,4',5',6'\}, \{5\},\{6\},\{2',3'\}\big\} , \\
\be &= \big\{\{1,4',6'\},\{2,3\},\{4,5,6,1',2',3'\},\{5'\}\big\}.
\intertext{Figure \ref{fig:PX1} illustrates (graphs representing) $\al$ and $\be$, their product}
\al\be &= \big\{\{1,2,3,4',6'\},\{4,1',2',3'\},\{5\},\{6\},\{5'\}\big\},
\end{align*}
as well as the product graph $\Pi(\al,\be)$.  Here and elsewhere, vertices are assumed to increase from left to right, $1,\ldots,n$, and similarly for (double) dashed vertices.
\end{eg}

\begin{figure}[ht]
\begin{center}
\begin{tikzpicture}[scale=.6]
\begin{scope}[shift={(0,0)}]	
\uvs{1,...,20}
\lvs{1,...,20}
\uarcs{1/2,2/3,4/5,7/8,8/9,9/10,11/12,12/13,13/14,14/15,15/16,16/17,17/18,18/19,19/20}
\darcs{1/2,3/4,7/8,8/9,12/13,14/15,16/17,18/19,19/20}
\stlines{3/5,6/6,10/11}
\draw(0.6,1)node[left]{$\al=$};
\draw(21,0)node[right]{$\Pi(\al,\be)$};
\draw[|-|] (21,2)--(21,-2);
\end{scope}
\begin{scope}[shift={(0,-2)}]	
\uvs{1,...,20}
\lvs{1,...,20}
\uarcs{2/3,5/6,6/7,9/10,11/12,13/14,17/18}
\darcs{2/3,5/6,6/7,7/8,11/12,12/13,13/14,14/15,16/17,17/18,18/19,19/20}
\stlines{1/1,4/4,8/8,15/15}
\draw(0.6,1)node[left]{$\be=$};
\end{scope}
\begin{scope}[shift={(0,-6)}]	
\uvs{1,...,20}
\lvs{1,...,20}
\uarcs{1/2,2/3,4/5,7/8,8/9,9/10,11/12,12/13,13/14,14/15,15/16,16/17,17/18,18/19,19/20}
\uarcx36{.7}
\darcs{2/3,5/6,6/7,7/8,11/12,12/13,13/14,14/15,16/17,17/18,18/19,19/20}
\darcx14{.7}
\stlines{6/6,10/11}
\draw(0.6,1)node[left]{$\al\be=$};
\end{scope}
\end{tikzpicture}
\caption{Top: partitions $\al,\be\in\PP_{20}$, already connected to create the product graph $\Pi(\al,\be)$.  Bottom: the product $\al\be\in\PP_{20}$.  For more information, see Example \ref{eg:PX2}.}
\label{fig:PX2}
\end{center}
\end{figure}
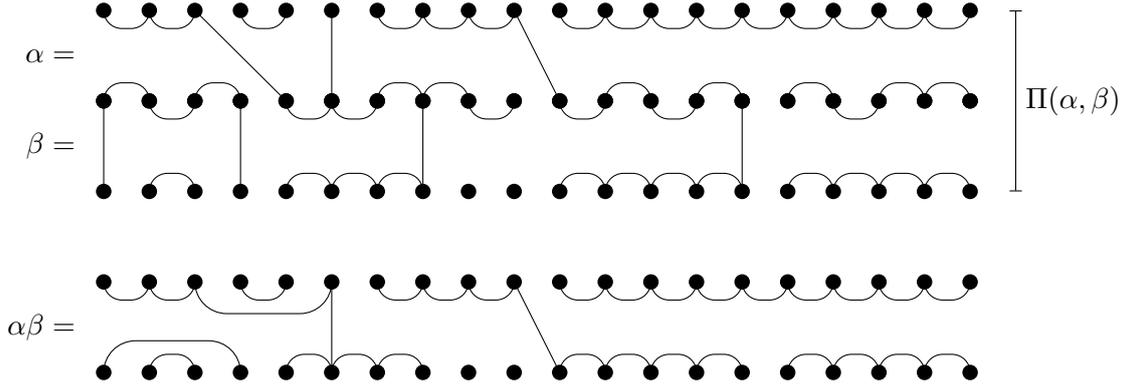

A reader familiar with partition monoids might protest that we have `cheated' in Example~\ref{eg:PX1}, and that our choice of $\al$ and $\be$ was too `easy' to fully illustrate the complexities of calculating products of partitions.  This is indeed the case, as the bottom half of $\al$ `matches' the top half of $\be$ in a way that can hopefully be understood by examining Figure \ref{fig:PX1}, but which will be made precise shortly.  Before this, we briefly consider a more `difficult' product:

\begin{eg}\label{eg:PX2}
Figure \ref{fig:PX2} gives another product, this time with $\al,\be\in\PP_{20}$.  One can immediately see that there is no such `matching' between the bottom of $\al$ and the top of $\be$.  Rather, to calculate the product $\al\be$, one needs to follow paths in the product graph $\Pi(\al,\be)$, often alternating several times between edges coming from $\al$ or $\be$.  For example, to see that $\{1',4'\}$ is a block of $\al\be$, one needs to trace the following path (or its reverse):
\[
1' \xrightarrow{\ \be\ } 1'' \xrightarrow{\ \al\ } 2'' \xrightarrow{\ \be\ } 3'' \xrightarrow{\ \al\ } 4'' \xrightarrow{\ \be\ } 4'.
\]
\end{eg}

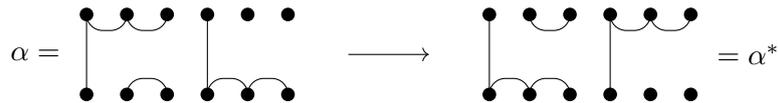
\begin{figure}[ht]
\begin{center}
\begin{tikzpicture}[scale=.53]
\begin{scope}[shift={(0,0)}]	
\uvs{1,...,6}
\lvs{1,...,6}
\uarc12
\uarc23
\darc23
\darc45
\darc56
\stline11
\stline44
\draw(0.6,1)node[left]{$\al=$};
\draw[->](7.5,1)--(9.5,1);
\end{scope}
\begin{scope}[shift={(10,0)}]	
\uvs{1,...,6}
\lvs{1,...,6}
\darc12
\darc23
\uarc23
\uarc45
\uarc56
\stline11
\stline44
\draw(6.4,1)node[right]{$=\al^*$};
\end{scope}
\end{tikzpicture}
\caption{A partition $\al\in\PP_6$ (left) and its image $\al^*$ (right) under the involution of $\PP_6$.}
\label{fig:PX3}
\end{center}
\end{figure}

We hope that Example \ref{eg:PX2} convinces the reader that products in $\PP_X$ can be `messy'.  Of course things get worse when we increase the number of vertices, even more so when $X$ is infinite~\cite{JE2014}.  Nevertheless, it is actually quite easy to see that $\PP_X$ is a regular $*$-semigroup.  The involution
\[
{}^*:\PP_X\to\PP_X:\al\mt\al^*
\]
can be defined diagrammatically as a reflection in a horizontal axis; see Figure \ref{fig:PX3}.  Formally, $\al^*$ is obtained from $\al$ by swapping dashed and undashed vertices, $x\leftrightarrow x'$.  It is not hard to see that
\[
(\al^*)^* = \al = \al\al^*\al \AND (\al\be)^*=\be^*\al^* \qquad\text{for all $\al,\be\in\PP_X$.}
\]

\newpage

At this point it is necessary to recall some further notation and terminology.  First, we say a non-empty subset $\es\not=A\sub X\cup X'$ is:
\bit
\item a \emph{transversal} if $A\cap X\not=\es$ and $A\cap X'\not=\es$ (i.e.~if $A$ contains both upper and lower vertices),
\item an \emph{upper non-transversal} if $A\sub X$ (i.e.~if $A$ contains only upper vertices), or
\item a \emph{lower non-transversal} if $A\sub X'$ (i.e.~if $A$ contains only lower vertices).
\eit
We can then describe a partition $\al\in\PP_X$ using a convenient two-line notation from \cite{EF2012}.  Specifically, we write
\[
\al = 	\begin{partn}{2} A_i&C_j\\ \hhline{~|-} B_i&D_k\end{partn}_{i\in I,\ j\in J,\ k\in K}
\]
to indicate that $\al$ has transversals $A_i\cup B_i'$ ($i\in I$), upper non-transversals $C_j$ ($j\in J$) and lower non-transversals $D_k$ ($k\in K$).  We often abbreviate this to $\al = \begin{partn}{2} A_i&C_j\\ \hhline{~|-} B_i&D_k\end{partn}$, with the indexing sets~$I$,~$J$ and~$K$ being implied, rather than explicitly listed.  When $X$ is finite we list the blocks of partitions, \emph{viz.}~$\al = \begin{partn}{6} A_1&\dots&A_q&C_1&\dots&C_s\\ \hhline{~|~|~|-|-|-} B_1&\dots&B_q&D_1&\dots&D_t\end{partn}$.  Thus, for example, the partition $\be\in\PP_6$ from Figure~\ref{fig:PX1} can be expressed as
\[
\be = \begin{partn}{3} 1&4,5,6&2,3\\ \hhline{~|~|-} 4,6&1,2,3&5\end{partn}.
\]
With this notation, the identity of $\PP_X$ is the partition $\id_X = \begin{partn}{1} x\\ \hhline{~} x\end{partn}_{x\in X}$, and the involution is given by
\[
\al = \begin{partn}{2} A_i&C_j\\ \hhline{~|-} B_i&D_k\end{partn}
\mt
\al^* = \begin{partn}{2} B_i&D_k\\ \hhline{~|-} A_i&C_j\end{partn}.
\]
Moreover, one can easily see that with $\al=\begin{partn}{2} A_i&C_j\\ \hhline{~|-} B_i&D_k\end{partn}$, and with the notation of Section \ref{sect:RSS}, we have
\[
\bd(\al) = \al\al^* = \begin{partn}{2} A_i&C_j\\ \hhline{~|-} A_i&C_j\end{partn}
\AND
\br(\al) = \al^*\al = \begin{partn}{2} B_i&D_k\\ \hhline{~|-} B_i&D_k\end{partn}.
\]
Thus, consulting Lemma \ref{lem:Green}, it follows that partitions $\al,\be\in\PP_X$ are $\R$-related (or $\L$-related) precisely when the `top halves' (or `bottom halves') of $\al,\be$ `match' in a sense that we again hope is clear.  Nevertheless, this `matching' can be formalised by defining some further notation.

For $\al\in\PP_X$, we define the \emph{(co)domain} and \emph{(co)kernel} of $\al$ by:
\begin{align*}
\dom(\al) &= \set{x\in X}{x \text{ belongs to a transversal of }\al}, \\
\codom(\al) &= \set{x\in X}{x' \text{ belongs to a transversal of }\al}, \\
\ker(\al) &= \set{(x,y)\in X\times X}{x\text{ and }y \text{ belong to the same block of }\al},\\
\coker(\al) &= \set{(x,y)\in X\times X}{x'\text{ and }y' \text{ belong to the same block of }\al}.
\end{align*}
The \emph{rank} of $\al$, denoted $\rank(\al)$, is the number of transversals of $\al$.  
Thus, $\dom(\al)$ and $\codom(\al)$ are subsets of $X$, $\ker(\al)$ and $\coker(\al)$ are equivalences on $X$, and $\rank(\al)$ is a cardinal between~$0$ and~$|X|$.  For example, with $\al\in\PP_6$ from Figure \ref{fig:PX1}, we have
\begin{align*}
\dom(\al) &= \{1,2,3,4\} , & \ker(\al) &= (1,2,3\mr\vert4\mr\vert5\mr\vert6) , \\
\codom(\al) &= \{1,4,5,6\} , & \coker(\al) &= (1\mr\vert2,3\mr\vert4,5,6) , & \rank(\al) &= 2,
\end{align*}
where we have indicated equivalences by listing their classes.
The various parts of the following result are contained in \cite{Wilcox2007,FL2011,ER2022}, though some of those papers use different terminology.

\newpage

\begin{lemma}\label{lem:Green_PX}
For $\al,\be\in\PP_X$, we have
\ben
\item $\al \R \be \iff [\dom(\al)=\dom(\be)$ and $\ker(\al)=\ker(\be)]$,
\item $\al \L \be \iff [\codom(\al)=\codom(\be)$ and $\coker(\al)=\coker(\be)]$,
\item $\al \D \be \iff \al \J \be \iff \rank(\al)=\rank(\be)$.  
\een
The ${\D}={\J}$-classes of $\PP_X$ are the sets
\[
D_\mu = D_\mu(\PP_X) = \set{\al\in\PP_X}{\rank(\al)=\mu}   \qquad\text{for cardinals $0\leq \mu\leq|X|$.}  
\]
Group $\H$-classes contained in $D_\mu$ are isomorphic to the symmetric group $\mathcal S_\mu$.  \epfres
\end{lemma}

Let us now return to the example products considered earlier.

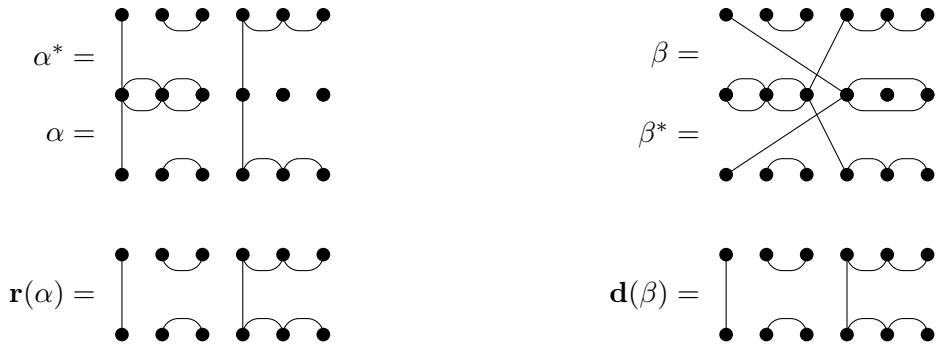
\begin{figure}[ht]
\begin{center}
\begin{tikzpicture}[scale=.53]
\begin{scope}[shift={(0,0)}]	
\uvs{1,...,6}
\lvs{1,...,6}
\darc12
\darc23
\uarc23
\uarc45
\uarc56
\stline11
\stline44
\draw(0.6,1)node[left]{$\al^*=$};
\end{scope}
\begin{scope}[shift={(0,-2)}]	
\uvs{1,...,6}
\lvs{1,...,6}
\uarc12
\uarc23
\darc23
\darc45
\darc56
\stline11
\stline44
\draw(0.6,1)node[left]{$\al=$};
\end{scope}
\begin{scope}[shift={(0,-6)}]	
\uvs{1,...,6}
\lvs{1,...,6}
\uarc23
\uarc45
\uarc56
\darc23
\darc45
\darc56
\stline11
\stline44
\draw(0.6,1)node[left]{$\br(\al)=$};
\end{scope}
\begin{scope}[shift={(15,0)}]	
\uvs{1,...,6}
\lvs{1,...,6}
\uarc23
\uarc45
\uarc56
\darc12
\darc23
\darc46
\stline14
\stline43
\draw(0.6,1)node[left]{$\be=$};
\end{scope}
\begin{scope}[shift={(15,-2)}]	
\uvs{1,...,6}
\lvs{1,...,6}
\darc23
\darc45
\darc56
\uarc12
\uarc23
\uarc46
\stline41
\stline34
\draw(0.6,1)node[left]{$\be^*=$};
\end{scope}
\begin{scope}[shift={(15,-6)}]	
\uvs{1,...,6}
\lvs{1,...,6}
\uarc23
\uarc45
\uarc56
\darc23
\darc45
\darc56
\stline11
\stline44
\draw(0.6,1)node[left]{$\bd(\be)=$};
\end{scope}
\end{tikzpicture}
\caption{The projections $\br(\al)=\al^*\al$ and $\bd(\be)=\be\be^*$, where $\al,\be\in\PP_6$ are as in Figure \ref{fig:PX1}.  For more information, see Examples \ref{eg:PX1} and \ref{eg:PX3}.}
\label{fig:PX4}
\end{center}
\end{figure}

\begin{eg}\label{eg:PX3}
Consider again the partitions $\al,\be\in\PP_6$ from Example \ref{eg:PX1} and Figure \ref{fig:PX1}.  The projections
\[
\br(\al) = \al^*\al \AND \bd(\be) = \be\be^*
\]
are calculated in Figure \ref{fig:PX4}.  Since we have $\br(\al)=\bd(\be)$, we see that $\al$ and $\be$ are composable in the groupoid $\G(\PP_6)$ from Definition \ref{defn:GS}, so that in fact $\al\be=\al\circ\be$.  This explains the `matching' phenomenon discussed above, and gives the reason for the `ease' of forming the product~$\al\be$.
\end{eg}

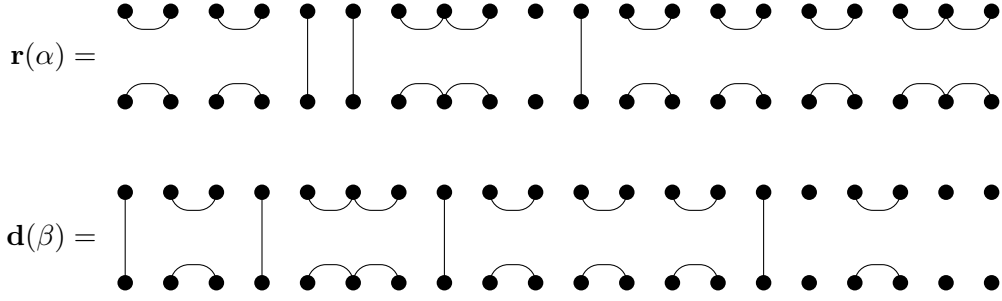
\begin{figure}[ht]
\begin{center}
\begin{tikzpicture}[scale=.6]
\begin{scope}[shift={(0,0)}]	
\uvs{1,...,20}
\lvs{1,...,20}
\uarcs{1/2,3/4,7/8,8/9,12/13,14/15,16/17,18/19,19/20}
\darcs{1/2,3/4,7/8,8/9,12/13,14/15,16/17,18/19,19/20}
\stlines{5/5,6/6,11/11}
\draw(0.6,1)node[left]{$\br(\al)=$};
\end{scope}
\begin{scope}[shift={(0,-4)}]	
\uvs{1,...,20}
\lvs{1,...,20}
\uarcs{2/3,5/6,6/7,9/10,11/12,13/14,17/18}
\darcs{2/3,5/6,6/7,9/10,11/12,13/14,17/18}
\stlines{1/1,4/4,8/8,15/15}
\draw(0.6,1)node[left]{$\bd(\be)=$};
\end{scope}
\end{tikzpicture}
\caption{The projections $\br(\al)=\al^*\al$ and $\bd(\be)=\be\be^*$, where $\al,\be\in\PP_{20}$ are as in Figure \ref{fig:PX2}.  For more information, see Examples \ref{eg:PX2} and \ref{eg:PX4}.}
\label{fig:PX5}
\end{center}
\end{figure}

\begin{figure}[ht]
\begin{center}
\begin{tikzpicture}[scale=.6]
\begin{scope}[shift={(0,-2)}]	
\darccols{2/3,5/6,6/7,9/10,11/12,13/14,17/18}{red}
\uarccols{1/2,3/4,5/6,7/8,8/9,12/13,14/15,16/17,18/19,19/20}{red}
\stlinecols{6/8,11/15}{red}
\darcxcol14{.7}{red}
\draw(0.6,1)node[left]{$e=$};
\end{scope}
\begin{scope}[shift={(0,0)}]	
\uvs{1,...,20}
\lvs{1,...,20}
\uarcs{1/2,2/3,4/5,7/8,8/9,9/10,11/12,12/13,13/14,14/15,15/16,16/17,17/18,18/19,19/20}
\uarcx36{.7}
\darcs{1/2,3/4,5/6,7/8,8/9,12/13,14/15,16/17,18/19,19/20}
\stlines{6/6,10/11}
\draw(0.6,1)node[left]{$\al'=$};
\draw(21,-1)node[right]{$\Pi(\al',e,\be')$};
\draw[|-|] (21,2)--(21,-4);
\end{scope}
\begin{scope}[shift={(0,-4)}]	
\uvs{1,...,20}
\lvs{1,...,20}
\uarcs{2/3,5/6,6/7,9/10,11/12,13/14,17/18}
\darcs{2/3,5/6,6/7,7/8,11/12,12/13,13/14,14/15,16/17,17/18,18/19,19/20}
\stlines{8/8,15/15}
\uarcx14{.7}
\darcx14{.7}
\draw(0.6,1)node[left]{$\be'=$};
\end{scope}
\begin{scope}[shift={(0,-8)}]	
\uvs{1,...,20}
\lvs{1,...,20}
\uarcs{1/2,2/3,4/5,7/8,8/9,9/10,11/12,12/13,13/14,14/15,15/16,16/17,17/18,18/19,19/20}
\uarcx36{.7}
\darcs{2/3,5/6,6/7,7/8,11/12,12/13,13/14,14/15,16/17,17/18,18/19,19/20}
\darcx14{.7}
\stlines{6/6,10/11}
\draw(0.6,1)node[left]{$\al'\circ e\circ\be'=$};
\end{scope}
\end{tikzpicture}
\caption{Top: partitions $\al',e,\be'\in\PP_{20}$, with the edges of $e$ shown in red.  Bottom: the composition $\al'\circ e\circ\be'$ in the groupoid $\G(\PP_{20})$.  Note that $\al'\circ e\circ\be' = \al\be$, where $\al,\be\in\PP_{20}$ are as in Figure \ref{fig:PX2}.  For more information, see Examples \ref{eg:PX2} and~\ref{eg:PX4}.}
\label{fig:PX6}
\end{center}
\end{figure}
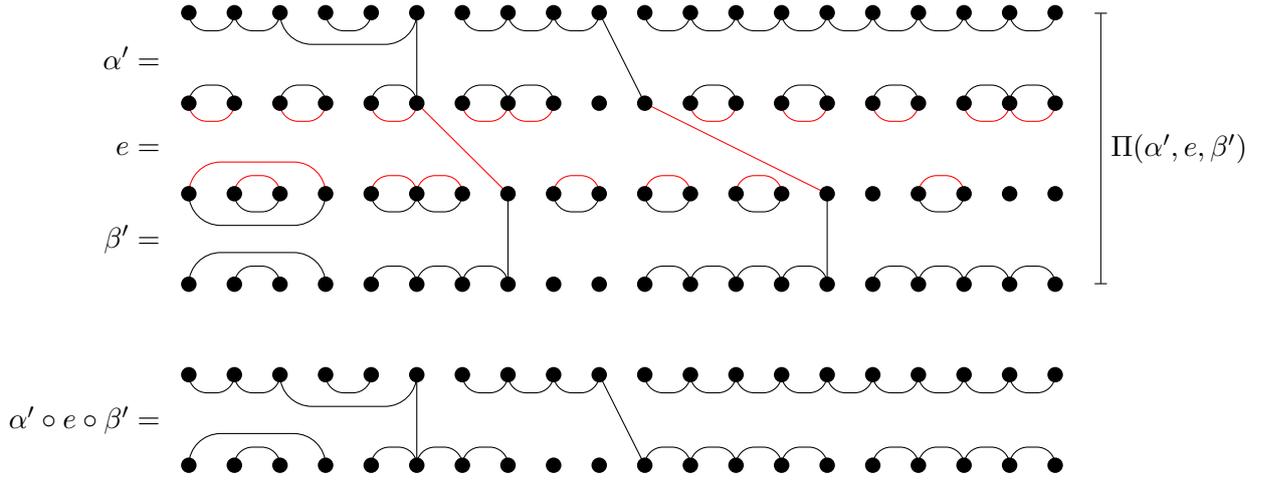

\begin{eg}\label{eg:PX4}
On the other hand, the partitions $\al,\be\in\PP_{20}$ from Example \ref{eg:PX2} and Figure \ref{fig:PX2} are not composable in the groupoid $\G(\PP_{20})$.  Indeed, one can check that the projections
\[
p = \br(\al) = \al^*\al \AND q = \bd(\be) = \be\be^*
\]
are as shown in Figure \ref{fig:PX5}, and we clearly do not have $\br(\al)=\bd(\be)$.  On the other hand, we can reduce the product $\al\be$ in $\PP_{20}$ to a composition $\al\be = \al' \circ e \circ \be'$ in $\G(\PP_{20})$, as in Remark \ref{rem:GS}.  Here we take $p' = q\th_p$ and $q'=p\th_q$, and then $\al'=\al\rest_{p'} = \al p'$, $\be'={}_{q'}\corest\be=q'\be$ and $e=p'q'=pq$ are as shown in Figure \ref{fig:PX6}.
\end{eg}

We conclude this section with a brief discussion of the relations $\leq$, $\leqF$ and $\F$ on the set $P=P(\PP_X)$ of projections of a partition monoid $\PP_X$, beginning with the case $|X|=2$.  %These relations were defined in \eqref{eq:leq} and~\eqref{eq:leqFF}.  

\begin{figure}[ht]
\begin{center}
\begin{tikzpicture}[scale=.4]
\begin{scope}[shift={(0,0)}]
\tuv1\tuv2\tlv1\tlv2\stline11\stline22 
\end{scope}
\begin{scope}[shift={(5,0)}]
\tuv1\tuv2\tlv1\tlv2\stline11
\end{scope}
\begin{scope}[shift={(10,0)}]
\tuv1\tuv2\tlv1\tlv2\stline11\stline22\uarc12\darc12 
\end{scope}
\begin{scope}[shift={(15,0)}]
\tuv1\tuv2\tlv1\tlv2\stline22 
\end{scope}
\begin{scope}[shift={(20,0)}]
\tuv1\tuv2\tlv1\tlv2 
\end{scope}
\begin{scope}[shift={(25,0)}]
\tuv1\tuv2\tlv1\tlv2\uarc12\darc12 
\end{scope}
\end{tikzpicture} 
\caption{The projections of the partition monoid $\PP_2$.}
\label{fig:PX11}
\end{center}
\end{figure}
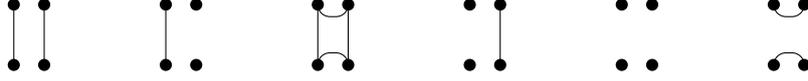

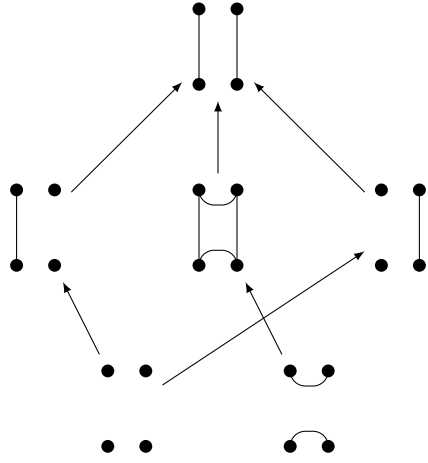
\begin{figure}[ht]
\begin{center}
\begin{tikzpicture}[scale=1.2]
\node (a) at (0,4) {\begin{tikzpicture}[scale=.5]\uv1\uv2\lv1\lv2\stline11\stline22\end{tikzpicture}};
\node (b) at (-2,2) {\begin{tikzpicture}[scale=.5]\uv1\uv2\lv1\lv2\stline11\end{tikzpicture}};
\node (c) at (0,2) {\begin{tikzpicture}[scale=.5]\uv1\uv2\lv1\lv2\stline11\stline22\uarc12\darc12\end{tikzpicture}};
\node (d) at (2,2) {\begin{tikzpicture}[scale=.5]\uv1\uv2\lv1\lv2\stline22\end{tikzpicture}};
\node (e) at (-1,0) {\begin{tikzpicture}[scale=.5]\uv1\uv2\lv1\lv2\end{tikzpicture}};
\node (f) at (1,0) {\begin{tikzpicture}[scale=.5]\uv1\uv2\lv1\lv2\uarc12\darc12\end{tikzpicture}};
\foreach \x/\y in {e/b,e/d,f/c,b/a,c/a,d/a} {\draw[-{latex}] (\x)--(\y);}
\end{tikzpicture}
\caption{Hasse diagram of the poset $P(\PP_2)$, with respect to the order $\leq$ given in \eqref{eq:leq}.  An arrow $p\to q$ means $p\leq q$.}
\label{fig:PX7}
\end{center}
\end{figure}

\begin{figure}[ht]
\begin{center}
\begin{tikzpicture}[scale=1.2]
\node (a) at (0,4) {\begin{tikzpicture}[scale=.5]\uv1\uv2\lv1\lv2\stline11\stline22\end{tikzpicture}};
\node (b) at (-2,2) {\begin{tikzpicture}[scale=.5]\uv1\uv2\lv1\lv2\stline11\end{tikzpicture}};
\node (c) at (0,2) {\begin{tikzpicture}[scale=.5]\uv1\uv2\lv1\lv2\stline11\stline22\uarc12\darc12\end{tikzpicture}};
\node (d) at (2,2) {\begin{tikzpicture}[scale=.5]\uv1\uv2\lv1\lv2\stline22\end{tikzpicture}};
\node (e) at (-1,0) {\begin{tikzpicture}[scale=.5]\uv1\uv2\lv1\lv2\end{tikzpicture}};
\node (f) at (1,0) {\begin{tikzpicture}[scale=.5]\uv1\uv2\lv1\lv2\uarc12\darc12\end{tikzpicture}};
\foreach \x/\y in {b/a,c/a,d/a,e/a,f/a,b/c,c/b,c/d,d/c,e/f,f/e,e/b,e/c,e/d,f/b,f/c,f/d} {\draw[-{latex}] (\x)--(\y);}
\begin{scope}[shift={(8,0)}]
\node (a) at (0,4) {\begin{tikzpicture}[scale=.5]\uv1\uv2\lv1\lv2\stline11\stline22\end{tikzpicture}};
\node (b) at (-2,2) {\begin{tikzpicture}[scale=.5]\uv1\uv2\lv1\lv2\stline11\end{tikzpicture}};
\node (c) at (0,2) {\begin{tikzpicture}[scale=.5]\uv1\uv2\lv1\lv2\stline11\stline22\uarc12\darc12\end{tikzpicture}};
\node (d) at (2,2) {\begin{tikzpicture}[scale=.5]\uv1\uv2\lv1\lv2\stline22\end{tikzpicture}};
\node (e) at (-1,0) {\begin{tikzpicture}[scale=.5]\uv1\uv2\lv1\lv2\end{tikzpicture}};
\node (f) at (1,0) {\begin{tikzpicture}[scale=.5]\uv1\uv2\lv1\lv2\uarc12\darc12\end{tikzpicture}};
\draw (b)--(c)--(d) (e)--(f);
\end{scope}
\end{tikzpicture}
\caption{Left: the relation $\leqF$ on $P(\PP_2)$; an arrow $p\to q$ means $p\leqF q$.  Right: the relation $\F$ on $P(\PP_2)$; an edge $p-q$ means $p\F q$.  These relations are given in \eqref{eq:leqFF}.  In both diagrams, loops are omitted.}
\label{fig:PX8}
\end{center}
\end{figure}
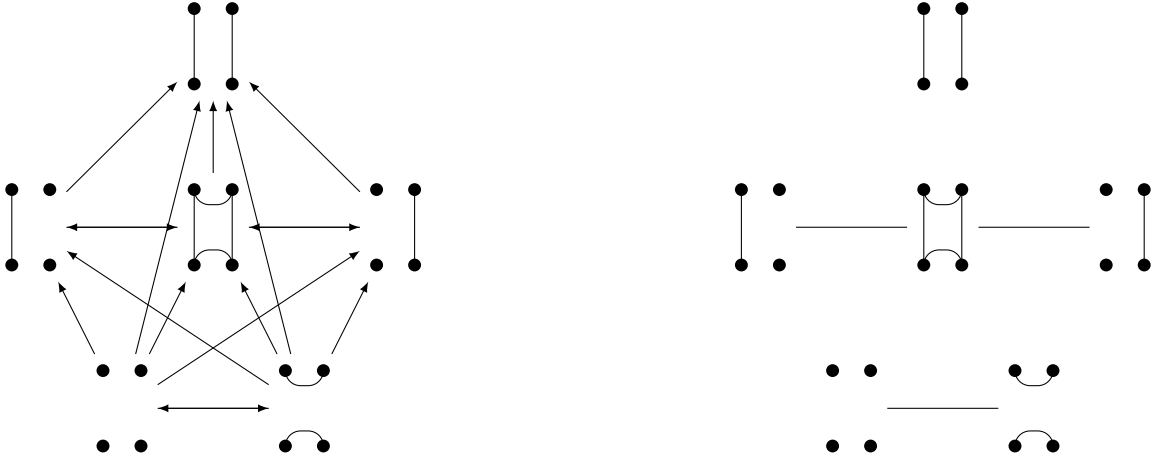

\begin{eg}\label{eg:P2}
Let $P=P(\PP_2)$ be the projection algebra of the partition monoid $\PP_2$.  The six elements of $P$ are shown in Figure \ref{fig:PX11}.
Figure \ref{fig:PX7} shows the partial order $\leq$ on $P(\PP_2)$; as usual for Hasse diagrams, only covering relationships are shown, and the rest can be deduced from transitivity (and reflexivity).  Figure \ref{fig:PX8} shows the relations~$\leqF$ and $\F$ on $P(\PP_2)$, and again one can see that these are not transitive.

In both Figures \ref{fig:PX7} and \ref{fig:PX8}, the elements of $P=P(\PP_2)$ have been arranged so that each row consists of $\D$-related elements; cf.~Lemma \ref{lem:Green_PX}.  So from top to bottom, the rows correspond to the projections from the $\D$-classes $D_2$, $D_1$ and $D_0$, where
\[
D_i=D_i(\PP_2)=\set{\al\in\PP_2}{\rank(\al)=i} \qquad\text{for $i=0,1,2$.}
\]
For $i=0,1,2$, we write $P_i = P\cap D_i$ for the set of projections from $D_i$.  The right-hand graph in Figure \ref{fig:PX8} is the graph $\Ga(\PP_2)$, as defined at the end of Section \ref{sect:RSS}.  We denote this by $\Ga=\Ga(\PP_2)$; so the vertex set of $\Ga$ is $P=P(\PP_2)$, and $\Ga$ has an undirected edge $\{p,q\}$ precisely when $p\not=q$ and $p\F q$.  
As in \eqref{eq:FD}, $\F$-related projections are $\D$-related.  Thus,
\[
\Ga=\Ga_0\sqcup\Ga_1\sqcup\Ga_2
\]
decomposes as the disjoint union of the induced subgraphs $\Ga_i$ on the vertex sets $P_i$, for $i=0,1,2$.  It is clear that each $\Ga_i$ is connected, but this is not necessarily the case for an arbitrary regular $*$-semigroup~$S$.
\end{eg}

\begin{figure}[ht]
\begin{center}
\includegraphics[width=0.9\textwidth]{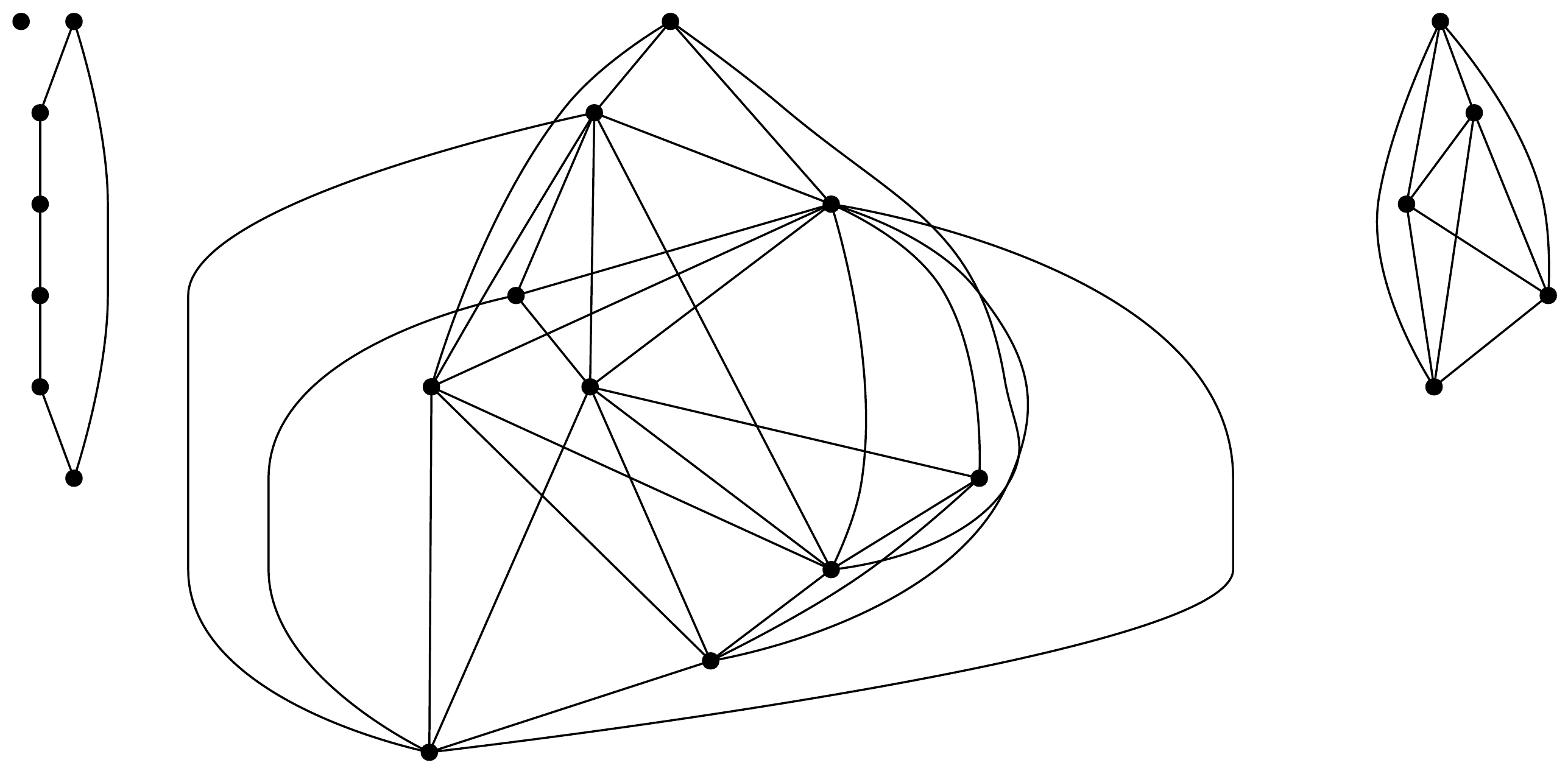}
\caption{The graph $\Ga(\PP_3)$, produced by GAP.  }
\label{fig:PX9}
\end{center}
\end{figure}

In general, for a regular $*$-semigroup $S$, we still have the graph $\Ga(S)$, and we still have the decomposition
\[
\Ga(S) = \bigsqcup_{D\in S/{\D}} \Ga(D),
\]
where $\Ga(D)$ is the induced subgraph of $\Ga(S)$ on vertex set $P(S)\cap D$, for each $\D$-class $D$ of $S$.  However, the subgraphs $\Ga(D)$ need not be connected in general.  For example, if $S$ is inverse, then the $\F$-relation is trivial; it follows that each $\Ga(D)$ is discrete (has empty edge set), and hence is disconnected if $D$ contains more than one idempotent.  As another example, if $S=A(\Ga)$ is an adjacency semigroup, as in Example \ref{eg:AG}, then the graph $\Ga(S)$ is simply $\Ga$ with a new vertex $0$ added; the $\D$-classes of $S$ are $D_0=\{0\}$ and $D_1=S\sm\{0\}=P\times P$, and we have $\Ga(D_1)=\Ga$.
On the other hand, it follows from results of \cite{EG2017} that the graphs $\Ga(D)$ are connected for any $\D$-class $D$ of a finite partition monoid $\PP_n$, and this turns out to be equivalent to certain facts about minimal idempotent-generation of the proper ideals of $\PP_n$.  We have seen this connectivity property in the case $n=2$, above.  Figure~\ref{fig:PX9} shows the graph $\Ga(\PP_3)$, produced using the Semigroups package for GAP \cite{Semigroups,GAP}.  From left to right, the connected components of this graph are the induced subgraphs $\Ga(D_3)$, $\Ga(D_2)$, $\Ga(D_1)$ and $\Ga(D_0)$.

As a special case, the induced subgraph $\Ga(D_{n-1})$ of the graph $\Ga(\PP_n)$ corresponding to the $\D$-class $D_{n-1}=D_{n-1}(\PP_n)$ has a very interesting structure.  First, one can easily check that the projections of $\PP_n$ contained in $D_{n-1}$ have the form
\[
\begin{tikzpicture}[scale=.5]
\begin{scope}[shift={(0,0)}]
\foreach \x in {1,3,4,5,7} {\uv\x\lv\x}
\foreach \x in {1,3,5,7} {\stline\x\x}
\foreach \x in {2,6} {\node()at(\x,2){$\cdots$};\node()at(\x,0){$\cdots$};}
\foreach \x/\y in {1/1,4/i,7/n} {\node()at(\x,2.5){\footnotesize $\y$};}
\draw(0.6,1)node[left]{$\pi_i=$};
\end{scope}
\begin{scope}[shift={(13,0)}]
\foreach \x in {1,3,4,5,7,8,9,11} {\uv\x\lv\x\stline\x\x}
\foreach \x in {2,6,10} {\node()at(\x,2){$\cdots$};\node()at(\x,0){$\cdots$};}
\foreach \x/\y in {1/1,4/j,8/k,11/n} {\node()at(\x,2.5){\footnotesize $\y$};}
\darc48
\uarc48
\draw(0.6,1)node[left]{and \qquad $\pi_{jk}=$};
\end{scope}
\end{tikzpicture} 
\]
for each $1\leq i\leq n$ and $1\leq j<k\leq n$.  It is also easy to check that the only non-identical $\F$-relations among these projections are
\[
\pi_i\F\pi_{ij}\F\pi_j.
\]
The graph $\Ga(D_{n-1})$ is shown in Figure \ref{fig:PX10} in the case $n=5$.  In the figure, vertices representing projections $\pi_i$ or $\pi_{jk}$ are labelled simply with the subscripts $i$ or $jk$.

\begin{figure}[ht]
\begin{center}
\begin{tikzpicture}[scale=3]
\tikzstyle{vertex}=[circle,draw=black, fill=white, inner sep = 0.07cm]
\draw (0,0) circle (1);
\node[vertex] (1) at (90:1) {  $1$ };
\node[vertex] (2) at (90-72:1) {  $2$ };
\node[vertex] (3) at (90-72-72:1) {  $3$ };
\node[vertex] (4) at (90-72-72-72:1) {  $4$ };
\node[vertex] (5) at (90-72-72-72-72:1) {  $5$ };
\node[vertex] (12) at (90-36:1) {\footnotesize  $12$ };
\node[vertex] (13) at (0.293890801	,0.095490177) {\footnotesize  $13$ };
\node[vertex] (14) at (-0.293894022	,0.095492517) {\footnotesize  $14$ };
\node[vertex] (15) at (90-36-72-72-72-72:1) {\footnotesize  $15$ };
\node[vertex] (23) at (90-36-72:1) {\footnotesize  $23$ };
\node[vertex] (24) at (0.181635098	,-0.250001636) {\footnotesize  $24$ };
\node[vertex] (25) at (0	,0.30901548) {\footnotesize  $25$ };
\node[vertex] (34) at (90-36-72-72:1) {\footnotesize  $34$ };
\node[vertex] (35) at (-0.181637088	,-0.25000019) {\footnotesize  $35$ };
\node[vertex] (45) at (90-36-72-72-72:1) {\footnotesize  $45$ };
\foreach \x/\y in {1/3,1/4,2/4,2/5,3/5} {\draw(\x)--(\x\y)--(\y);}
\end{tikzpicture}
\end{center}
\caption{The graph $\Ga(D_4)$, where $D_4=D_4(\PP_5)$.}
\label{fig:PX10}
\end{figure}

\newpage

The idempotent-generated subsemigroup $\E(\PP_X)=\la E(\PP_X)\ra=\la P(\PP_X)\ra$ of a partition monoid~$\PP_X$ was described for finite and infinite $X$ in \cite{JEpnsn} and \cite{EF2012}, respectively.  One of the main results of~\cite{JEpnsn} is that finite~$\E(\PP_n)$ is (minimally) generated as a monoid by the set
\[
\Om = \set{\pi_i,\pi_{jk}}{1\leq i\leq n,\ 1\leq j<k\leq n}
\]
of all projections of rank $n-1$, and that in fact
\[
\E(\PP_n) = \{\id_n\}\cup\Sing(\PP_n),
\]
where here ${\Sing(\PP_n)=\PP_n\setminus\mathcal S_n}$ is the \emph{singular ideal} of $\PP_n$, consisting of all non-invertible elements.  The second main result was a presentation for $\Sing(\PP_n)$ in terms of the above generating set $\Om$.  For example, one of the defining relations was
\begin{equation}\label{eq:triangle}
\pi_i \pi_{ij} \pi_j \pi_{jk} \pi_k \pi_{ki} \pi_i = \pi_i \pi_{ik} \pi_k \pi_{kj} \pi_j \pi_{ji} \pi_i \qquad\text{for distinct $1\leq i,j,k\leq n$,}
\end{equation}
where here we use symmetrical notation $\pi_{st}=\pi_{ts}$.  Examining Figure \ref{fig:PX10}, one can see that the two words in this equation represent two different triangular paths in the graph $\Ga(D_{n-1})$:
\[
i\rightsquigarrow j\rightsquigarrow k\rightsquigarrow i \AND i\rightsquigarrow k\rightsquigarrow j\rightsquigarrow i.
\]
(In representing the above paths we have omitted the intermediate vertices, so $i\rightsquigarrow j$ is shorthand for $i\to ij\to j$, and so on.)  Since paths in this graph correspond to lists of sequentially $\F$-related projections, the above words represent factorisations of the form described in Proposition~\ref{prop:ERSS}, and the corresponding tuples $(\pi_i,\pi_{ij},\pi_j,\ldots)$ and $(\pi_i,\pi_{ik},\pi_k,\ldots)$ are examples of what we will later call \emph{$P$-paths}; see Section \ref{subsect:PP}.  The two factorisations in \eqref{eq:triangle} represent the following partition, pictured in the case~${i<j<k}$:
\[
\begin{tikzpicture}[scale=.5]
\foreach \x in {1,3,4,5,7,8,9,11,12,13,15} {\uv\x\lv\x}
\foreach \x in {1,3,5,7,9,11,13,15} {\stline\x\x}
\stline8{12}
\stline{12}8
\foreach \x in {2,6,10,14} {\node()at(\x,2){$\cdots$};\node()at(\x,0){$\cdots$};}
\foreach \x/\y in {1/1,4/i,8/j,12/k,15/n} {\node()at(\x,2.5){\footnotesize $\y$};}
\end{tikzpicture} 
\]
%Much more could be said, but we will conclude our preliminary discussion of partition monoids here.  We will return to them at various stages throughout the text.

\section{Projection algebras}\label{sect:P}

One of the key ideas in this paper is that of a \emph{projection algebra}.  Such an algebra consists of a set $P$, along with a family $\th_p$ ($p\in P$) of unary operations, one for each element of $P$.  These operations are required to satisfy certain axioms, as listed in Definition \ref{defn:P} below, and are meant to abstractly model the unary algebras of projections of regular $*$-semigroups encountered in Section \ref{subsect:Sprelim}.  This is in fact not a new concept.  Indeed, projection algebras have appeared in a number of settings, under a variety of names, including the \emph{$P$-groupoids} of Imaoka \cite{Imaoka1983}, the \emph{$P$-sets} of Yamada \cite{Yamada1981}, (certain special) \emph{$P$-sets} of Nambooripad and Pastijn \cite{NP1985}, and the \emph{(left and right) projection algebras} of Jones~\cite{Jones2012}.  (We also mention Yamada's \emph{$p$-systems} \cite{Yamada1982}; these are defined in a very different way, but for a similar purpose.  Strictly speaking, Yamada's $P$-sets are different to Imaoka's $P$-groupoids, but are ultimately equivalent.)  We prefer Jones' terminology `projection algebras', since these are indeed algebras, not just sets, though we note that Jones considered \emph{binary} algebras rather than Imaoka's \emph{unary} algebras, which we use here; see Remark~\ref{rem:diamond}.  We also only use the term `groupoid' in its usual categorical sense.  Imaoka used it in one of its less-common meanings, stemming from the fact that a projection algebra also has a partially defined binary operation; curiously this partial operation plays no role in the current work.

Although the concept of a projection algebra is not new, here we build several new categorical structures on top of such algebras.  The most general situation is covered in Section \ref{sect:G} below, where we have an ordered groupoid $\G$ whose object set $v\G=P$ is a projection algebra (with certain algebraic and order-theoretic ties between $\G$ and $P$).  Section \ref{sect:CP} lays the foundation for this, by showing how to construct various categories directly from a projection algebra.  The main such construction is the \emph{chain groupoid} $\C=\C(P)$; this will be used in the definition of the groupoids in Section \ref{sect:G}, but we will uncover its deeper categorical significance in the sequel~\cite{Paper3}.

Before all of this, however, we need to discuss projection algebras themselves.
We begin in Section \ref{subsect:P} by giving/recalling their definition, and establishing some of their important properties.  %In Sections~\ref{subsect:PP} and~\ref{subsect:CP}, respectively, we introduce the \emph{path category} $\P=\P(P)$ and the \emph{chain groupoid} $\C=\C(P)$ of a projection algebra~$P$ (see Definitions \ref{defn:PP} and \ref{defn:CP}).  %The groupoid $\C$ is defined as a quotient of the category~$\P$ by a certain congruence $\approx$, whose definition requires the somewhat-technical notion of \emph{linked pairs}, which are the subject of Section \ref{sect:LP}.
We discuss the category $\PA$ of all projection algebras in Section \ref{subsect:PA}, where in Proposition~\ref{prop:Pfunctor} we construct a functor $\RSS\to\PA$ from the category of regular $*$-semigroups.  Finally, in Section \ref{subsect:PtoS} we give a new proof of the known fact \cite{Imaoka1983} that each (abstract) projection algebra is the unary algebra of projections of some regular $*$-semigroup.

\subsection{Definitions and basic properties}\label{subsect:P}

Here then is the key definition:

\begin{defn}\label{defn:P}
A \emph{projection algebra} is a set $P$, together with a collection of unary operations~$\th_p$~($p\in P$) satisfying the following axioms, for all $p,q\in P$:
\begin{enumerate}[label=\textup{\textsf{(P\arabic*)}},leftmargin=10mm]
\item \label{P1} $p\theta_p = p$,
\item \label{P2} $\theta_p\theta_p=\theta_p$,
\item \label{P3} $p\theta_q\theta_p=q\theta_p$,
\item \label{P4} $\theta_p\theta_q\theta_p=\theta_{q\theta_p}$,
\item \label{P5} $\theta_p\theta_q\theta_p\theta_q=\theta_p\theta_q$.
\end{enumerate}
The elements of a projection algebra are called \emph{projections}.
\end{defn}

Strictly speaking, one should refer to `a projection algebra $(P,\th)$', but we will almost always use the symbol $\th$ to denote the unary operations of such an algebra, and so typically refer to `a projection algebra $P$'.  When we need to refer simultaneously to more than one projection algebra (as for example in Definition \ref{defn:PA}), we will be careful to distinguish the operations.  

\begin{rem}\label{rem:diamond}
It is worth noting at the outset that the collection of unary operations~$\th_p$~($p\in P$) of a projection algebra $P$ could be replaced by a single binary operation $\diamond$, defined by $q\diamond p = q\th_p$, and resulting in what we will here call a \emph{binary projection algebra}.  
(For example, Jones took this binary approach in \cite{Jones2012}, in his study of the more general $P$-restriction semigroups, though his preferred symbol was $\star$, which we reserve for other uses throughout the paper.)
The binary form of the axioms \ref{P1}--\ref{P5} are as follows:
\begin{enumerate}[label=\textup{\textsf{(P\arabic*)$'$}},leftmargin=10mm]
\item \label{P1'} $p\diamond p = p$,
\item \label{P2'} $(q\diamond p)\diamond p=q\diamond p$,
\item \label{P3'} $(p\diamond q)\diamond p=q\diamond p$,
\item \label{P4'} $((r\diamond p)\diamond q)\diamond p=r\diamond (q\diamond p)$,
\item \label{P5'} $(((r\diamond p)\diamond q)\diamond p)\diamond q=(r\diamond p)\diamond q$.
\end{enumerate}
(Jones' axioms are slightly different to these, though they are equivalent, as he explains in~\cite[Section~7]{Jones2012}.)
Note that \ref{P2'} requires one more variable than \ref{P2} for its quantification; indeed, the latter is expressed in terms of equality of the term function $\th_p\th_p$ and the basic operation $\th_p$, which of course means that these functions agree on every input, i.e.~that $q\th_p\th_p=q\th_p$ for all $p,q\in P$.  Similar comments apply to \ref{P4'} and \ref{P5'}.  In this way, axioms \ref{P1'}--\ref{P5'} are all \emph{identities} in the signature of a single binary operation, meaning that binary projection algebras form a \emph{variety}.  On the other hand, axioms \ref{P1} and \ref{P3} do not have the form of identities, as~$p\th_p$ and~$p\th_q\th_p$ are not \emph{terms}.  In fact, unary projection algebras (as in Definition \ref{defn:P}) \emph{cannot} be defined by any set of identities, as they do not form varieties; indeed, the elements of a projection algebra are in one-one correspondence with its (unary) operations, and this property is not generally inherited by subalgebras.

Nevertheless, we do prefer the unary approach of Imaoka \cite{Imaoka1983}, for three main reasons:
\ben
%\item As mentioned above, the unary approach allows us to quantify some of the axioms over fewer variables.  
%For example, the binary form of \ref{P2} is $(q\diamond p)\diamond p = q\diamond p$, whereas the unary form simply says that the $\th$ operations are all idempotents.  Similar comments apply to~\ref{P4} and~\ref{P5}, which require three variables in binary form.
\item We feel that some of axioms \ref{P1}--\ref{P5} are clearer than their binary counterparts \mbox{\ref{P1'}--\ref{P5'}}.  
%For example, the binary form of~\ref{P4} is
%\[
%r\diamond(q\diamond p) = ((r\diamond p)\diamond q)\diamond p \qquad\text{for all $p,q,r\in P$.}
%\]
For example, and as discussed in more detail in Remark \ref{rem:G1}, axiom \ref{P4} can be thought of as a rule for `iterating' the unary operations, and we feel that this intuition is somewhat lost in the binary form \ref{P4'}, although this is of course entirely subjective.  %(We also note that Jones' axioms are not the direct translations \ref{P1'}--\ref{P5'}, as listed above, though they are equivalent, as he explains in~\cite[Section~7]{Jones2012}.)
\item The non-associativity of $\diamond$ means that brackets are necessary when working with $\diamond$-terms.  On the other hand, associativity of function composition allows us to dispense with bracketing when using the $\th$ maps, and this is a significant advantage when dealing with lengthy terms.
\item In later sections, we will consider groupoids $\G$ whose object sets are projection algebras,~${v\G=P}$.  One of the key tools when studying such groupoids are certain maps
\[
\Th_a=\th_{\bd(a)}\vt_a:P\to P
\]
that are built from the unary operations of $P$ and the maps $\vt_a:\bd(a)^\da\to\br(a)^\da$ defined in~\eqref{eq:vta}.  These $\Th$ maps could certainly be defined in terms of $\diamond$ and the $\vt$ maps, \emph{viz.}
\[
p\Th_a=(p\diamond\bd(a))\vt_a,
\]
but their definition as a composition $\Th_a=\th_{\bd(a)}\vt_a$ is more direct, and leads to more succinct statements and proofs, as for example with Proposition \ref{prop:G1} and Definition \ref{defn:PG}.
\een
\end{rem}

As we have already observed, axioms \ref{P1}--\ref{P5} are abstractions of the properties of projections of regular $*$-semigroups.  Specifically, we can make the following definition because of Lemma~\ref{lem:PS2}.  

\begin{defn}\label{defn:PS}
The projection algebra of a regular $*$-semigroup $S$ has:
\bit
\item underlying set $P(S) = \set{p\in S}{p^2=p=p^*}$, i.e.~the set of all projections of $S$, and
\item unary operations $\th_p$ $(p\in P)$, defined by $q\th_p=pqp$ for all $q\in P$.
\eit
\end{defn}

%`The projection algebra' of a regular $*$-semigroup always means the one constructed as in the above proposition.  It turns out that any (abstract) projection algebra (as in Definition \ref{defn:P}) is the projection algebra of some regular $*$-semigroup.  Indeed, this is well known (see for example \cite{Imaoka1983} or \cite{Jones2012}), but for convenience we will also give a proof in Section \ref{subsect:PtoS}; see Proposition \ref{prop:PtoS}.

For the rest of this section we fix a projection algebra $P$, as in Definition \ref{defn:P}.  In what follows, for $x,y\in P$, we write $x=_1y$ to indicate that $x=y$ by an application of \ref{P1}, and similarly for~$=_2$, and so on.

A number of relations on $P$ will play a crucial role in all that follows.  The first of these is defined by
\begin{align}
\label{eq:leqP} p\leq q &\Iff p=p\th_q.
\intertext{By \ref{P2}, it quickly follows that}
\label{eq:leqP2} p\leq q &\Iff p = r\th_q \qquad\text{for some $r\in P$.}
\end{align}

\begin{lemma}\label{lem:leqP}
$\leq$ is a partial order on $P$.
\end{lemma}

\pf
Reflexivity follows from \ref{P1}.  For anti-symmetry, suppose $p\leq q$ and $q\leq p$, so that $p=p\th_q$ and $q=q\th_p$.  Then
\[
p =_1 p\th_p = p\th_q\th_p =_3 q\th_p = q.
\]
For transitivity, suppose $p\leq q$ and $q\leq r$, so that $p=p\th_q$ and $q=q\th_r$.  Then
\[
p = p\th_q = p\th_{q\th_r} =_4 p\th_r\th_q\th_r.
\]
By \eqref{eq:leqP2}, this implies $p\leq r$.
\epf

The next lemma gives a simple criteria for the meet $p\wedge q$ of two projections $p,q\in P$ to exist, with respect to the partial order $\leq$.

\begin{lemma}\label{lem:meet}
If $p,q\in P$ satisfy $p\th_q=q\th_p$, then $p\wedge q$ exists in $(P,{\leq})$, and $p\wedge q=p\th_q=q\th_p$.
\end{lemma}

\pf
Let $r=p\th_q=q\th_p$.  Certainly $r\leq p,q$ by \eqref{eq:leqP2}.  Now suppose $s\leq p,q$; we must show that $s\leq r$, i.e.~that $s=s\th_r$.  Since $s\leq p,q$ it follows that $s=s\th_p=s\th_q$, and so
\[
s\th_r = s\th_{p\th_q} =_4s\th_q\th_p\th_q = s,
\]
as required.
\epf

\begin{rem}\label{rem:pqqp}
Lemma \ref{lem:meet} gave a sufficient condition for the meet $p\wedge q$ to exist, but this condition is not necessary.  Indeed, let $S=A(\Ga)$ be the adjacency semigroup of a symmetric, reflexive digraph $\Ga=(P,E)$, as in Example \ref{eg:AG}, and suppose distinct vertices $p,q\in P$ are connected by an edge in $\Ga$.  Then $p\wedge q=0$ (cf.~\eqref{eq:AGleq}), yet $p\th_q=q\not=p=q\th_p$ (cf.~\eqref{eq:AGth}).
\end{rem}

By \eqref{eq:leqP2}, the image of the basic operation $\th_p:P\to P$ is precisely
\begin{equation}\label{eq:imthp}
\im(\th_p) = p^\da = \set{q\in P}{q\leq p},
\end{equation}
the down-set of $p$ in the poset $(P,\leq)$.

\begin{lemma}\label{lem:thpthq}
If $p\leq q$, then $\th_p=\th_p\th_q = \th_q\th_p$.
\end{lemma}

\pf
Since $p=p\th_q$, we have $\th_p = \th_{p\th_q} =_4 \th_q\th_p\th_q$.  The claim then quickly follows from~\ref{P2}.
\epf

An equally important role will be played by two further relations, $\leqF$ and $\F$, defined as follows.  For $p,q\in P$, we say that
\begin{equation}\label{eq:leqF}
p \leqF q \Iff p = q\th_p .
\end{equation}
By \ref{P1}, $\leqF$ is reflexive, but it need not be symmetric or transitive.  We also define
\begin{equation}\label{eq:F}
{\F} = {\leqF}\cap{\geqF},
\end{equation}
which is the largest symmetric (and reflexive) relation contained in $\leqF$.  So
\[
p \F q \Iff p = q\th_p \text{ and } q=p\th_q.
\]

\newpage

\begin{lemma}\label{lem:pqp}
For any $p,q\in P$,
\ben\bmc2
\item \label{pqP1} $p\th_q\leqF p$,
\item \label{pqP2} $p\leq q\implies p\leqF q$,
\item \label{pqp1} $p\leqF q \implies \th_p=\th_p\th_q\th_p$.  
\item \label{pqp2} $p\F q\implies\th_p=\th_p\th_q\th_p$ and $\th_q=\th_q\th_p\th_q$.
\emc\een
\end{lemma}

\pf
\firstpfitem{\ref{pqP1}}  We have $p\th_{p\th_q} =_4 p\th_q\th_p\th_q =_3 q\th_p\th_q =_3 p\th_q$.

\pfitem{\ref{pqP2}}  If $p\leq q$, then $p=p\th_q$, so $q\th_p = q\th_{p\th_q} =_4 q\th_q\th_p\th_q =_1 q\th_p\th_q =_3 p\th_q = p$.

\pfitem{\ref{pqp1} and \ref{pqp2}}  These follow immediately from \ref{P4}.
\epf

\begin{rem}
The partial order $\leq$ from \eqref{eq:leqP} was used by Imaoka in \cite{Imaoka1983}, where it was defined by
\[
p\leq q \Iff p=p\th_q=q\th_p.
\]
Note that part \ref{pqP2} of the lemma just proved means that the `$=q\th_p$' part of Imaoka's definition is superfluous.  Jones also used the order $\leq$ in \cite{Jones2012}, defined exactly as in \eqref{eq:leqP}, albeit in binary form, $p\leq q \iff p=p\diamond q$ (cf.~Remark~\ref{rem:diamond}).  We are not aware of any previous use of the $\leqF$ relation in the literature, however, despite its central importance in the current work.  (In binary form we have $p\leqF q \iff p=q\diamond p$).
\end{rem}

\begin{rem}\label{rem:rev}
Some of the above implications are reversible.  For example, we have
\ben
\item \label{rev1} $p\leq q \iff \th_p=\th_p\th_q = \th_q\th_p \iff \th_p=\th_p\th_q$ (cf.~Lemma \ref{lem:thpthq}), and
\item \label{rev2} $p\leqF q \iff \th_p=\th_p\th_q\th_p$ (cf.~Lemma \ref{lem:pqp}\ref{pqp1}).
\een
Indeed, to establish \ref{rev1}, it suffices (by Lemma \ref{lem:thpthq}) to note that
\[
\th_p=\th_p\th_q \Implies p =_1 p\th_p = p\th_p\th_q =_1 p\th_q \Implies p\leq q.
\]
For the backward implication in \ref{rev2}, we have
\[
\th_p=\th_p\th_q\th_p \Implies  p =_1 p\th_p = p\th_p\th_q\th_p =_1 p\th_q\th_p =_3 q\th_p \Implies p\leqF q.
\]
\end{rem}

The following simple result will be crucial in later constructions.  Roughly speaking, it says that any pair of projections $p,q\in P$ sit above an $\F$-related pair $p',q'\in P$.

\begin{lemma}\label{lem:p'q'}
Let $p,q\in P$ be arbitrary, and let $p'=q\th_p$ and $q'=p\th_q$.  Then
\[
p'\leq p \COMMA q'\leq q \AND p'\F q'.
\]
\end{lemma}

\pf
We obtain $p'\leq p$ and $q'\leq q$ directly from \eqref{eq:leqP2}.  We also have
\[
%p'\th_{q'} = q\th_p\th_{p\th_q} =_4 q\th_p\th_q\th_p\th_q =_5 q\th_p\th_q =_3 p\th_q = q' \ANDSIM q'\th_{p'} = p'.  \qedhere
p'\th_{q'} = q\th_p\th_{p\th_q} =_4 q\th_p\th_q\th_p\th_q =_3 p\th_q\th_p\th_q =_3 q\th_p\th_q =_3 p\th_q = q' \ANDSIM q'\th_{p'} = p'.  \qedhere
\]
\epf

Although $\leqF$ need not be transitive, $\leq$ and $\leqF$ have some `transitivity-like' properties in combination:

\begin{lemma}\label{lem:pqr}
If $p,q,r\in P$ are such that $p\leq q\leqF r$ or $p\leqF q\leq r$, then 
\ben\bmc2
\item \label{pqr1} $p\leqF r$,
\item \label{pqr2} $p \F p\th_r$.
\emc\een
\end{lemma}

\pf
Throughout the proof we assume that $p\leq q\leqF r$.  The arguments for the case of $p\leqF q\leq r$ are similar, and are omitted.

Since $p\leq q$ and $q\leqF r$, we have $p=p\th_q$ and $q=r\th_q$.  Since also $p\leqF q$ by Lemma \ref{lem:pqp}\ref{pqP2}, we also have $p=q\th_p$.

\pfitem{\ref{pqr1}}  We must show that $p=r\th_p$.  Since $p\leq q$, Lemma \ref{lem:thpthq} gives $\th_p=\th_q\th_p$.  Combining all of the above gives $r\th_p = r\th_q\th_p = q\th_p = p$.

\pfitem{\ref{pqr2}}  We have $p\th_r\leqF p$ by Lemma \ref{lem:pqp}\ref{pqP1}, so we just need to show that $p\leqF p\th_r$, i.e.~that $p = (p\th_r)\th_p$.  But $p\th_r\th_p = r\th_p$ by \ref{P3}, and we saw in the previous part that $r\th_p=p$.
% This follows from the axioms and the above-mentioned consequences of $p\leq q\leqF r$:
%\[
%p\th_r\th_p =_3 r\th_p = r\th_{p\th_q} =_4 r\th_q\th_p\th_q = q\th_p\th_q =_3 p\th_q = p.  \qedhere
%\]
\epf

\begin{rem}\label{rem:rels}
Most of the results of this section have concerned properties of the $\leq$, $\leqF$ and $\F$ relations on a projection algebra $P$.  These relations were defined in \eqref{eq:leqP}, \eqref{eq:leqF} and \eqref{eq:F}, from the algebraic structure of $P$, i.e.~in terms of the $\th$ operations.  One might wonder if, conversely, the $\leq$, $\leqF$ and $\F$ relations uniquely determine the algebraic structure of $P$.  This turns out not to be the case, however.  Indeed, Example \ref{eg:Kinyon} and Remark \ref{rem:Kinyon2} give two projection algebras on the same underlying set, with the same $\leq$, $\leqF$ and $\F$ relations, but with different $\th$ operations.
\end{rem}

We now record two simple results that will be used in Section \ref{subsect:PtoS}.

\begin{lemma}\label{lem:thp=thq}
For $p,q\in P$, we have $\th_p=\th_q \iff p=q$.
\end{lemma}

\pf
If $\th_p=\th_q$, then $p =_1 p\th_p\th_p = p\th_q\th_p =_3 q\th_p = q\th_q =_1 q$.  The converse is trivial.
%
%Only the forward direction needs proof, so suppose $\th_p=\th_q$.  Then
%\[
%p =_1 p\th_p = p\th_q \implies p\leq q \ANDSIM q\leq p.  \qedhere
%\]
\epf

%\begin{lemma}\label{lem:P}
%For any $p_1,\ldots,p_k,q\in P$ we have
%\[
%\th_{q\th_{p_1}\cdots\th_{p_k}} = \th_{p_k}\cdots\th_{p_1}\th_q\th_{p_1}\cdots\th_{p_k}.
%\]
%\end{lemma}
%
%\pf
%This follows by iterating \ref{P4}.
%\epf

\begin{lemma}\label{lem:albeal}
If $\al = \th_{p_1}\cdots\th_{p_k}$ and $\be = \th_{p_k}\cdots\th_{p_1}$, for some $p_1,\ldots,p_k\in P$, then $\al=\al\be\al$ and~$\be=\be\al\be$.
%we have
%%and with $\al=\th_{p_1}\cdots\th_{p_k}$ and $\be = \th_{p_k}\cdots\th_{p_1}$, we have $\al\be\al=\al$.
%\[
%\th_{p_1}\cdots\th_{p_k} \cdot \th_{p_k}\cdots\th_{p_1} \cdot \th_{p_1}\cdots\th_{p_k} = \th_{p_1}\cdots\th_{p_k} .
%\]
\end{lemma}

\pf
By symmetry it suffices to show that $\al=\al\be\al$, and we do this by induction on $k$.  The $k=1$ case follows immediately from \ref{P2}.  %and $k=2$ cases follow from \ref{P2} and \ref{P5}, respectively.  For $k\geq3$ we have
Now suppose $k\geq2$, and write
\[
%\al = \th_{p_1}\cdots\th_{p_k} \COMMA \be = \th_{p_k}\cdots\th_{p_1} \COMMA 
\al' = \th_{p_1\th_{p_2}}\th_{p_3}\cdots\th_{p_k} \AND \be' = \th_{p_k}\cdots\th_{p_3}\th_{p_1\th_{p_2}}.
\]
Then $\al = \th_{p_1}\th_{p_2}\th_{p_3}\cdots\th_{p_k} =_5 \th_{p_1}\th_{p_2}\th_{p_1}\th_{p_2}\th_{p_3}\cdots\th_{p_k} =_4 \th_{p_1}\cdot \th_{p_1\th_{p_2}}\th_{p_3}\cdots\th_{p_k} = \th_{p_1} \cdot \al'$, and
\begin{align*}
\be\al = \th_{p_k}\cdots\th_{p_3}\th_{p_2}\th_{p_1} \cdot \th_{p_1}\th_{p_2}\th_{p_3}\cdots\th_{p_k} 
&=_2 \th_{p_k}\cdots\th_{p_3} \cdot\th_{p_2} \th_{p_1}\th_{p_2} \cdot\th_{p_3}\cdots\th_{p_k} \\
&=_4 \th_{p_k}\cdots\th_{p_3} \cdot \th_{p_1\th_{p_2}} \cdot\th_{p_3}\cdots\th_{p_k} \\
&=_2 \th_{p_k}\cdots\th_{p_3} \cdot \th_{p_1\th_{p_2}}\th_{p_1\th_{p_2}} \cdot\th_{p_3}\cdots\th_{p_k} = \be'\al'.
\end{align*}
Combining the last two calculations with $\al'\be'\al'=\al'$ (which holds by induction), we obtain
\[
\al\be\al = \th_{p_1}\cdot \al'\be'\al' = \th_{p_1}\cdot\al' = \al.  \qedhere
\]
\epf

\begin{rem}
It is worth noting that the proof of Lemma \ref{lem:albeal} was the only place in Section~\ref{subsect:P} that we used \ref{P5}.
\end{rem}

%It remains an open problem to determine the extent to which a projection algebra $P$ (including its $\th$ operations) can be somehow characterised by the purely order-theoretic properties of the relations $\leq$ and $\leqF$.  The next result is included in case it is of use in such a characterisation, but we note that the relations $\leq$ and $\leqF$ alone do not uniquely determine the operations of a projection algebra; see Example \ref{eg:Kinyon} and Remark \ref{rem:Kinyon}.
%
%\begin{lemma}
%If $p,q\in P$ are such that $p\leqF q$, then we have $p\F p'\leq q$ for some $p'\in P$.
%\end{lemma}
%
%\pf
%It is a routine matter to verify that $p'=p\th_q$ has the stated properties.
%\epf

%\end{document}

\subsection{The category of projection algebras}\label{subsect:PA}

We have now proved a number of results on projection algebras themselves, and we now turn to the category they form.

\begin{defn}\label{defn:PA}
We denote by $\PA$ the (large) category of projection algebras.  A morphism in $\PA$ is a map $\phi:P\to P'$ (for projection algebras $P$ and $P'$) satisfying
\[
(p\th_q)\phi = (p\phi)\th'_{q\phi} \qquad\text{for all $p,q\in P$.}
\]
Here we use $\th$ and $\th'$ to denote the unary operations on $P$ and $P'$, respectively.  
\end{defn}

One can also think of a projection algebra morphism as a morphism of binary algebras, in the sense discussed in Remark \ref{rem:diamond}.  Indeed, using $\diamond$ and $\diamond'$ to denote the binary operations of $P$ and $P'$ (as in that remark), $\phi:P\to P'$ is a projection algebra morphism if and only if
\[
(p\diamond q)\phi = (p\phi) \diamond' (q\phi) \qquad\text{for all $p,q\in P$.}
\]
In the next result we also write~$\leq$ and~$\leq'$ for the partial orders on $P$ and $P'$, as in \eqref{eq:leqP}, and similarly for the $\leqF$ and $\F$ relations from \eqref{eq:leqF} and \eqref{eq:F}.

\begin{lemma}\label{lem:PP'}
If $\phi:P\to P'$ is a projection algebra morphism, then for any $p,q\in P$,
\[
p\leq q \implies p\phi \leq' q\phi
\COMMA
p\leqF q \implies p\phi \leqF' q\phi
\AND
p\F q \implies p\phi \F' q\phi.
\]
\end{lemma}

\pf
The proofs are all essentially the same, so we just prove the first.  For this we have
\[
p\leq q \Implies p = p\th_q \Implies p\phi = (p\th_q)\phi = (p\phi)\th'_{q\phi} \Implies p\phi \leq' q\phi.  \qedhere
\]
\epf

\begin{rem}
The first implication in Lemma \ref{lem:PP'} says that any projection algebra morphism $P\to P'$ is a poset morphism $(P,{\leq})\to(P',{\leq'})$.  The converse of this is not true, however, as we will see in Example \ref{eg:Kinyon} and Remark \ref{rem:Kinyon1}.
\end{rem}

The construction of the projection algebra $P(S)$ from a regular $*$-semigroup $S$ (cf.~Definition~\ref{defn:PS}) can be thought of as an object map from the category~$\RSS$ of regular $*$-semigroups to the category $\PA$ of projection algebras.  It is easy to see that any $*$-morphism $\phi:S\to S'$ maps projections to projections, so we can define $P(\phi)=\phi|_{P(S)}:P(S)\to P(S')$.  The next result shows that $P$ (interpreted in this way) is a functor:

\begin{prop}\label{prop:Pfunctor}
$P$ is a functor $\RSS\to\PA$.  
\end{prop}

\pf
It remains to check that $P(\phi):P(S)\to P(S')$, as above, is a projection algebra morphism.  And indeed, for any $p,q\in P(S)$ we have
\[
(q\th_p) \phi = (pqp)\phi = (p\phi)(q\phi)(p\phi) = (q\phi) \th_{p\phi}.  \qedhere
\]
\epf

\subsection[From projection algebras to regular $*$-semigroups]{\boldmath From projection algebras to regular $*$-semigroups}\label{subsect:PtoS}

We have already seen that any regular $*$-semigroup $S$ naturally defines a projection algebra $P(S)$.  We now give a construction to show that we can also go in the reverse direction; Proposition~\ref{prop:PtoS} below shows that any (abstract) projection algebra is the projection algebra of some regular $*$-semigroup.  As we have already noted, this is well known, and follows from results in \cite{Imaoka1983,Jones2012}; see also \cite{Paper2,Paper3} for other proofs.  We include the present construction for completeness, and because it is considerably simpler than those from \cite{Imaoka1983,Jones2012,Paper2,Paper3}, as those papers concern larger problems than simply exhibiting a regular $*$-semigroup with a given projection algebra.

Consider an arbitrary semigroup $S$, and let $S^\opp$ be the \emph{opposite semigroup}.  So the underlying set of $S^\opp$ is $S$, and the product $\star$ in $S$ is defined by $x\star y = yx$.  Multiplication in the direct product $S\times S^\opp$ is given by the rule $(x,y)(u,v) = (xu,vy)$.  It is well known, and easy to see, that the operation $(x,y)^* = (y,x)$ is an involution on $S\times S^\opp$, meaning that the laws $(a^*)^*=a$ and $(ab)^*=b^*a^*$ are satisfied.

Now consider a projection algebra $P$.  Since an operation $\th_p$ is a map $P\to P$, we can think of it as an element of the \emph{full transformation monoid} $\T_P$.  (This is the monoid of \emph{all} maps~$P\to P$, under composition.)  We define $\FF_P$ to be the subsemigroup of the direct product $\T_P\times\T_P^\opp$ generated by all pairs $(\th_p,\th_p)$:
\[
\FF_P = \pres{(\th_p,\th_p)}{p\in P} = \set{(\th_{p_1}\cdots\th_{p_k},\th_{p_k}\cdots\th_{p_1})}{k\geq1,\ p_1,\ldots,p_k\in P} \leq \T_P\times\T_P^\opp.
\]
Since $\FF_P$ is closed under the ${}^*$ operation on $\T_P\times\T_P^\opp$ (defined above), it is an involution semigroup.

\begin{prop}\label{prop:PtoS}
If $P$ is an (abstract) projection algebra, then $\FF_P$ is a regular $*$-semigroup with projection algebra (isomorphic to) $P$.
\end{prop}

\pf
We have already observed that the laws $(a^*)^*=a$ and $(ab)^*=b^*a^*$ are satisfied.  For the remaining law, let $a=(\al,\be)\in\FF_P$, so that $\al=\th_{p_1}\cdots\th_{p_k}$ and $\be=\th_{p_k}\cdots\th_{p_1}$ for some $p_1,\ldots,p_k\in P$.  Using Lemma \ref{lem:albeal}, we have $aa^*a = (\al\be\al,\be\al\be) = (\al,\be) = a$.

Now that we know $\FF_P$ is a regular $*$-semigroup, we claim that its set of projections is
\begin{equation}\label{eq:PFP}
P(\FF_P) = \set{\ol p}{p\in P}, \qquad\text{where for simplicity we write $\ol p =(\th_p,\th_p)$.}
\end{equation}
Clearly each $\ol p$ is a projection (cf.~\ref{P2}).  Conversely, Lemma \ref{lem:PS1}\ref{PS11} tells us that each projection is of the form $aa^*$ for some $a\in\FF_P$.  Again writing $a=(\th_{p_1}\cdots\th_{p_k},\th_{p_k}\cdots\th_{p_1})$, we see that
%it follows from \ref{P2} and several applications of \ref{P4} that
\begin{align*}
aa^* &= (\th_{p_1}\cdots\th_{p_k}\cdot\th_{p_k}\cdots\th_{p_1},\th_{p_1}\cdots\th_{p_k}\cdot\th_{p_k}\cdots\th_{p_1})\\
&=_2 (\th_{p_1}\cdots\th_{p_{k-1}}\th_{p_k}\th_{p_{k-1}}\cdots\th_{p_1},\th_{p_1}\cdots\th_{p_{k-1}}\th_{p_k}\th_{p_{k-1}}\cdots\th_{p_1}) \\
&=_4 (\th_{p_k\th_{p_{k-1}}\cdots\th_{p_1}},\th_{p_k\th_{p_{k-1}}\cdots\th_{p_1}}) = \ol{p_k\th_{p_{k-1}}\cdots\th_{p_1}}
\end{align*}
%\[
%aa^* = (\th_{p_1}\cdots\th_{p_{k-1}}\th_{p_k}\th_{p_{k-1}}\cdots\th_{p_1},\th_{p_1}\cdots\th_{p_{k-1}}\th_{p_k}\th_{p_{k-1}}\cdots\th_{p_1}) = (\th_{p_k\th_{p_{k-1}}\cdots\th_{p_1}},\th_{p_k\th_{p_{k-1}}\cdots\th_{p_1}})
%\]
has the claimed form.  It remains to check that the map
\[
\phi:P\to P(\FF_P) : p \mt \ol p
\]
is a projection algebra isomorphism.  Surjectivity of $\phi$ follows from \eqref{eq:PFP}, and injectivity from Lemma~\ref{lem:thp=thq}.  It remains to check that $\phi$ is a projection algebra morphism, i.e.~that
$\ol{q\th_p} = \ol q\th_{\ol p}$ for all $p,q\in P$,
and for this we have
\[
\ol{q\th_p} = (\th_{q\th_p},\th_{q\th_p}) =_4 (\th_p\th_q\th_p,\th_p\th_q\th_p) = (\th_p,\th_p)(\th_q,\th_q)(\th_p,\th_p) = \ol p\ \ol q\ \ol p = \ol q\th_{\ol p}.  \qedhere
\]
\epf

\begin{rem}
We will see in \cite{Paper2} that the semigroup $\FF_P$ constructed above is (up to isomorphism) the unique fundamental, idempotent-generated regular $*$-semigroup with projection algebra (isomorphic to)~$P$.
\end{rem}

\begin{rem}\label{rem:notFfunctor}
One might wonder if the construction $P\mt\FF_P$ is the object part of a functor $\PA\to\RSS$.  This turns out not to be the case.  Indeed, if there was such a functor, any morphism $\phi:P\to P'$ in $\PA$ would need to induce a $*$-morphism $\Phi:\FF_P\to\FF_{P'}$ under which $(\th_p,\th_p)\Phi=(\th'_{p\phi},\th'_{p\phi})$ for all $p\in P$.  This map being well defined is equivalent to the implication
\begin{equation}\label{eq:pq=>p'q'}
\th_{p_1}\cdots\th_{p_k} = \th_{q_1}\cdots\th_{q_l} \text{ \ (in $P$)} \Implies \th'_{p_1\phi}\cdots\th'_{p_k\phi} = \th'_{q_1\phi}\cdots\th'_{q_l\phi} \text{ \ (in $P'$)} \qquad\text{for all $p_i,q_i\in P$.}
\end{equation}
To see that this does not always hold, let $P=\{p,q\}$ be the projection algebra in which $\th_p$ and $\th_q$ are both constant maps, and let $P'=P(\PP_2)$ be the projection algebra of the partition monoid~$\PP_2$; cf.~Example \ref{eg:P2}.  One can easily check that the mapping
\[
p \mt p' = \custpartn{1,2}{1,2}{\stline11}
\AND
q \mt q' = \custpartn{1,2}{1,2}{\stline11\stline22\uarc12\darc12}
\]
defines a projection algebra morphism $P\to P'$.  Of course we have $\th_p\th_q=\th_q$ in $P$, but we claim that $\th'_{p'}\th'_{q'}\not=\th'_{q'}$ in $P'$, witnessing the failure of \eqref{eq:pq=>p'q'} in this case.  Indeed, with $r = \custpartn{1,2}{1,2}{\stline22}\in P'$, we have
\[
r\th'_{p'}\th'_{q'} = \custpartn{1,2}{1,2}{\uarc12\darc12} \qquad\text{but}\qquad r\th'_{q'} = \custpartn{1,2}{1,2}{\stline11\stline22\uarc12\darc12} (= q').
\]
In the sequel \cite{Paper3}, we will exhibit a functor $\PA\to\RSS$ that is in fact a left adjoint to the (forgetful) functor $P:\RSS\to\PA$ from Proposition \ref{prop:Pfunctor}.  This shows that $\PA$ is (isomorphic to) a coreflective subcategory of $\RSS$, and demonstrates the existence of a free regular $*$-semigroup with a given projection algebra.  We will say more about this in Section \ref{subsect:PGP}.
\end{rem}

\subsection{Another example}\label{subsect:Kinyon}

We conclude Section \ref{sect:P} with another interesting projection algebra, again included to illustrate the above ideas, and highlight some subtleties; see Remarks \ref{rem:Kinyon1}, \ref{rem:Kinyon2} and \ref{rem:Cfunctor}.  The example was communicated to us by Michael Kinyon, who discovered it using Mace4 \cite{Mace4}.  

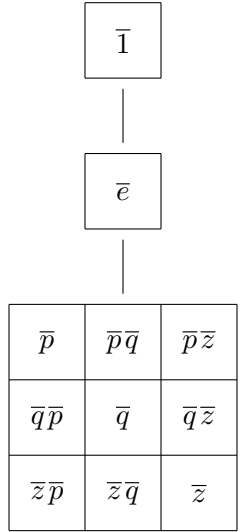
\begin{figure}[ht]
\begin{center}
\begin{tikzpicture}
\node (1) at (0,5) {
\begin{tikzpicture}
\foreach \x in {0,1} {\draw(0,\x)--(1,\x); \draw(\x,0)--(\x,1);}
\node()at(.5,.5){$\ol 1$};
\end{tikzpicture}
};
\node (2) at (0,3) {
\begin{tikzpicture}
\foreach \x in {0,1} {\draw(0,\x)--(1,\x); \draw(\x,0)--(\x,1);}
\node()at(.5,.5){$\ol e$};
\end{tikzpicture}
};
\node (3) at (0,0) {
\begin{tikzpicture}
\foreach \x in {0,1,2,3} {\draw(0,\x)--(3,\x); \draw(\x,0)--(\x,3);}
\foreach \x/\y/\z in {
0/0/\ol z\;\! \ol p,
0/1/\ol q\;\! \ol p,
0/2/\ol p,
1/0/\ol z\;\! \ol q,
1/1/\ol q,
1/2/\ol p\;\! \ol q,
2/0/\ol z,
2/1/\ol q\;\! \ol z,
2/2/\ol p\;\! \ol z,
} {\node()at(\x+.5,\y+.5){$\z$};}
\end{tikzpicture}
};
\draw(1)--(2)--(3);
\end{tikzpicture}
\caption{Egg-box diagram for the regular $*$-semigroup $\FF_P$, for the projection algebra $P=\{1,e,p,q,z\}$ from Example \ref{eg:Kinyon}.}
\label{fig:Kinyon2}
\end{center}
\end{figure}

\begin{eg}\label{eg:Kinyon}
Consider the five-element set $P=\{1,e,p,q,z\}$, and define the following unary operations on $P$, each written using standard two-line notation for mappings (for now the reader can ignore the colouring of $\boldsymbol{\red p}=z\th_e$):
\begin{equation}\label{eq:Kinyon1}
\th_1 = \trans{1&e&p&q&z\\1&e&p&q&z} \COMMa 
\th_e = \trans{1&e&p&q&z\\e&e&p&q&\boldsymbol{\red p}} \COMMa 
\th_p = \trans{1&e&p&q&z\\p&p&p&p&p} \COMMa 
\th_q = \trans{1&e&p&q&z\\q&q&q&q&q} \ANd
\th_z = \trans{1&e&p&q&z\\z&z&z&z&z} .
\end{equation}
One can easily check (by hand, or with a computer) that these give $P$ the structure of a projection algebra.  Note that $\th_1$ is the identity map, while $\th_p$, $\th_q$ and $\th_z$ are all constant maps.  The multiplication table for $P$, considered as a binary algebra under $u\diamond v = u\th_v$ (cf.~Remark \ref{rem:diamond}), is as follows (and again the reader can ignore the colouring of $\boldsymbol{\red p}=z\diamond e$):
\begin{equation}\label{eq:Kinyon2}
\begin{array}{c|ccccc}
\diamond & 1 & e & p & q & z \\
\hline
1 & 1 & e & p & q & z \\
e & e & e & p & q & z \\
p & p & p & p & q & z \\
q & q & q & p & q & z \\
z & z & \boldsymbol{\red p} & p & q & z 
\end{array}
\end{equation}
Figure \ref{fig:Kinyon} shows the Hasse diagram of the $\leq$ order on $P$, and also the $\leqF$ and $\F$ relations; cf.~\eqref{eq:leqP}, \eqref{eq:leqF} and \eqref{eq:F}.  The diagrams for $\leqF$ and $\F$ omit the loops at each vertex.

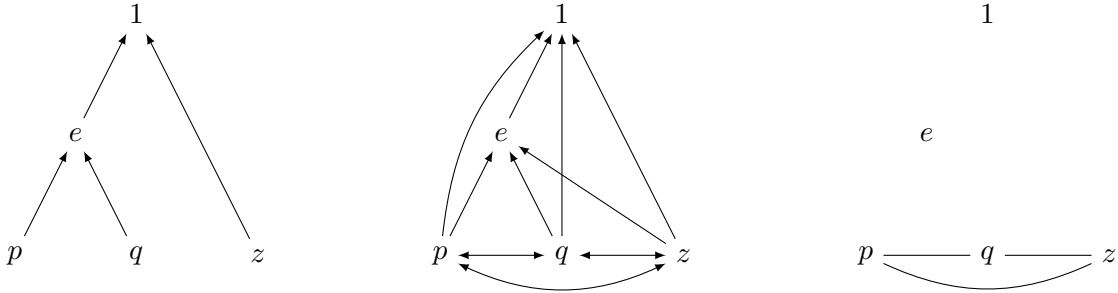
\begin{figure}[ht]
\begin{center}
\begin{tikzpicture}[scale=0.8]
\node (1) at (0,4) {$1$};
\node (e) at (-1,2) {$e$};
\node (p) at (-2,0) {$p$};
\node (q) at (0,0) {$q$};
\node (z) at (2,0) {$z$};
\foreach \x/\y in {e/1,p/e,q/e,z/1} {\draw[-{latex}] (\x)--(\y);}
\begin{scope}[shift = {(7,0)}]
\node (1) at (0,4) {$1$};
\node (e) at (-1,2) {$e$};
\node (p) at (-2,0) {$p$};
\node (q) at (0,0) {$q$};
\node (z) at (2,0) {$z$};
\foreach \x/\y in {e/1,p/e,q/e,z/e,q/1,z/1} {\draw[-{latex}] (\x)--(\y);}
\draw[-{latex}] (p) to [bend left = 20] (1);
\draw[{latex}-{latex}] (p)--(q);
\draw[{latex}-{latex}] (q)--(z);
\draw[{latex}-{latex}] (z) to [bend left = 25] (p);
\end{scope}
\begin{scope}[shift = {(14,0)}]
\node (1) at (0,4) {$1$};
\node (e) at (-1,2) {$e$};
\node (p) at (-2,0) {$p$};
\node (q) at (0,0) {$q$};
\node (z) at (2,0) {$z$};
\draw (p)--(q);
\draw (q)--(z);
\draw (z) to [bend left = 25] (p);
\end{scope}
\end{tikzpicture}
\caption{The relations $\leq$ (left), $\leqF$ (middle) and $\F$ (right) on the projection algebra $P=\{1,e,p,q,z\}$ from Example \ref{eg:Kinyon}.}
\label{fig:Kinyon}
\end{center}
\end{figure}

The construction outlined in Section \ref{subsect:PtoS} (cf.~Proposition \ref{prop:PtoS}) yields the regular $*$-semigroup
\[
\FF_P = \la \ol 1, \ol e, \ol p, \ol q, \ol z\ra \leq \T_P\times\T_P^\opp
\]
with projection algebra (isomorphic to) $P$, where again $\ol u = (\th_u,\th_u)\in\T_P\times\T_P^\opp$ for $u\in P$.  Simple computations with GAP \cite{GAP,Semigroups} show that this semigroup $\FF_P$ is in fact a band of size~$11$.  The element~$\ol1=(\id_P,\id_P)$ is of course the identity of $\FF_P$, while $\ol p$, $\ol q$ and $\ol z$ generate a $3\times3$ rectangular sub-band of $\FF_P$.  Figure \ref{fig:Kinyon2} shows an egg-box diagram for $\FF_P$.  Egg-box diagrams are standard visualisation tools in semigroup theory; see for example \cite{CPbook}.  Each cell of the diagram represents an $\H$-class (all singletons in this case); rows and columns are $\R$- and $\L$-classes, respectively; entire boxes are $\D$-classes; the ordering on ${\D}({}={\J})$-classes is given by $D_a\leq D_b \iff a \in S^1bS^1$.  
\end{eg}

\begin{rem}\label{rem:Kinyon1}
It follows from Lemma \ref{lem:PP'} that any projection algebra automorphism of $P$ (as in Example \ref{eg:Kinyon}) is an order isomorphism of the poset $(P,{\leq})$.  Examining Figure \ref{fig:Kinyon} (left), it is easy to see that there are two of the latter; as well as the identity map, there is the transposition $\phi:P\to P$ swapping $p$ and $q$.  However, this $\phi$ is not a projection algebra morphism, as for example we have
\[
(z\diamond e)\phi = p\phi = q \qquad\text{while}\qquad z\phi\diamond e\phi = z\diamond e = p. 
\]
\end{rem}

\begin{rem}\label{rem:Kinyon2}
One can modify the definition of the projection algebra $P=\{1,e,p,q,z\}$ from Example \ref{eg:Kinyon} by defining $z\th_e$ to be $q$ rather than $p$, i.e.~by changing the red $\boldsymbol{\red p}$ to $\boldsymbol{\red q}$ in \eqref{eq:Kinyon1} and~\eqref{eq:Kinyon2}.  Denote the resulting projection algebra by $P'$.  It is easy to see that the $\leq$ and $\leqF$ (and $\F$) relations on $P$ and $P'$ are exactly the same, even though the operations on $P$ and~$P'$ are (slightly) different.  It follows from this that the relations $\leq$ and $\leqF$ do not uniquely determine the operations on a projection algebra (cf.~Remark \ref{rem:rels}).  On the other hand, $P$ and~$P'$ are clearly isomorphic, so we are left with the question of whether the relations $\leq$ and $\leqF$ determine the algebraic structure of a projection algebra up to isomorphism.
\end{rem}

\section{The chain groupoid of a projection algebra}\label{sect:CP}

We have seen that a regular $*$-semigroup $S$ induces a projection algebra $P=P(S)$ and an ordered groupoid $\G=\G(S)$.  The pair $(P,\G)$ provides \emph{almost} enough information to construct a total invariant of $S$.  The missing ingredient is a certain `link' between $P$ and $\G$, which is provided by a functor $\C\to\G$, where $\C=\C(P)$ is the so-called \emph{chain groupoid} of~$P$, the subject of this section.  We begin in Section \ref{subsect:PP} by defining the \emph{path category} $\P=\P(P)$, and then in Section~\ref{subsect:CP} construct $\C=\P/{\approx}$ as a quotient of $\P$ by a certain congruence $\approx$.  The groupoid $\C$ is a regular $*$-analogue of Nambooripad's groupoid of $E$-chains \cite{Nambooripad1979}.  As we explain in Remark~\ref{rem:Cfunctor}, we can think of $\C$ as a functor $\PA\to\OG$, from projection algebras to ordered groupoids.

\subsection{The path category}\label{subsect:PP}

For the duration of this section we fix a projection algebra $P$, as in Definition \ref{defn:P}, including the operations $\th_p$ ($p\in P$), and the relations $\leq$, $\leqF$ and $\F$ from \eqref{eq:leqP}, \eqref{eq:leqF} and \eqref{eq:F}.
Roughly speaking, we wish to provide an abstract setting in which we can think about `products' of projections $p_1\cdots p_k$, even though~$P$ itself does not have a binary operation.  Guided by the (as-yet unproved) Proposition~\ref{prop:ERSS}, we are (for now) solely interested in such `products' in the case that $p_1\F\cdots\F p_k$.  The \emph{path category}~$\P=\P(P)$ is the first attempt at doing all of this; it represents such a `product' as a tuple $(p_1,\ldots,p_k)$, which is considered as a morphism~$p_1\to p_k$.  %In Section \ref{subsect:CP}, we will define the \emph{chain groupoid} $\C$ (see Definition~\ref{defn:CP}) as a quotient~$\C=\P/{\approx}$ by a certain congruence $\approx$.  This groupoid will be used extensively in 
%The intuition here is that the chain groupoid records more information about such products than the factors alone, but only such information that is present in \emph{any} regular $*$-semigroup with projection algebra $P$.  For example,~$(p,p)$ and~$(p)$ should always represent the same `product', since projections are idempotents; so too should~$(p,q,p)$ and~$(p)$ when $p\F q$ (cf.~\eqref{eq:leqFF}).  These pairs of paths (and one other family of pairs) are taken as generators for the congruence $\approx$.  But before we get ahead of ourselves, here is the main definition for the current section:

\begin{defn}\label{defn:PP}
A \emph{$P$-path} in a projection algebra $P$ (cf.~Definition \ref{defn:P}) is a tuple
\[
\p = (p_1,p_2,\ldots,p_k)\in P^k \qquad\text{for some $k\geq1$, such that $p_1\F p_2\F\cdots\F p_k$.}
\]
(Since $\F$ need not be transitive, this does not imply that the $p_i$ are \emph{all} $\F$-related.)  We say that $\p$ is a $P$-path from $p_1$ to $p_k$, and write $\bd(\p)=p_1$ and $\br(\p)=p_k$.  We identify each $p\in P$ with the path $(p)$ of length $1$.  

The \emph{path category} of $P$ is the $*$-category $\P=\P(P)$ of all $P$-paths, with object set $v\P=P$, under the following operations:
\bit
\item For $\p=(p_1,\ldots,p_k)$ and $\q=(q_1,\ldots,q_l)$ with $p_k=q_1$, we define
\[
\p\circ\q = (p_1,\ldots,p_{k-1},p_k=q_1,q_2,\ldots,q_l).
\]
\item For $\p=(p_1,\ldots,p_k)\in\P$, we define
\[
\p^\rev=(p_k,\ldots,p_1),
\]
the reverse of $\p$.  (It is convenient to write $\p^\rev$ instead of $\p^*$, and it is clear that conditions \ref{I1}--\ref{I3} from Definition \ref{defn:*cat} all hold.)
\eit
\end{defn}

So a morphism set $\P(p,q)$ consists of all $P$-paths from $p$ to $q$.  The identities are the paths of the form $p\equiv(p)$.  Although $\P$ is a $*$-category, it is not a groupoid, as when~$\p$ has length $k\geq2$, $\p\circ\p^\rev$ has length $2k-1>k$.  However, and as noted above, a very important role will be played by a certain groupoid quotient of $\P$; see Definition \ref{defn:CP}.

\begin{rem}\label{rem:Pfree}
It is also worth noting that $\P=\P(P)$
%, considered as an unordered $*$-category, 
is (isomorphic to) the free $*$-category over the relation~$\F$ in the following sense.  We define~${\Ga=\Ga_P}$ to be the digraph with vertex set~$P$, and an edge $x_{pq}:p\to q$ for each $(p,q)\in{\F}$.  As in \cite[p.~49]{MacLane1998}, the \emph{free category} $\CC=\CC(\Ga)$ has object set $v\CC=P$, and its morphisms are the paths in $\Ga$, including the empty path at each vertex $p\in P$ (which are identified with $p$, as usual).  Non-empty paths can be thought of as words $x_{p_1p_2}x_{p_2p_3}\cdots x_{p_{k-1}p_k}$ (note the matching subscripts between successive letters), and composition of paths is given by concatenation when the endpoints match.  The involution of~$\CC$ is given by reversal of paths; this is well-defined since~$\F$ is symmetric.  It is then clear that
\[
x_{p_1p_2}x_{p_2p_3}\cdots x_{p_{k-1}p_k}\to(p_1,p_2,\ldots,p_k)
\]
defines a $*$-isomorphism $\CC\to\P$.
\end{rem}

We now wish to show that $\P$ is an \emph{ordered} $*$-category, using Lemma \ref{lem:C}.  To apply the lemma, we need a partial order on the object set $v\P=P$, and a collection of (left) restrictions~${}_q\corest\p$.  We already have the order $\leq$ on $P$ given in \eqref{eq:leqP}.  To define the restrictions, consider a $P$-path $\p=(p_1,\ldots,p_k)$, and suppose $q\in P$ is such that $q\leq \bd(\p) = p_1$.  We define a tuple
\begin{equation}\label{eq:rest}
{}_q \corest \p = (q_1,\ldots,q_k) \WHERE q_1=q \ANd q_i = q_{i-1}\th_{p_i} \text{ for $2\leq i\leq k$.}
\end{equation}
Note that $q_i = q\th_{p_2}\cdots\th_{p_i}$ for all $1\leq i\leq k$.  In fact, since $q\leq p_1$ we have $q=q\th_{p_1}$, so that
\begin{equation}\label{eq:qi}
q_i = q\th_{p_2}\cdots\th_{p_i} = q\th_{p_1}\cdots\th_{p_i} \qquad\text{for all $1\leq i\leq k$.}
\end{equation}
It follows immediately that $q_i\leq p_i$ for all $i$.  We now prove a sequence of lemmas that will verify axioms \ref{O1'}--\ref{O5'} from Lemma \ref{lem:C}.

\begin{lemma}\label{lem:q|p}
If $\p\in\P$, and if $q\leq\bd(\p)$, then ${}_q\corest\p\in\P$.  Moreover, $\bd({}_q\corest\p)=q$ and ${\br({}_q\corest\p)\leq\br(\p)}$.
\end{lemma}

\pf
Write $\p=(p_1,\ldots,p_k)$ and ${}_q\corest\p=(q_1,\ldots,q_k)$, as above.  Since $\p\in\P$, we have $p_i\F p_{i+1}$ for all $1\leq i<k$, and we must show that $q_i\F q_{i+1}$ for all such $i$, i.e.~that
\[
q_i\th_{q_{i+1}}=q_{i+1} \AND q_{i+1}\th_{q_i} = q_i.
\]
For the first, we have
\[
q_i \th_{q_{i+1}} = q_i \th_{q_i\th_{p_{i+1}}} =_4 q_i\th_{p_{i+1}}\th_{q_i}\th_{p_{i+1}} =_3 p_{i+1}\th_{q_i}\th_{p_{i+1}} =_3 q_i\th_{p_{i+1}} = q_{i+1}.
\]
For the second, it follows from \eqref{eq:qi} that $q_i=_1q_i\th_{p_i}$.  Combining this with $p_i=p_{i+1}\th_{p_i}$ (as~${p_i\F p_{i+1}}$), we have
\[
q_{i+1} \th_{q_i} = q_i\th_{p_{i+1}}\th_{q_i} =_3 p_{i+1}\th_{q_i} = p_{i+1}\th_{q_i\th_{p_i}} =_4 p_{i+1}\th_{p_i}\th_{q_i}\th_{p_i} = p_i\th_{q_i}\th_{p_i} =_3 q_i\th_{p_i} = q_i.
\]
This all shows that ${}_q\corest\p\in\P$.

By definition, we have $\bd({}_q\corest\p)=q_1=q$, and $\br({}_q\corest\p)=q_k\leq p_k=\br(\p)$, where the latter follows from the fact, noted above, that $q_i\leq p_i$ for all $i$.
\epf

\begin{lemma}\label{lem:qp=pr}
If $\p\in\P$ and $q\leq\bd(\p)$, and if $r=\br({}_q\corest\p)$, then $({}_q\corest\p)^\rev={}_r\corest\p^\rev$.
\end{lemma}

\pf
Write $\p=(p_1,\ldots,p_k)$ and ${}_q\corest\p=(q_1,\ldots,q_k)$, so that $r=q_k$.  We also of course have $\p^\rev=(p_k,\ldots,p_1)$ and $({}_q\corest\p)^\rev=(q_k,\ldots,q_1)$.
To make the subscripts match up, it is convenient to write ${}_r\corest\p^\rev=(r_k,\ldots,r_1)$.  So
\begin{equation}\label{eq:qiri}
q_i = q\th_{p_1}\cdots\th_{p_i} \AND r_i = r\th_{p_k}\cdots\th_{p_i} \qquad\text{for all $1\leq i\leq k$.}
\end{equation}
We show by descending induction that $q_i=r_i$ for all $i$.  The $i=k$ case is trivial, as $r_k=r=q_k$.  For $1\leq i<k$, we have
\begin{align*}
r_i &= r_{i+1}\th_{p_i} &&\text{by \eqref{eq:qiri}}\\
&= q_{i+1}\th_{p_i} &&\text{by induction}\\
&= q\th_{p_1}\cdots\th_{p_i}\th_{p_{i+1}}\th_{p_i} &&\text{by \eqref{eq:qiri}}\\
&= q\th_{p_1}\cdots\th_{p_i} &&\text{by Lemma \ref{lem:pqp}\ref{pqp2}}\\
&= q_i &&\text{by \eqref{eq:qiri}.}  \qedhere
\end{align*}
\epf

\begin{lemma}\label{lem:reflexive}
If $\p\in\P$ and $q=\bd(\p)$, then ${}_q\corest\p=\p$.
\end{lemma}

\pf
Write $\p=(p_1,\ldots,p_k)$, so that $q=p_1$.  Then ${}_q\corest\p=(q_1,\ldots,q_k)$, where each $q_i=q\th_{p_1}\cdots\th_{p_i}$, and we must show that $p_i=q_i$ for all $i$.  This is clear for $i=1$.  For $i\geq2$, we have
\[
q_i = q_{i-1}\th_{p_i} = p_{i-1}\th_{p_i} = p_i,
\]
by definition, induction, and the fact that $p_i\F p_{i-1}$.
\epf

\begin{lemma}\label{lem:transitive}
If $\p\in\P$ and $r\leq q\leq\bd(\p)$, then ${}_r\corest{}_q\corest\p={}_r\corest\p$.
\end{lemma}

\pf
Write
\[
\p=(p_1,\ldots,p_k) \COMMA {}_q\corest\p=(q_1,\ldots,q_k) \COMMA {}_r\corest\p=(r_1,\ldots,r_k) \AND {}_r\corest{}_q\corest\p=(s_1,\ldots,s_k).
\]
So
\begin{equation}\label{eq:pqrsi}
q_i=q\th_{p_1}\cdots\th_{p_i} \COMMA r_i=r\th_{p_1}\cdots\th_{p_i} \AND s_i=r\th_{q_1}\cdots\th_{q_i} \qquad\text{for all $1\leq i\leq k$,}
\end{equation}
and we must show that $r_i=s_i$ for all $i$.  For $i=1$ we have $r_1=r=s_1$.  For $i\geq2$,
\begin{align*}
s_i &= s_{i-1}\th_{q_i}   &&\text{by \eqref{eq:pqrsi}}\\
&= (s_{i-1}\th_{q_{i-1}})\th_{q_{i-1}\th_{p_i}}  &&\text{by \eqref{eq:pqrsi}, noting that $s_{i-1}\leq q_{i-1}$}\\
&= s_{i-1}\th_{q_{i-1}}\th_{p_i}\th_{q_{i-1}}\th_{p_i}  &&\text{by \ref{P4}}\\
&= s_{i-1}\th_{q_{i-1}}\th_{p_i}  &&\text{by \ref{P5}}\\
&= s_{i-1}\th_{p_i}  &&\text{by \eqref{eq:pqrsi}}\\
&= r_{i-1}\th_{p_i}  &&\text{by induction}\\
&= r_i  &&\text{by \eqref{eq:pqrsi}.}  \qedhere
\end{align*}
\epf

\begin{lemma}\label{lem:O2}
If $\p,\q\in\P$ with $\br(\p) = \bd(\q)$, and if $r\leq\bd(\p)$ and $s=\br({}_r\corest\p)$, then
\[
{}_r\corest(\p\circ\q) = {}_r\corest\p \circ {}_{s}\corest\q.
\]
\end{lemma}

\pf
Write $\p=(p_1,\ldots,p_k)$ and $\q=(q_1,\ldots,q_l)$, noting that $p_k=q_1$.  Then
\begin{align*}
{}_r\corest\p&=(r,r\th_{p_2},r\th_{p_2}\th_{p_3},\ldots,r\th_{p_2}\cdots\th_{p_k}), \qquad\text{so}\qquad s=r\th_{p_2}\cdots\th_{p_k}.
\intertext{It follows that}
{}_r\corest(\p\circ\q) &= {}_r\corest(p_1,\ldots,p_k,q_2,\ldots,q_l) \\
&= (r,r\th_{p_2},r\th_{p_2}\th_{p_3},\ldots,r\th_{p_2}\cdots\th_{p_k}=s,s\th_{q_2},s\th_{q_2}\th_{q_3},\ldots,s\th_{q_2}\cdots\th_{q_l}) \\
&= (r,r\th_{p_2},r\th_{p_2}\th_{p_3},\ldots,r\th_{p_2}\cdots\th_{p_k})\circ(s,s\th_{q_2},s\th_{q_2}\th_{q_3},\ldots,s\th_{q_2}\cdots\th_{q_l})  = {}_r\corest\p \circ {}_{s}\corest\q. \qedhere
\end{align*}
\epf

\begin{prop}
For any projection algebra $P$, the path category $\P=\P(P)$ is an ordered $*$-category, with ordering given by
\[
\p\leq\q \IFF \p={}_r\corest\q  \qquad \text{for some $r\leq\bd(\q)$.}
\]
\end{prop}

\pf
This follows from an application of Lemma \ref{lem:C}.  The ordering on $v\P=P$ is given in~\eqref{eq:leqP}.  Properties \ref{O1'}--\ref{O5'} were established in Lemmas \ref{lem:q|p}--\ref{lem:O2}.
\epf

As usual (cf.~Remark \ref{rem:dual}), we can use the involution to define a right-handed restriction:
\begin{equation}\label{eq:revrev}
\p\rest_q = ( {}_q\corest\p^\rev)^\rev \qquad\text{for $\p\in\P$ and $q\leq\br(\p)$.}
\end{equation}
Explicitly, if $\p=(p_1,\ldots,p_k)$ and $q\leq\br(\p)=p_k$, then
\begin{equation}\label{eq:corest}
\p\rest_q = (q_1,\ldots,q_k) \WHERE q_i = q\th_{p_k}\cdots\th_{p_i} \qquad\text{for all $1\leq i\leq k$.}
\end{equation}

\subsection{The chain groupoid}\label{subsect:CP}

We are now almost ready to define the \emph{chain groupoid} $\C=\C(P)$ associated to a projection algebra $P$.  This groupoid is defined below as a certain quotient $\C=\P/{\approx}$ of the path category $\P=\P(P)$ from Definition \ref{defn:PP}.  The congruence $\approx$ is defined by specifying a generating set:

\begin{defn}\label{defn:approx}
Given a projection algebra $P$ (cf.~Definition \ref{defn:P}), let $\Om=\Om(P)$ be the set of all pairs $(\s,\t)\in\P\times\P$ of the following two forms:
\begin{enumerate}[label=\textup{\textsf{($\mathsf{\Om}$\arabic*)}},leftmargin=10mm]
\item \label{Om1} $\s=(p,p)$ and $\t=(p)\equiv p$, for some $p\in P$,
\item \label{Om2} $\s=(p,q,p)$ and $\t=(p)\equiv p$, for some $(p,q)\in {\F}$.
\end{enumerate}
We define ${\approx}=\Om^\sharp$ to be the congruence on $\P$ generated by $\Om$.  
\end{defn}

Since $\bd(\s)=\bd(\t)$ and $\br(\s)=\br(\t)$ for all $(\s,\t)\in\Om$, it follows that $\approx$ is a $v$-congruence.

\begin{lemma}\label{lem:OmP}
$\approx$ is an ordered $*$-congruence.
\end{lemma}

\pf
We must show that $\approx$ satisfies conditions~\ref{C4} and~\ref{C5} from Definition \ref{defn:cong}.  By Lemma~\ref{lem:Om}, it suffices to show that the generating set $\Om$ satisfies \eqref{eq:Om*} and \eqref{eq:Om}.  

For \eqref{eq:Om*}, we must show that $(\s,\t)\in\Om \implies \s^\rev\approx\t^\rev$ for all $(\s,\t)\in\Om$.  But this is immediate, as $\s^\rev=\s$ and $\t^\rev=\t$ for all such $(\s,\t)$.  

For \eqref{eq:Om}, let $(\s,\t)\in\Om$.  We must show that
\[
{}_r\corest\s \approx {}_r\corest\t \qquad\text{for all $r\leq\bd(\s)$.}
\]
We do this separately for the two forms the pair $(\s,\t)$ can take.  

\pfitem{\ref{Om1}}  Here we have ${}_r\corest\s=(r,r)$ and ${}_r\corest\t=(r)$.  In particular, $({}_r\corest\s,{}_r\corest\t)\in\Om$, so certainly ${}_r\corest\s \approx {}_r\corest\t$.

\pfitem{\ref{Om2}}  To calculate ${}_r\corest\s = {}_r\corest(p,q,p)$, we first note that $r\th_q\th_p = r\th_p\th_q\th_p = r\th_p = r$, where we used $r\leq p$ in the first and third steps, and Lemma \ref{lem:pqp}\ref{pqp2} in the second.  We then have
\[
{}_r\corest\s = (r,r\th_q,r\th_q\th_p) = (r,r\th_q,r) \qquad\text{and again}\qquad {}_r\corest\t=(r).
\]
So once again we have $({}_r\corest\s,{}_r\corest\t)\in\Om$.
\epf

It follows quickly from iterating \ref{Om2} that
\[
\p\circ\p^\rev\approx\bd(\p) \qquad\text{for all $\p\in\P$.}
\]
Combining this with Lemmas \ref{lem:approx} and \ref{lem:OmP}, it follows that the quotient $\P/{\approx}$ is an ordered groupoid.

\begin{defn}\label{defn:CP}
The \emph{chain groupoid} of a projection algebra $P$ (cf.~Definition \ref{defn:P}) is the quotient
\[
\C = \C(P) = \P/{\approx},
\]
where $\P=\P(P)$ is the path category of $P$ (cf.~Definition \ref{defn:PP}), and where $\approx$ is the congruence given in Definition \ref{defn:approx}.  
\bit
\item The elements of $\C$, which are $\approx$-classes of $P$-paths, are called \emph{$P$-chains}.  For $\p\in\P$, we write $[\p]\in\C$ for the $\approx$-class of $\p$.  If $\p=(p_1,\ldots,p_k)$, then we write $[\p] = [p_1,\ldots,p_k]$.  We then have $\bd[\p]=\bd(\p)=p_1$ and $\br[\p]=\br(\p)=p_k$.
\item For $\p=(p_1,\ldots,p_k)$ and $\q=(q_1,\ldots,q_l)$ with $p_k=q_1$, we have
\[
[\p]\circ[\q] = [\p\circ\q] = [p_1,\ldots,p_{k-1},p_k=q_1,q_2,\ldots,q_l].
\]
\item For $\p=(p_1,\ldots,p_k)\in\P$, we have
\[
[\p]^{-1}=[\p^\rev] = [p_k,\ldots,p_1].
\]
\item The order in $\C$ is given by
\[
\c\leq\d \IFF \p\leq\q \qquad\text{for some $\p\in\c$ and $\q\in\d$.}
\]
\eit
\end{defn}

%\begin{rem}
%As in Remark \ref{rem:Pfree}, it is not hard to see that $\C$ (considered as an unordered category) is (isomorphic to) the free groupoid over the relation ${\F}\sm\De_P$, which is in turn isomorphic to the fundamental groupoid of the graph obtained from $\Ga_P$ (as defined in Remark \ref{rem:Pfree}).  
%
%%It is also worth commenting on the definition of the congruence $\approx$.  The chain groupoid $\C=\C(P)$ will be used in Section \ref{sect:G} to provide an environment in which to `interpret' products of projections in certain abstract groupoids $\G$ with object set $P$.  Since distinct projections/objects have distinct co/domains, they cannot be composed in such a groupoid.  But a $P$-chain $[p_1,\ldots,p_k]\in\C$ is meant to be thought of as a `product' $p_1\cdots p_k$, and the `interpretation' mentioned above is a functor $\C\to\G$.  Thus, $\approx$ is meant to equate such `products' when they `should be equal' in any regular $*$-semigroup with projection algebra $P$.  
%%%
%%In this way, it is clear why $[p,p]$ and $[p]$ `should' be equal for any $p\in P$; cf.~\ref{Om1}.  Similar comments apply to $[p,q,p]$ and $[p]$ when $p\F q$; cf.~\ref{Om2} and \eqref{eq:leqFF}.  
%\end{rem}

The next result shows that any projection algebra morphism naturally induces a functor between the corresponding chain groupoids.

\begin{prop}\label{prop:CC'}
If $\phi:P\to P'$ is a projection algebra morphism, then there is a well-defined ordered groupoid functor
\[
\Phi = \C(\phi) :\C(P)\to\C(P') \GIVENBY [p_1,\ldots,p_k]\Phi = [p_1\phi,\ldots,p_k\phi].
\]
\end{prop}

\pf
During the proof we write $\C=\C(P)$ and $\C'=\C(P')$, and similarly for $\P$ and $\P'$.  By Lemma \ref{lem:PP'}, we have a well-defined functor
\[
\varphi:\P\to\C' \GIVENBY (p_1,\ldots,p_k)\varphi=[p_1\phi,\ldots,p_k\phi].
\]
We first note that ${\approx}\sub\ker(\varphi)$.  Indeed, for any pair $(\s,\t)\in\Om$ we have $(\s\varphi,\t\varphi)\in\Om'$ (the generating set for the corresponding congruence on $\P'=\P(P')$).  It follows that we have a well-defined functor
\[
\Phi:\C=\P/{\approx} \to \C' \GIVENBY [\p]\Phi = \p\varphi,
\]
and this is of course the map in the statement of the proposition.  It is easy to check that $\phi$ being a morphism implies $({}_q\corest\c)\Phi={}_{q\phi}\corest(\c\Phi)$ for all $\c\in\C$ and $q\leq\bd(\c)$, so that $\Phi$ is ordered.
\epf

\begin{rem}\label{rem:Cfunctor}
As in Propositions \ref{prop:calGfunctor} and \ref{prop:Pfunctor}, it follows from Proposition \ref{prop:CC'} that we can think of $\C$ as a functor $\PA\to\OG$.  (Verification that $\C(\phi\circ\phi')=\C(\phi)\circ\C(\phi')$ for composable morphisms $\phi,\phi'$ in $\PA$ is routine.)

Interestingly, the functor $\C$ is not injective (on objects).  For example, let $P=\{1,e,p,q,z\}$ be the projection algebra from Example \ref{eg:Kinyon}, and let $P'$ be the closely related projection algebra from Remark \ref{rem:Kinyon2}, defined on the same underlying set.  Since $P$ and $P'$ have exactly the same~$\F$ relation, $\C(P)$ and $\C(P')$ are identical as groupoids.  
%(To see this, we note that the path categories~$\P(P)$ and~$\P(P')$ contain exactly the same paths, and that the $\approx$ congruences on $\P(P)$ and $\P(P')$ are generated by the same sets of pairs.)  
In fact, they are identical as \emph{ordered} groupoids.  Indeed, $P$ and $P'$ have the same $\leq$ order, as pictured in Figure \ref{fig:Kinyon} (left), and it follows that the only non-trivial restrictions in the two groupoids are ${}_s\corest [1]=[1]\rest_s=[s]$ ($s=e,p,q,z$) and ${}_t\corest [e]=[e]\rest_t=[t]$ ($t=p,q$).  

The functor $\C$ is also not surjective on morphisms, meaning that the image of $\PA$ under $\C$ is not a \emph{full} category of $\OG$.  Indeed, consider again the projection algebra $P$ from Example~\ref{eg:Kinyon}, and let $\phi$ be the order-isomorphism of $(P,{\leq})$ swapping $p\leftrightarrow q$, and fixing $1,e,z$ pointwise.  Then it is easy to see that there is a well-defined ordered functor $\Phi:\C(P)\to\C(P)$ given by $[p_1,\ldots,p_k]\Phi=[p_1\phi,\ldots,p_k\phi]$.  (That this is a functor follows from the form of chains in $\C(P)$, and that it is ordered again follows from considering the non-trivial restrictions in $\C(P)$.)  However, this $\Phi$ is not in the image of $\C$, as $\phi=v\Phi$ is not a projection algebra morphism, as explained in Remark \ref{rem:Kinyon2}.
\end{rem}

\begin{rem}\label{rem:Cfree}
As in Remark \ref{rem:Pfree}, we can think of the unordered reduct of a chain groupoid~$\C(P)$ as the free groupoid over the relation ${\F}\sm\De_P$.  This time we take $\Ga'_P$ to be the digraph with vertex set $P$, and an edge $x_{pq}:p\to q$ for each $(p,q)\in{\F}\sm\De_P$.  The free groupoid over ${\F}\sm\De_P$ consists of the paths in $\Ga'_P$ (as in Remark \ref{rem:Pfree}), modulo the relations $x_{pq}x_{qp}=p$ (where the latter is identified with the empty path at $p$).  Free categories and groupoids over graphs need not be ordered in general, but the projection algebra axioms ensure that those of the form~$\P(P)$ and~$\C(P)$ are.
\end{rem}

\section{Chained projection groupoids}\label{sect:G}

In Sections \ref{sect:RSS} and \ref{sect:P} we constructed functors
\[
\G:\RSS\to\OG \AND P:\RSS\to\PA.
\]
At the object level, these functors map a regular $*$-semigroup $S$ to the ordered groupoid~$\G(S)$ and projection algebra $P(S)$ from Definitions \ref{defn:GS} and \ref{defn:PS}.  It follows from the examples considered in Section \ref{subsect:Seg} that neither $\G$ nor $P$ is injective (on objects).
%
%Recall from Section \ref{sect:RSS} (see Definition \ref{defn:GS} and Proposition \ref{prop:GS}) that to any regular $*$-semigroup~$S$ we can associate an ordered groupoid $\G=\G(S)$.  The object set of this groupoid is $v\G=P=P(S)$, the projection algebra of~$S$, and the (partial) composition in $\G$ is a restriction of the (total) product in $S$.  We showed in Proposition \ref{prop:calGfunctor} that $\G$ can be thought of as a functor $\RSS\to\OG$ from the category of regular $*$-semigroups to the category of ordered groupoids.  In Section \ref{subsect:Seg} we gave some examples to show that $\G$ is not injective, in the sense that distinct regular $*$-semigroups $S$ and $S'$ (on the same underlying set) can give exactly the same ordered groupoid $\G(S)=\G(S')$.
%
The purpose of the current section is to provide an appropriate enrichment of ordered groupoids in order to obtain an injective functor from regular $*$-semigroups to a suitable category of groupoids.  We do this in stages, beginning with \emph{weak projection groupoids} (Section~\ref{subsect:WPG}), \emph{projection groupoids} (Section \ref{subsect:PG}) and \emph{weak chained projection groupoids} (Section \ref{subsect:ChPG}), eventually culminating in the category $\CPG$ of \emph{chained projection groupoids} (Section \ref{subsect:CoPG}).  The main result of the current section is Theorem~\ref{thm:Gfunctor}, which establishes the existence of a functor $\bG:\RSS\to\CPG$.  In Section \ref{sect:GtoS} we construct a functor $\bS:\CPG\to\RSS$ in the opposite direction, and in Section~\ref{sect:iso} we show that~$\bG$ and~$\bS$ are mutually inverse isomorphisms between the categories $\RSS$ and~$\CPG$.

\subsection{Weak projection groupoids}\label{subsect:WPG}

In all that follows, we are concerned with ordered groupoids whose object sets form projection algebras.  The most basic situation is as follows:

\begin{defn}\label{defn:WPG}
A \emph{weak projection groupoid} is a pair $(P,\G)$, consisting of an ordered groupoid~$\G$ and a projection algebra $P=v\G$, for which
%\begin{enumerate}[label=\textup{\textsf{(G\arabic*)}},leftmargin=10mm]%\addtocounter{enumi}{0}
%\item \label{G1} 
the restriction to $P$ of the order on $\G$ is the relation~$\leq$ from \eqref{eq:leqP}.
%\end{enumerate}
\end{defn}

%\begin{rem}
Given a weak projection groupoid $(P,\G)$, the ordering on $\G$ is inherited from that of $P$ in the sense described in Lemma~\ref{lem:C}.  That is, for each $a\in\G$ and each $p\leq\bd(a)$ we have a left restriction~${{}_p\corest a\in\G}$, with respect to which conditions \ref{O1'}--\ref{O5'} hold.  We also have the right restrictions $a\rest_q=({}_q\corest a^{-1})^{-1}$, defined for $q\leq\br(a)$, and these satisfy the duals of \ref{O1'}--\ref{O5'}.  We typically use these conditions without explicit reference, and also \ref{O6'} and its dual, which follow from the others.
%\end{rem}

In a weak projection groupoid $(P,\G)$, the only `link' between the structures of the groupoid~$\G$ and the projection algebra $P=v\G$ is at the order-theoretic level.  In particular, the link does not incorporate the algebraic structure of $P$.  Indeed, while the order on $P$ is defined \emph{from} the $\th$ operations (cf.~\eqref{eq:leqP}), the order does not \emph{determine} the~$\th$ operations.  (We saw in Examples~\ref{eg:AG},~\ref{eg:Rees} and \ref{eg:Kinyon} that the same order can arise from different projection algebra structures defined on a common set.)  In what follows, we will be particularly interested in groupoid/projection algebra pairs with a stronger, \emph{algebraic} link.  Such a link manifests itself in a number of equivalent properties listed in Proposition \ref{prop:G1} below, and stated in terms of certain additional maps $P\to P$ whose definition we now give.

Consider a weak projection groupoid $(P,\G)$, and a morphism $a\in\G$.  As in~\eqref{eq:vta}, we have the map
\[
\vt_a : \bd(a)^\da\to\br(a)^\da \GIVENBY p\vt_a = \br({}_p\corest a).
\]
Since $\bd(a)\in P$, we also have the unary operation $\th_{\bd(a)}$, and by \eqref{eq:imthp} its image is $\bd(a)^\da$.  It follows that we can compose $\th_{\bd(a)}$ with $\vt_a$, and we denote this composition by
\begin{equation}\label{eq:Tha}
\Th_a = \th_{\bd(a)}\vt_a : P\to\br(a)^\da.
\end{equation}
As a special case, note that by Lemma \ref{lem:vt}\ref{vt2}, we have
\begin{equation}\label{eq:Thp}
\Th_p = \th_{\bd(p)}\vt_p = \th_p\id_{p^\da} = \th_p \qquad\text{for any projection $p\in P$.}
\end{equation}
We begin by recording some basic properties of the $\Th$ maps.

\begin{lemma}\label{lem:G2}
If $(P,\G)$ is a weak projection groupoid, then for any $a\in\G$ we have:
\ben
\item \label{G21} $\vt_a\th_{\br(a)} = \vt_a$,
\item \label{G22} $\Th_a\th_{\br(a)} = \Th_a = \th_{\bd(a)}\Th_a$,
\item \label{G23} $\Th_a\vt_{a^{-1}} = \th_{\bd(a)}$,
%\item \label{G24} $\Th_a = \th_{\bd(a)}\vt_b$ for any $a\leq b$,
\item \label{G25} $\Th_{{}_p\corest a} = \th_p\Th_a$ for any $p\leq\bd(a)$.
\een
\end{lemma}

\pf
\firstpfitem{\ref{G21}}  Since $\im(\vt_a)=\br(a)^\da$, this follows immediately from the fact that each operation $\th_p$ fixes each element of $\im(\th_p)=p^\da$.

\pfitem{\ref{G22}}  Since $\Th_a = \th_{\bd(a)}\vt_a$, we obtain $\Th_a\th_{\br(a)} = \Th_a$ from part \ref{G21}, and $\th_{\bd(a)}\Th_a = \Th_a$ from \ref{P2}.

\pfitem{\ref{G23}}  Since $\im(\th_{\bd(a)})=\bd(a)^\da$, it follows from Lemma \ref{lem:vtavta*} that
\[
\Th_a\vt_{a^{-1}} = \th_{\bd(a)}\vt_a\vt_{a^{-1}} = \th_{\bd(a)}\id_{\bd(a)^\da} = \th_{\bd(a)}.
\]

%\pfitem{\ref{G24}}  Since $a\leq b$, we have $a={}_p\corest b$, where $p=\bd(a)$, and Lemma \ref{lem:vt}\ref{vt3} gives $\vt_a=\vt_b|_{p^\da}$.  It follows that
%\[
%\Th_a = \th_{\bd(a)}\vt_a = \th_p \vt_b|_{p^\da} = \th_p\vt_b.
%\]
%
%%Now let $s\in P$ be arbitrary; we must show that $s\Th_a = s\th_p\vt_b$.  For this, we write $t=s\th_p$, and use \ref{O4'} to calculate
%%\[
%%s\Th_a = s\th_p\vt_a = t\vt_a = \br({}_t\corest a) = \br({}_t\corest {}_p\corest b) = \br({}_t\corest b) = t\vt_b = s\th_p\vt_b.
%%\]

\pfitem{\ref{G25}}  We have
\begin{align*}
\Th_{{}_p\corest a} = \th_{\bd({}_p\corest a)}\vt_{{}_p\corest a} &= \th_p\vt_a|_{p^\da} &&\text{by Lemma \ref{lem:vt}\ref{vt3}}\\
&=\th_p\vt_a &&\text{as $p^\da=\im(\th_p)$ by \eqref{eq:imthp}}\\
&=\th_p\th_{\bd(a)}\vt_a &&\text{by Lemma \ref{lem:thpthq}, as $p\leq\bd(a)$}\\
&= \th_p\Th_a. &&\qedhere
\end{align*}
%
%\begin{align*}
%\Th_{{}_p\corest a} = \th_{\bd({}_p\corest a)}\vt_a &= \th_p\vt_a &&\text{by part \ref{G24}, as ${}_p\corest a\leq a$}\\
%&= \th_p\th_{\bd(a)}\vt_a &&\text{by Lemma \ref{lem:thpthq}, as $p\leq\bd(a)$}\\
%&= \th_p\Th_a. &&\qedhere
%\end{align*}
\epf

Another important property is that when $a$ and~$b$ are composable in $\G$, the map~$\Th_{a\circ b}$ is the composition of~$\Th_a$ and $\Th_b$.

\begin{lemma}\label{lem:Thab}
If $(P,\G)$ is a weak projection groupoid, and if $a,b\in\G$ are such that $\br(a)=\bd(b)$, then ${\Th_{a\circ b}=\Th_a\Th_b}$.
\end{lemma}

\pf
Write $p=\bd(a)$ and $q=\br(a)=\bd(b)$.  Then by Lemmas \ref{lem:vt}\ref{vt1} and \ref{lem:G2}\ref{G21} we have
\[
\Th_a\Th_b = \th_p\vt_a\th_q\vt_b = \th_p\vt_a\th_{\br(a)}\vt_b = \th_p\vt_a\vt_b = \th_p\vt_{a\circ b} = \Th_{a\circ b}.  \qedhere
\]
\epf

Part \ref{G25} of Lemma \ref{lem:G2} above shows how the $\Th$ maps interact with left restrictions, but the lemma did not include a corresponding statement concerning right restrictions.  Of course if $q\leq\br(a)$, then $a\rest_q = {}_p\corest a$ where $p=q\vt_{a^{-1}}$ (cf.~\eqref{eq:vta*}), and then
\begin{equation}\label{eq:Thaq}
\Th_{a\rest_q} = \Th_{{}_p\corest a} = \th_p\Th_a = \th_{q\vt_{a^{-1}}}\Th_a.
\end{equation}
However, the groupoids we will be concerned with satisfy the neater identity $\Th_{a\rest_q}=\Th_a\th_q$.  This is in fact equivalent to a number of additional properties linking the groupoid structure of~$\G$ to the projection algebra structure of $P=v\G$:

\begin{prop}\label{prop:G1}
If $(P,\G)$ is a weak projection groupoid, then the following are equivalent:
\begin{enumerate}[label=\textup{\textsf{(G1\alph*)}},leftmargin=12mm]
\item \label{G2a} $\th_{p\vt_a} = \Th_{a^{-1}}\th_p\Th_a$ for all $a\in\G$ and $p\leq\bd(a)$,
\item \label{G2b} $\th_{p\Th_a} = \Th_{a^{-1}}\th_p\Th_a$ for all $a\in\G$ and $p\in P$,
\item \label{G2c} $\Th_{a\rest_q}=\Th_a\th_q$ for all $a\in\G$ and $q\leq\br(a)$,
\item \label{G2d} $\vt_a$ is a projection algebra morphism (and hence isomorphism) $\bd(a)^\da\to\br(a)^\da$ for all $a\in\G$.
\end{enumerate}
\end{prop}

\pf
For the duration of the proof we fix $a\in\G$, and we write $s=\bd(a)$ and $t=\br(a)$.  On a number of occasions we will use the fact that
\begin{equation}\label{eq:Tha*Tha}
\Th_{a^{-1}}\Th_a = \Th_{a^{-1}\circ a} = \Th_{\br(a)} = \Th_t = \th_t.
\end{equation}
In the above calculation we used Lemma \ref{lem:Thab} and \eqref{eq:Thp}.  In what follows, we show that \ref{G2a} implies each of \ref{G2b}--\ref{G2d}, and conversely that any of \ref{G2b}--\ref{G2d} implies \ref{G2a}.

\pfitem{\ref{G2a}$\implies$\ref{G2b}--\ref{G2d}}  Suppose first that \ref{G2a} holds.  To verify \ref{G2b}, let $p\in P$ be arbitrary.  Then
\begin{align*}
\th_{p\Th_a} &= \th_{(p\th_s)\vt_a} &&\text{by definition}\\
&= \Th_{a^{-1}}\th_{p\th_s}\Th_a &&\text{by \ref{G2a}, as $p\th_s\leq s=\bd(a)$}\\
&= \Th_{a^{-1}}\th_s\th_p\th_s\Th_a &&\text{by \ref{P4}}\\
&= \Th_{a^{-1}}\th_{\br(a^{-1})}\th_p\th_{\bd(a)}\Th_a \\
&= \Th_{a^{-1}}\th_p\Th_a &&\text{by Lemma \ref{lem:G2}\ref{G22}.}
\intertext{For \ref{G2c}, if $q\leq \br(a)$ then}
\Th_{a\rest_q} &= \th_{q\vt_{a^{-1}}}\Th_a &&\text{by \eqref{eq:Thaq}}\\
&= \Th_a\th_q\Th_{a^{-1}} \Th_a &&\text{by \ref{G2a}}\\
&= \Th_a\th_q\th_{\br(a)} &&\text{by \eqref{eq:Tha*Tha}}\\
&= \Th_a\th_q &&\text{by Lemma \ref{lem:thpthq}, as $q\leq\br(a)$.}  
\intertext{For \ref{G2d}, since $\vt_a$ is a bijection $s^\da\to t^\da$ by Lemma \ref{lem:vtavta*}, it is enough to show that $\vt_a$ is a projection algebra morphism, i.e.~that}
(p\th_q)\vt_a &= (p\vt_a)\th_{q\vt_a} &&\text{for all $p,q\leq s=\bd(a)$.}
\intertext{But for any such $p,q$, we have}
(p\vt_a)\th_{q\vt_a} &= p\vt_a\Th_{a^{-1}}\th_q\Th_a &&\text{by \ref{G2a}}\\
&= p\vt_a\th_t\vt_{a^{-1}}\th_q\th_s\vt_a &&\text{by definition of the $\Th$ maps}\\
&= p\vt_a\vt_{a^{-1}}\th_q\vt_a &&\text{as $p\vt_a\leq\br(a) = t$ and $q\leq s$ (cf.~Lemma \ref{lem:thpthq})}\\
&= p\th_q\vt_a &&\text{by Lemma \ref{lem:vtavta*}.}
\end{align*}

\pfitem{\ref{G2b}$\implies$\ref{G2a}}  This is immediate, since for any $p\leq\bd(a)=s$ we have $p\vt_a=(p\th_s)\vt_a=p\Th_a$ .

\pfitem{\ref{G2c}$\implies$\ref{G2a}}  If \ref{G2c} holds, then for any $p\leq \bd(a)=\br(a^{-1})$ we have
\begin{align*}
\Th_{a^{-1}}\th_p\Th_a &= \Th_{a^{-1}\rest_p}\Th_a &&\text{by \ref{G2c}}\\
&= \th_{p\vt_a}\Th_{a^{-1}}\Th_a &&\text{by \eqref{eq:Thaq}}\\
&= \th_{p\vt_a}\th_{\br(a)} &&\text{by \eqref{eq:Tha*Tha}}\\
&= \th_{p\vt_a} &&\text{by Lemma \ref{lem:thpthq}, as $p\vt_a\leq\br(a)$.}
\end{align*}

\pfitem{\ref{G2d}$\implies$\ref{G2a}}  Suppose \ref{G2d} holds, and let $p\leq s$.  We must show that 
\begin{align*}
e\th_{p\vt_a} &= e\Th_{a^{-1}}\th_p\Th_a &&\text{for any $e\in P$,}
\intertext{so fix some such $e$.  To make the following calculation easier to read, let $f=e\th_t$.  Since ${f\leq t=\bd(a^{-1})}$, we can also define $g=f\vt_{a^{-1}}$.  By Lemma \ref{lem:vtavta*} we have $f=g\vt_a$, and we note also that $g\leq\br(a^{-1})=s$.  We then have}
e\th_{p\vt_a} = e\th_t\th_{p\vt_a} &= f\th_{p\vt_a} &&\text{by Lemma \ref{lem:thpthq}, as $p\vt_a\leq\br(a)=t$}\\
&= (g\vt_a)\th_{p\vt_a} &&\text{as observed above}\\
&= (g\th_p)\vt_a &&\text{by \ref{G2d}, as $g,p\leq s=\bd(a)$}\\
&= e\th_t\vt_{a^{-1}}\th_p\vt_a &&\text{as $g=f\vt_{a^{-1}}$ and $f=e\th_t$}\\
&= e\th_t\vt_{a^{-1}}\th_p\th_s\vt_a &&\text{by Lemma \ref{lem:thpthq}, as $p\leq s$}\\
&= e\Th_{a^{-1}}\th_p\Th_a &&\text{by definition of the $\Th$ maps.} \qedhere
\end{align*}
\epf

\begin{rem}\label{rem:G1}
We briefly discussed property \ref{G2c} before the statement of Proposition~\ref{prop:G1}, and we now comment on the others.

First, we note that properties \ref{G2a} and \ref{G2b} bear some resemblance to the projection algebra axiom \ref{P4}.  In fact, this is more than just superficial.  Recall that (abstract) projection algebras are supposed to `model' projections of regular $*$-semigroups, and that the $\th_p$ operations model `conjugation', \emph{viz.}~$q\th_p \equiv pqp$.  In this way, axiom \ref{P4} can be thought of as a recipe for `iterating' conjugation by projections.  In the next section we will see that the $\vt$ and $\Th$ maps can also be interpreted as `conjugation' by \emph{morphisms} (see Remark \ref{rem:conj}).  In this way conditions~\ref{G2a} and~\ref{G2b} can then be thought of as recipes for iterating these kinds of conjugations.

By Lemma \ref{lem:vtavta*}, $\vt_a$ is always an order-isomorphism $\bd(a)^\da\to\br(a)^\da$ in a weak projection groupoid.  As we saw in Remark \ref{rem:Kinyon2}, this alone is not enough to guarantee that it is a projection algebra morphism, as in property \ref{G2d}.
\end{rem}

\subsection{Projection groupoids}\label{subsect:PG}

Proposition \ref{prop:G1} allows us to formalise the desired algebraic link between the components of a weak projection groupoid $(P,\G)$:

\begin{defn}\label{defn:PG}
A \emph{projection groupoid} is a weak projection groupoid $(P,\G)$ (cf.~Definition \ref{defn:WPG}), for which:
\begin{enumerate}[label=\textup{\textsf{(G1)}},leftmargin=10mm]%\addtocounter{enumi}{1}
\item \label{G2} any (and hence all) of the conditions \ref{G2a}--\ref{G2d} hold.
\end{enumerate}
We denote by $\PG$ the category of projection groupoids.  A morphism ${(P,\G)\to(P',\G')}$ in $\PG$ is an ordered groupoid functor $\phi:\G\to\G'$ whose object map $v\phi = \phi|_P$ is a projection algebra morphism $P\to P'$.
\end{defn}

The next result involves the groupoid $\G(S)$ and projection algebra $P(S)$, constructed from a regular $*$-semigroup $S$ in Definitions \ref{defn:GS} and \ref{defn:PS}.

\begin{prop}\label{prop:PGfunctor}
The assignment $S\mt(P(S),\G(S))$ is the object part of a functor ${\RSS\to\PG}$.
\end{prop}

\pf
Let $S$ be a regular $*$-semigroup.  We first check that $(P,\G)=(P(S),\G(S))$ is a projection groupoid.  
%Certainly $P=v\G$, and \ref{G1} holds by Definition \ref{defn:GS}.  
It follows quickly from Definition \ref{defn:GS} that $(P,\G)$ is a \emph{weak} projection groupoid, so it remains to verify \ref{G2}.
%Working towards \ref{G2}, 
For this, we first note that for any $a\in S$ and $p\leq\bd(a)$, we have
\begin{equation}
\label{eq:tb1} p\vt_a = \br({}_p\corest a) = \br(pa) = (pa)^*pa = a^*pa.
%\label{eq:tb1} p\vt_a = a^*pa \qquad\text{for all $a\in S$ and $p\leq\bd(a)$.}
%\label{eq:tb2}  t\Th_b &= b^*tb &&\text{for all $b\in S$ and $t\in P$.}
\end{equation}
%Indeed, we calculate $p\vt_a = \br({}_p\corest a) = \br(pa) = (pa)^*pa = a^*p^*pa = a^*pa$.  
%Using this, we have
%\[
%t\Th_b = t\th_{\bd(b)}\vt_b = b^*\cdot \bd(b)\cdot t\cdot \bd(b)\cdot b = b^*\cdot bb^*\cdot t\cdot bb^*\cdot b = b^*tb,
%\]
%giving \eqref{eq:tb2}.
To show that \ref{G2d} holds, let $a\in\G(=S)$.  We must show that $\vt_a$ is a projection algebra morphism,
%let $a\in\G(=S)$, and let $p,q\in\bd(a)^\da$.  We 
i.e.~that 
\[
(q\th_p)\vt_a = (q\vt_a)\th_{p\vt_a} \qquad\text{for all $p,q\in\bd(a)^\da$.}
\]
For this we use \eqref{eq:tb1} and $q=aa^*qaa^*$ (as $q\leq\bd(a)=aa^*$) to calculate
\[
(q\vt_a)\th_{p\vt_a} = (a^*qa)\th_{a^*pa} = a^*pa\cdot a^*qa\cdot a^*pa = a^*p(a a^*qa a^*)pa = a^*pqpa = (pqp)\vt_a = (q\th_p)\vt_a.
\]

%To show that \ref{G2a} holds, let $a\in\G(=S)$ and $p\leq\bd(a)=aa^*$.  Then for any $t\in P$, we use~\eqref{eq:tb1} and~\eqref{eq:tb2} to calculate
%\[
%t\Th_{a^{-1}}\th_p\Th_a = t\Th_{a^*}\th_p\Th_a = a^*(p(ata^*)p)a = a^*pa \cdot t \cdot a^*pa = t\th_{a^*pa} = t\th_{p\vt_a},
%\]
%so that $\Th_{a^{-1}}\th_p\Th_a = \th_{p\vt_a}$, as required.

Now let $S'$ be another regular $*$-semigroup, and write $P'=P(S')$ and $\G'=\G(S')$.  We must show that any $*$-morphism $\phi:S\to S'$ is a projection groupoid morphism $(P,\G)\to(P,\G')$.  But this follows from Propositions \ref{prop:calGfunctor} and \ref{prop:Pfunctor}, which respectively show that $\phi$ is an ordered groupoid functor $\G\to\G'$, and $\phi|_P$ is a projection algebra morphism $P\to P'$.
\epf

%\pf
%Let $S$ be a regular $*$-semigroup.  We first check that $(P,\G)=(P(S),\G(S))$ is a projection groupoid.  Certainly $P=v\G$, and \ref{G1} holds by Definition \ref{defn:GS}.  Working towards \ref{G2}, we first claim that
%\begin{align}
%\label{eq:tb1} p\vt_a &= a^*pa &&\text{for all $a\in S$ and $p\leq\bd(a)$,}\\
%\label{eq:tb2}  t\Th_b &= b^*tb &&\text{for all $b\in S$ and $t\in P$.}
%\end{align}
%Indeed, for \eqref{eq:tb1}, we calculate $p\vt_a = \br({}_p\corest a) = \br(pa) = (pa)^*pa = a^*p^*pa = a^*pa$.  Using this, we have
%\[
%t\Th_b = t\th_{\bd(b)}\vt_b = b^*\cdot \bd(b)\cdot t\cdot \bd(b)\cdot b = b^*\cdot bb^*\cdot t\cdot bb^*\cdot b = b^*tb,
%\]
%giving \eqref{eq:tb2}.
%
%To show that \ref{G2a} holds, let $a\in\G(=S)$ and $p\leq\bd(a)=aa^*$.  Then for any $t\in P$, we use~\eqref{eq:tb1} and~\eqref{eq:tb2} to calculate
%\[
%t\Th_{a^{-1}}\th_p\Th_a = t\Th_{a^*}\th_p\Th_a = a^*(p(ata^*)p)a = a^*pa \cdot t \cdot a^*pa = t\th_{a^*pa} = t\th_{p\vt_a},
%\]
%so that $\Th_{a^{-1}}\th_p\Th_a = \th_{p\vt_a}$, as required.
%
%Now let $S'$ be another regular $*$-semigroup, and write $P'=P(S')$ and $\G'=\G(S')$.  We must show that any $*$-morphism $\phi:S\to S'$ is a projection groupoid morphism $(P,\G)\to(P,\G')$.  But this follows from Propositions \ref{prop:calGfunctor} and \ref{prop:Pfunctor}, which respectively show that $\phi$ is an ordered groupoid functor $\G\to\G'$, and $\phi|_P$ is a projection algebra morphism $P\to P'$.
%\epf

\begin{rem}\label{rem:conj}
It follows from \eqref{eq:tb1} that $\vt_a$ acts as a kind of `conjugation' by the morphism $a\in\G=\G(S)$ on projections from $\bd(a)^\da$.  The $\Th$ maps can be thought of as conjugations as well, but on all of $P$.  Indeed, for any $a\in\G$ and $p\in P$, we use \eqref{eq:tb1} to calculate
\begin{equation}
\label{eq:tb2} p\Th_a = p\th_{\bd(a)}\vt_a = a^*\cdot \bd(a)\cdot p\cdot \bd(a)\cdot a = a^*\cdot aa^*\cdot p\cdot aa^*\cdot a = a^*pa.
\end{equation}
\end{rem}

\begin{rem}\label{rem:Rees1}
Certainly the functor $\RSS\to\PG$ from Proposition \ref{prop:PGfunctor} carries more information than the functor $\RSS\to\OG$ from Proposition \ref{prop:calGfunctor}.  For example, the latter did not allow us to distinguish adjacency semigroups $S=A(\Ga)$ and $S'=A(\Ga')$, where $\Ga$ and $\Ga'$ were different (reflexive, undirected) graphs on the same vertex set (cf.~Example \ref{eg:AG}).  However, since $\Ga$ and~$\Ga'$ have different edge sets, it follows from \eqref{eq:AGth} that the projection algebras $P(S)$ and $P(S')$ are different, and so $(P(S),\G(S)) \not= (P(S'),\G(S'))$.

On the other hand, the functor $\RSS\to\PG$ from Proposition \ref{prop:PGfunctor} is still not injective (on objects).  Indeed, if $S=\M^0(P,P,G,M)$ is a Rees 0-matrix regular $*$-semigroup as in Example~\ref{eg:Rees}, then the groupoid $\G(S)$ depends only on the indexing set $P$ and the group $G$, while the projection algebra~$P(S)$ depends only on the \emph{locations} of the non-zero entries of the sandwich matrix~$M$, and not on the \emph{values} of these entries; cf.~\eqref{eq:M0th}.  
\end{rem}

Before we move on, we record a simple technical result that will be used later.

\begin{lemma}\label{lem:pa'qap}
If $(P,\G)$ is a projection groupoid, and if $a\in\G$ and $p,q\in P$, then
\[
p\Th_{a^{-1}}\th_q\Th_a\th_p = q\Th_a\th_p.
\]
\end{lemma}

\pf
By \ref{G2b} and \ref{P3} we have $p\Th_{a^{-1}}\th_q\Th_a\th_p = p\th_{q\Th_a}\th_p = q\Th_a\th_p$.
\epf

\subsection{Weak chained projection groupoids}\label{subsect:ChPG}

The functor $\RSS\to\PG$ from Proposition \ref{prop:PGfunctor} allows us to represent a regular $*$-semigroup~$S$ by its corresponding projection groupoid $(P(S),\G(S))$, but we saw in Remark \ref{rem:Rees1} that this representation is not faithful.  We address this by adding an extra layer of structure:

\begin{defn}\label{defn:ve}
Let $(P,\G)$ be a projection groupoid (cf.~Definition \ref{defn:PG}), and let $\C=\C(P)$ be the chain groupoid of $P=v\G$ (cf.~Definition \ref{defn:CP}).  An \emph{evaluation map} is an ordered $v$-functor $\ve:\C\to\G$, meaning that the following hold:
\begin{enumerate}[label=\textup{\textsf{(E\arabic*)}},leftmargin=10mm]
\item \label{E1} $\ve(p)=p$ for all $p\in P$ (where as usual we identify $p\equiv[p]\in\C$),
\item \label{E2} $\ve(\c\circ\d)=\ve(\c)\circ\ve(\d)$ if $\br(\c)=\bd(\d)$,
\item \label{E3} $\c\leq\d \implies \ve(\c)\leq\ve(\d)$.
\end{enumerate}
\end{defn}

\begin{rem}\label{rem:ve}
It is worthwhile noting a number of simple consequences of \ref{E1}--\ref{E3}.  As with group homomorphisms, it is easy to see that
\begin{enumerate}[label=\textup{\textsf{(E\arabic*)}},leftmargin=10mm]\addtocounter{enumi}{3}
\item \label{E4} $\ve(\c^{-1})=\ve(\c)^{-1}$ for all $\c\in\C$.  
\end{enumerate}
It also follows that:
\begin{enumerate}[label=\textup{\textsf{(E\arabic*)}},leftmargin=10mm]\addtocounter{enumi}{4}
\item \label{E5} $\bd(\ve(\c))=\bd(\c)$ and $\br(\ve(\c))=\br(\c)$ for all $\c\in\C$.  
\end{enumerate}
For example, $\bd(\ve(\c)) = \ve(\c)\circ\ve(\c)^{-1} = \ve(\c\circ\c^{-1}) = \ve(\bd(\c))$.

It is also easy to see that \ref{E3} is equivalent (in the presence of the other axioms) to:
\begin{enumerate}[label=\textup{\textsf{(E\arabic*)}},leftmargin=10mm]\addtocounter{enumi}{5}
\item \label{E6} $\ve({}_q\corest\c)={}_q\corest\ve(\c)$ if $q\leq\bd(\c)$.
\end{enumerate}
For example, if \ref{E6} holds, and if $\c\leq\d$, then $\c={}_p\corest\d$ for some $p\leq\bd(\d)$; it follows from this that $\ve(\c) = \ve({}_p\corest\d) = {}_p\corest\ve(\d) \leq\ve(\d)$.
\end{rem}

\begin{rem}
The existence of an evaluation map $\C=\C(P)\to\G$ (for a projection groupoid $(P,\G)$) implies a certain degree of connectivity in the groupoid $\G$.  Indeed, since there is a $P$-chain $p\to q$ for every pair $(p,q)$ in the transitive closure of $\F$, the morphism set $\G(p,q)$ 
%contains the image of such a chain, and hence 
must be non-empty for every such pair.  For any groupoid $\G$ with (at least) this much connectivity, freeness of $\C$ (cf.~Remark \ref{rem:Cfree}) means that one can define a $v$-functor $\phi:\C\to\G$ by arbitrarily choosing morphisms $[p,q]\phi\in\G(p,q)$ for each $(p,q)\in{\F}\sm\De_P$, as long as one ensures that the identity $[q,p]\phi=([p,q]\phi)^{-1}$ is satisfied.  Whether or not this $\phi$ will be an \emph{ordered} $v$-functor is another question entirely.
\end{rem}

\begin{defn}\label{defn:ChPG}
A \emph{weak chained projection groupoid} is a triple $(P,\G,\ve)$, where $(P,\G)$ is a projection groupoid (cf.~Definition \ref{defn:PG}), and $\ve:\C\to\G$ is an evaluation map (cf.~Definition~\ref{defn:ve}), where $\C=\C(P)$.  

We write $\WCPG$ for the category of weak chained projection groupoids.  Morphisms in $\WCPG$ are called \emph{chained projection functors}.  Such a morphism~${(P,\G,\ve)\to(P',\G',\ve')}$ is a projection groupoid morphism $\phi:(P,\G)\to(P',\G')$ (cf.~Definition \ref{defn:PG}) that \emph{respects evaluation maps} in the sense that the following diagram commutes:
\[
\begin{tikzcd}
\C \arrow{rr}{\Phi} \arrow[swap]{dd}{\ve} & ~ & \C' \arrow{dd}{\ve'} \\%
~&~&~\\
\G \arrow{rr}{\phi}& ~ & \G',
\end{tikzcd}
\]
where here $\Phi = \C(v\phi)$ is the functor $\C\to\C'=\C(P')$ from Proposition \ref{prop:CC'}, given by $[p_1,\ldots,p_k]\Phi = [p_1\phi,\ldots,p_k\phi]$.  Commutativity of the diagram means that
\[
(\ve(\c))\phi = \ve'(\c\Phi) \qquad\text{for all $\c\in\C$.}
\]
(To check that the composition of chained projection functors is indeed a chained projection functor, one uses the functoriality of $\C$, as in Remark \ref{rem:Cfunctor}.)
\end{defn}

Consider a weak chained projection groupoid $(P,\G,\ve)$.  
For $p,q\in P$ with $p\F q$, we have the $P$-chain $[p,q]\in\C$, and the morphisms $\ve[p,q] \in \G$ will play a crucial role in all that follows.  For one thing, these morphisms generate the image of $\C$ under $\ve$ (as the $[p,q]$ generate $\C$).  Specifically, if $\c=[p_1,p_2,\ldots,p_k]\in\C$, then 
\begin{equation}\label{eq:vec}
\ve(\c) = \ve[p_1,p_2]\circ\ve[p_2,p_3]\circ\cdots\circ\ve[p_{k-1},p_k].
\end{equation}
The next lemma gathers some important basic properties of the $\ve[p,q]$; in what follows, we will typically use these without explicit reference.

\begin{lemma}\label{lem:ve}
Let $(P,\G,\ve)$ be a weak chained projection groupoid.
\ben
\item \label{ve1} For any $p\in P$ we have $\ve[p,p]=p$.
\item \label{ve2} If $p,q\in P$ are such that $p\F q$, then
\begin{align*}
\bd(\ve[p,q]) &= p ,& \ve[p,q]^{-1} &= \ve[q,p], &\vt_{\ve[p,q]} &= \th_q|_{p^\da},\\
\br(\ve[p,q]) &= q ,& && \Th_{\ve[p,q]} &= \th_p\th_q.
\end{align*}
\item \label{ve3} If $p,q,r,s\in P$ are such that $p\F q$, $r\leq p$ and $s\leq q$, then
\[
{}_r\corest\ve[p,q] = \ve[r,r\th_q] \AND \ve[p,q]\rest_s = \ve[s\th_p,s].
\]
\een
\end{lemma}

\pf
\firstpfitem{\ref{ve1}}  This follows quickly from \ref{E1} and $[p,p]=[p]$, keeping in mind the identification of $p\in P$ with the chain $[p]\in\C$.

\pfitem{\ref{ve3}}  This follows from \ref{E6}, together with \eqref{eq:rest} and \eqref{eq:corest}.

\pfitem{\ref{ve2}}  The claims concerning domains, ranges and inversion follow from \ref{E4} and \ref{E5}, and the fact that $[p,q]^{-1}=[q,p]$ in $\C$.  

Since the maps $\vt_{\ve[p,q]}$ and $\th_q|_{p^\da}$ both have domain $p^\da$, we can show the maps are equal by showing that
\[
r\vt_{\ve[p,q]} = r\th_q \qquad\text{for all $r\leq p$.}
\]
But for such $r$, we use the definition of the $\vt$ maps, and parts of the current lemma that have already been proved, to calculate
\[
r\vt_{\ve[p,q]} = \br({}_r\corest\ve[p,q]) = \br(\ve[r,r\th_q]) = r\th_q.
\]
Finally, again using already-proved parts of the current lemma, we have
\[
\Th_{\ve[p,q]} = \th_p\vt_{\ve[p,q]} = \th_p\circ \th_q|_{p^\da} = \th_p\th_q,
\]
where in the last step we used the fact that $p^\da=\im(\th_p)$.
\epf

Before we move on, it will also be convenient to record the following result.

\begin{lemma}\label{lem:aepq}
Let $(P,\G,\ve)$ be a weak chained projection groupoid.  If $p,q\in P$ are such that $p\F q$, and if $a\leq\ve[p,q]$, then 
\[
a=\ve[r,s] \WHERE r=\bd(a) \ANd s=\br(a).
\]
\end{lemma}

\pf
We use Lemma \ref{lem:ve} to calculate
\[
a = {}_r\corest\ve[p,q] = \ve[r,r\th_q] \AND r\th_q = \br(\ve[r,r\th_q]) = \br(a)=s.  \qedhere
\]
\epf

For the rest of Section \ref{subsect:ChPG} we fix a regular $*$-semigroup $S$, and let $(P,\G) = (P(S),\G(S))$ be the projection groupoid associated to $S$ (cf.~Proposition \ref{prop:PGfunctor}).  We wish to extend this to a \emph{weak chained} projection groupoid, for which we must define an evaluation map $\ve = \ve(S):\C\to\G$, where $\C=\C(P)$ is the chain groupoid of $P$.  % by $\ve[p_1,\ldots,p_k] = p_1\cdots p_k$.
To do so, it is first convenient to define a map 
\begin{equation}\label{eq:piS}
\pi=\pi(S):\P\to\G(=S) \BY \pi(p_1,\ldots,p_k) = p_1\cdots p_k,
\end{equation}
where $\P=\P(P)$ is the path category of $P$.  So $\pi(\p)$ is simply the product (in $S$) of the entries of $\p$ (taken in order).  We begin with the following fact:

\begin{lemma}\label{lem:dpip}
For any $\p\in\P$, we have $\bd(\pi(\p)) = \bd(\p)$ and $\br(\pi(\p)) = \br(\p)$.
\end{lemma}

\pf
Write $\p=(p_1,\ldots,p_k)$.  A routine calculation shows that
\[
\bd(\pi(\p)) = \pi(\p)\cdot\pi(\p)^* = p_1\cdots p_{k-1}\cdot p_k\cdot p_{k-1}\cdots p_1.
\]
Since $p_1\F\cdots\F p_k$, this reduces to $p_1=\bd(\p)$.  We obtain $\br(\pi(\p)) = p_k = \br(\p)$ by symmetry.
\epf

The next result concerns the congruence ${\approx}=\Om^\sharp$ on $\P$ from Definition \ref{defn:approx}.

\begin{lemma}\label{lem:GS4}
The map $\pi=\pi(S):\P\to\G$ is a $v$-functor, and ${\approx}\sub\ker(\pi)$.
\end{lemma}

\pf
For the first statement, it suffices by Lemma \ref{lem:dpip} to show that $\pi$ is a functor.  For this, consider composable paths $\p=(p_1,\ldots,p_k)$ and $\q=(p_k,\ldots,p_l)$.  Then since $p_k$ is an idempotent, and since $\br(\pi(\p)) = p_k = \bd(\pi(\q))$ by Lemma \ref{lem:dpip}, we have
\begin{align*}
\pi(\p\circ\q) = \pi(p_1,\ldots,p_l) &= p_1\cdots p_kp_{k+1}\cdots p_l \\
&= p_1\cdots p_k\cdot p_kp_{k+1}\cdots p_l = \pi(\p)\pi(\q) = \pi(\p)\circ\pi(\q).
\end{align*}
To show that ${\approx}\sub\ker(\pi)$, we must show that $\Om\sub\ker(\pi)$, i.e.~that $\pi(\s)=\pi(\t)$ for all $(\s,\t)\in\Om$.  This is clear if $(\s,\t)$ has the form \ref{Om1}.  For \ref{Om2}, recall that $p=q\th_p=pqp$ when $p\F q$.
\epf

\begin{defn}\label{defn:veS}
Since $\C=\P/{\approx}$, it follows from Lemma \ref{lem:GS4} that we have a well-defined $v$-functor
\[
\ve = \ve(S):\C \to \G \GIVENBY \ve[\p] = \pi(\p) \qquad\text{for $\p\in\P$.}
\]
That is, $\ve[p_1,\ldots,p_k] = p_1\cdots p_k$ whenever $p_1\F\cdots\F p_k$.  
\end{defn}

\begin{lemma}\label{lem:veS}
The functor $\ve=\ve(S):\C\to\G$ is an evaluation map.
\end{lemma}

\pf
Conditions \ref{E1} and \ref{E2} hold because $\ve$ is a $v$-functor.  As explained in Remark \ref{rem:ve}, we can prove \ref{E6} in place of \ref{E3}, and we note that \ref{E6} says
\begin{align}
\nonumber \ve( {}_q\corest[\p]) &= {}_q\corest\ve[\p] &&\text{for all $\p\in\P$ and $q\leq\bd[\p]=\bd(\p)$.}
\intertext{Since $\ve( {}_q\corest[\p]) = \ve[{}_q\corest\p] = \pi({}_q\corest\p)$ and ${}_q\corest\ve[\p]={}_q\corest\pi(\p)$, this is equivalent to}
\label{eq:prodqp}
\pi( {}_q\corest\p) &= {}_q\corest\pi(\p) &&\text{for all $\p\in\P$ and $q\leq\bd(\p)$.}
\end{align}
We prove this by induction on $k$, the length of the path $\p=(p_1,\ldots,p_k)$.  When $k=1$, we have
\[
\pi({}_q\corest\p) = \pi({}_q\corest(p_1)) = \pi (q) = q = qp_1 = q\pi(\p) = {}_q\corest\pi(\p),
\]
where we used the fact that $q\leq\bd(\p)=p_1 \implies q=qp_1$.  Now suppose $k\geq2$, and let ${}_q\corest\p = (q_1,\ldots,q_k)$ be as in \eqref{eq:rest}.  Let $\p'=(p_1,\ldots,p_{k-1})\in\P$, noting that ${}_q\corest\p'=(q_1,\ldots,q_{k-1})$.  Then
\begin{align*}
{}_q\corest\pi(\p) = q\cdot p_1\cdots p_k &= q p_1\cdots p_k \cdot (q p_1\cdots p_k)^*\cdot q p_1\cdots p_k \\
&= q p_1\cdots p_k \cdot p_k\cdots p_1q \cdot qp_1\cdots p_k \\
&= q\cdot (p_1\cdots p_{k-1}) \cdot (p_k\cdots p_1\cdot q\cdot p_1\cdots p_k) \\
&=  {}_q\corest\pi(\p') \cdot q\th_{p_1}\cdots\th_{p_k} \\
&=  \pi({}_q\corest\p') \cdot q_k &&\hspace{-2cm}\text{by induction and \eqref{eq:qi}}\\
&=  q_1\cdots q_{k-1}\cdot q_k = \pi({}_q\corest \p).
\end{align*}
This completes the proof of \eqref{eq:prodqp}, and hence of the lemma.
\epf

\begin{defn}\label{defn:PGveS}
It follows from Proposition \ref{prop:PGfunctor} and Lemma \ref{lem:veS} that $(P(S),\G(S),\ve(S))$ is a weak chained projection groupoid for any regular $*$-semigroup $S$ (cf.~Definitions \ref{defn:GS}, \ref{defn:PS} and~\ref{defn:veS}).  
\end{defn}

\begin{prop}\label{prop:PGvefunctor}
%The assignment $S\mt\bG(S) = (P(S),\G(S),\ve(S))$ is the object part of a functor ${\RSS\to\WCPG}$.
%The assignment $S\mt\bG(S)$ is the object part of a functor ${\RSS\to\WCPG}$.
The assignment $S\mt(P(S),\G(S),\ve(S))$ is the object part of a functor ${\RSS\to\WCPG}$.
\end{prop}

\pf
It remains to check that any morphism $\phi:S\to S'$ in $\RSS$ is a chained projection functor from $(P,\G,\ve)=(P(S),\G(S),\ve(S))$ to $(P',\G',\ve')=(P(S'),\G(S'),\ve(S'))$.  Since $\phi$ is a projection groupoid morphism (cf.~Proposition \ref{prop:PGfunctor}), it remains to check that it respects the evaluation maps (cf.~Definition \ref{defn:ChPG}).  For this, fix a chain ${\c=[p_1,\ldots,p_k]\in\C=\C(P)}$.  Writing $\Phi=\C(\phi):\C\to\C'=\C(P')$ as in Proposition \ref{prop:CC'}, we have
\[
(\ve(\c))\phi = (p_1\cdots p_k)\phi = (p_1\phi)\cdots(p_k\phi) = \ve'[p_1\phi,\ldots,p_k\phi] = \ve'(\c\Phi).  \qedhere
\]
\epf

We will ultimately see that the functor $\RSS\to\WCPG$ from Proposition \ref{prop:PGvefunctor} is injective (on objects), although it is still not surjective.  Rather than proving injectivity directly, we will first identify the image of the functor, and then construct an inverse functor in the opposite direction.  The image will consist of the \emph{chained projection groupoids} defined in Section \ref{subsect:CoPG}, and the inverse functor will be constructed in Section \ref{sect:GtoS}.

Before we do this, we can already check that the functor from Proposition \ref{prop:PGvefunctor} allows us to distinguish between Rees 0-matrix semigroups.

\begin{rem}\label{rem:Rees2}
Consider again the Rees 0-matrix regular $*$-semigroup $S=\M^0(P,P,G,M)$, as in Example \ref{eg:Rees}.  We saw previously that the projection algebra of $S$ is $P_0=P\cup\{0\}$, with $\th$ operations as in \eqref{eq:M0th}.  As in \eqref{eq:vec}, the evaluation map $\ve=\ve(S):\C(P_0)\to S$ is completely determined by the elements $\ve[p,q]$, for $(p,q)\in{\F}$.  Of course we have $\ve[0,0]=0$.  Consulting~\eqref{eq:M0F}, the other $\F$-related pairs in $P_0$ have the form $(p,q)$, for $p,q\in P$ with $m_{pq}\not=0$, and for these we have
\begin{equation}\label{eq:Rees2}
\ve[p,q] = pq \equiv(p,1,p)(q,1,q) = (p,m_{pq},q).
\end{equation}
In this way, the evaluation map $\ve=\ve(S)$ specifies the entries of the sandwich matrix $M$, and hence the entire structure of $S$, so that sandwich matrices and evaluation maps are two sides of the same coin.
\end{rem}

\subsection{Chained projection groupoids}\label{subsect:CoPG}

We are almost ready to give the definition of chained projection groupoids, as the weak chained projection groupoids satisfying a certain coherence condition.  A key ingredient is the concept of a \emph{linked pair} of projections, the definition of which does not require the existence of an evaluation map.

\begin{defn}\label{defn:LP1}
Let $(P,\G)$ be a projection groupoid (cf.~Definition \ref{defn:PG}), and consider a morphism $b\in\G$.  A pair of projections $(e,f)\in P\times P$ is said to be \emph{$b$-linked} if 
\begin{equation}\label{eq:LP}
f = e\Th_b\th_f \AND e = f\Th_{b^{-1}}\th_e.
\end{equation}
\end{defn}

We will soon associate two morphisms $\lam(e,b,f)$ and $\rho(e,b,f)$ to the $b$-linked pair $(e,f)$.  But first we prove the next result, which will ensure these morphisms are well defined.

\newpage

\begin{lemma}\label{lem:LP}
Let $(P,\G)$ be a projection groupoid, and suppose $(e,f)$ is $b$-linked, where ${b\in\G(q,r)}$.  Define further projections
\begin{equation}\label{eq:e1e2f1f2}
e_1 = e\th_q \COMMA e_2 = f\Th_{b^{-1}} \COMMA f_1 = e\Th_b \AND f_2 = f\th_r.
\end{equation}
Then
\ben\bmc2
\item \label{LP1} $e\leqF q$ and $f\leqF r$,
\item \label{LP2} $e_i\leq q$ and $f_i\leq r$ for $i=1,2$,
\item \label{LP3} $e\F e_i$ and $f\F f_i$ for $i=1,2$,
\item \label{LP4} ${}_{e_i}\corest b=b\rest_{f_i}$ for $i=1,2$.
\emc\een
\end{lemma}

\pf
It is clear from Definition \ref{defn:LP1} that $(e,f)$ is $b$-linked if and only if $(f,e)$ is $b^{-1}$-linked.  If we replace $b\leftrightarrow b^{-1}$ and $(e,f)\leftrightarrow(f,e)$, then the projections defined in \eqref{eq:e1e2f1f2} are interchanged accordingly, \emph{viz.}~$e_1\leftrightarrow f_2$ and $e_2\leftrightarrow f_1$.
Because of this symmetry, for parts \ref{LP1}--\ref{LP3}, it suffices to prove the claims concerning $e,e_1,e_2$.  (Alternatively, the arguments below can be easily adapted to prove the claims regarding $f,f_1,f_2$.)

\pfitem{\ref{LP2}}  This follows from $\im(\Th_{b^{-1}})=\im(\th_q)=q^\da$.

\pfitem{\ref{LP3}}  From $e\leq e\leqF q$, it follows from Lemma \ref{lem:pqr}\ref{pqr2} that $e\F e\th_q=e_1$.  By Lemma \ref{lem:pqp}\ref{pqP1}, we have $e=e_2\th_e\leqF e_2$.  It remains to show that $e_2 \leqF e$, and for this we use \ref{G2b}, \eqref{eq:LP} and~\eqref{eq:e1e2f1f2}:
\[
e\th_{e_2} = e\th_{f\Th_{b^{-1}}} = e\Th_b\th_f\Th_{b^{-1}} = f\Th_{b^{-1}} = e_2.
\]

\pfitem{\ref{LP4}}  We must show that $f_i = e_i\vt_b$ ($i=1,2$).  Keeping $q=\bd(b)$ and $r=\bd(b^{-1})$ in mind, we have
\[
e_1\vt_b = e\th_q\vt_b = e\Th_b = f_1 \AND e_2\vt_b = f\Th_{b^{-1}}\vt_b = f\th_r = f_2,
\]
where we used Lemma \ref{lem:G2}\ref{G23} in the second calculation.

\pfitem{\ref{LP1}}  Combining parts \ref{LP2} and \ref{LP3}, we have $e\leqF e_1\leq q$, and Lemma~\ref{lem:pqr}\ref{pqr1} then gives ${e\leqF q}$.
\epf

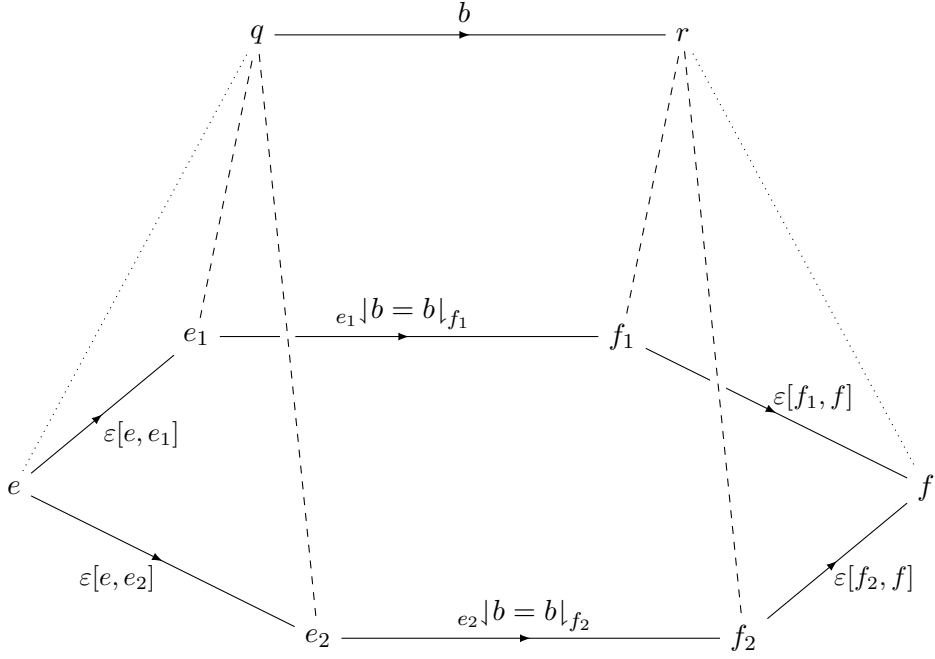
\begin{figure}[h]
\begin{center}
\begin{tikzpicture}[xscale=0.8]
\tikzstyle{vertex}=[circle,draw=black, fill=white, inner sep = 0.07cm]
\node (e) at (0,0){$e$};
\node (e1) at (3,2){$e_1$};
\node (e2) at (5,-2){$e_2$};
\node (f) at (15,0){$f$};
\node (f1) at (10,2){$f_1$};
\node (f2) at (12,-2){$f_2$};
\node (q) at (4,6){$q$};
\node (r) at (11,6){$r$};
\draw[->-=0.5] (q)--(r);
\draw[->-=0.5] (e)--(e1);
\draw[->-=0.5] (e1)--(f1);
\draw[->-=0.5] (f1)--(f);
\draw[->-=0.5] (e)--(e2);
\draw[->-=0.5] (e2)--(f2);
\draw[->-=0.5] (f2)--(f);
\draw[white,line width=2mm] (q)--(e2);
\draw[white,line width=2mm] (r)--(f2);
\draw[dashed] (e1)--(q)--(e2) (f1)--(r)--(f2);
\draw[dotted] (e)--(q) (f)--(r);
\node () at (7.4,6.3) {$b$};
\node () at (6.4,2.3) {${}_{e_1}\corest b = b\rest_{f_1}$};
\node () at (8.4,-1.7) {${}_{e_2}\corest b = b\rest_{f_2}$};
\node () at (2.1,0.7) {{\small $\ve[e,e_1]$}};
\node () at (1.7,-1.2) {{\small $\ve[e,e_2]$}};
\node () at (13.15,1.2) {{\small $\ve[f_1,f]$}};
\node () at (14.15,-1.2) {{\small $\ve[f_2,f]$}};
\end{tikzpicture}
\caption{The projections and morphisms associated to a $b$-linked pair $(e,f)$, for $b\in\G(q,r)$; see Definitions \ref{defn:LP1}, \ref{defn:LP2} and \ref{defn:CoPG}, and Lemma \ref{lem:LP}.  Dotted and dashed lines indicate $\leqF$ and~$\leq$ relationships, respectively.  Axiom \ref{G3} says that the hexagon at the bottom of the diagram commutes.}
\label{fig:G3}
\end{center}
\end{figure}

\begin{defn}\label{defn:LP2}
Let $(P,\G,\ve)$ be a weak chained projection groupoid (cf.~Definition \ref{defn:ChPG}).  Let~$(e,f)$ be a $b$-linked pair, where $b\in\G(q,r)$, and let $e_1,e_2,f_1,f_2\in P$ be as in \eqref{eq:e1e2f1f2}.  By Lemma~\ref{lem:LP}, we have two well-defined morphisms
\begin{align}
\nonumber \lam(e,b,f) &= \ve[e,e_1]\circ {}_{e_1}\corest b \circ \ve[f_1,f] &\text{and}&& \rho(e,b,f) &= \ve[e,e_2]\circ {}_{e_2}\corest b \circ \ve[f_2,f]\\
\label{eq:G3} &= \ve[e,e_1]\circ b\rest_{f_1} \circ \ve[f_1,f] &&&  &= \ve[e,e_2]\circ b\rest_{f_2} \circ \ve[f_2,f].
\end{align}
These morphisms are shown in Figure \ref{fig:G3}.
\end{defn}

The above-mentioned coherence property states that the two terms in \eqref{eq:G3} must be equal:

\begin{defn}\label{defn:CoPG}
A \emph{chained projection groupoid} is a weak chained projection groupoid $(P,\G,\ve)$ (cf.~Definition \ref{defn:ChPG}) satisfying the following coherence condition:
\begin{enumerate}[label=\textup{\textsf{(G2)}},leftmargin=10mm]%\addtocounter{enumi}{2}
\item \label{G3} For every $b\in\G$, and for every $b$-linked pair $(e,f)$, we have $\lam(e,b,f) = \rho(e,b,f)$, where these morphisms are as in \eqref{eq:G3}.
\end{enumerate}
We denote by $\CPG$ the full subcategory of $\WCPG$ consisting of all chained projection groupoids (and all chained projection functors between them).
\end{defn}

\begin{rem}\label{rem:G3'}
One might wonder if the concept of linked pairs could be avoided when defining chained projection groupoids.  Specifically, one might wonder if \ref{G3} is equivalent to:
\begin{enumerate}[label=\textup{\textsf{(G2)$'$}},leftmargin=10mm]%\addtocounter{enumi}{2}
\item \label{G3'} For any morphism $b\in\G(q,r)$, and for any projections $e,e_1,e_2,f,f_1,f_2\in P$ satisfying conditions~\ref{LP1}--\ref{LP4} from Lemma \ref{lem:LP}, we have $\ve[e,e_1]\circ {}_{e_1}\corest b \circ \ve[f_1,f] = \ve[e,e_2]\circ {}_{e_2}\corest b \circ \ve[f_2,f]$.
\end{enumerate}
While \ref{G3'} of course implies \ref{G3}, the converse does not hold, as we will see in Remark \ref{rem:G3'2}.
%.  Most significantly, the stronger condition \ref{G3'} does not (generally) hold in the all-important case that $\G=\G(S)$ is the groupoid associated to a regular $*$-semigroup~$S$, as in Definition \ref{defn:GS}.  We will say more about this in Remark \ref{rem:G3'2}.
\end{rem}

For a regular $*$-semigroup $S$, we denote by $\bG(S) = (P(S),\G(S),\ve(S))$ the weak chained projection groupoid from Definition \ref{defn:PGveS}.  Proposition \ref{prop:PGvefunctor} says that the assignment $S\mt\bG(S)$ is the object part of a functor $\RSS\to\WCPG$.  The next result says that $\bG$ maps into the coherent subcategory~$\CPG$.

\begin{thm}\label{thm:Gfunctor}
%The assignment $S\mt(P(S),\G(S),\ve(S))$ is the object part of 
$\bG$ is a functor ${\RSS\to\CPG}$.
\end{thm}

\pf
Given Proposition \ref{prop:PGvefunctor}, it remains to check that $(P,\G,\ve)=\bG(S)$ satisfies \ref{G3} for any regular $*$-semigroup $S$.  To do so, let $(e,f)$ be a $b$-linked pair, where ${b\in\G(=S)}$.  Write
\[
q = \bd(b) = bb^* \AND r = \br(b) = b^*b,
\]
and let the $e_i,f_i$ be as in \eqref{eq:e1e2f1f2}.  Keeping \eqref{eq:tb2} in mind, we have
\[
e_1 = qeq \COMMA e_2 = bfb^* \COMMA f_1 = b^*eb \AND f_2 = rfr.
\]
We then calculate $ee_1 = eqeq = eq$ and $bf_1 = bb^*eb = qeb$.  Since $e\leqF q$ by Lemma \ref{lem:LP}\ref{LP1}, we also have $e=q\th_e=eqe$.  Putting this all together, we have
\[
\lam(e,b,f) = \ve[e,e_1]\circ {}_{e_1}\corest b \circ \ve[f_1,f] = ee_1 \cdot e_1b \cdot f_1f = ee_1\cdot bf_1\cdot f = eq\cdot qeb \cdot f = eqe\cdot bf = ebf.
\]
A similar calculation gives $\rho(e,b,f) = ebf$, and the proof is complete.
\epf

\begin{rem}\label{rem:G3'2}
In the above proof, verification of \ref{G3} boils down to checking that
\[
ee_1bf_1f = ee_2bf_2f,
\]
where $b\in\G(q,r)$, $(e,f)$ is $b$-linked, and the $e_i,f_i$ are as in \eqref{eq:e1e2f1f2}.  We are now in a position to see that condition \ref{G3} cannot be replaced by the stronger \ref{G3'} discussed in Remark \ref{rem:G3'}.  Indeed, consider the partition monoid $\PP_4$, and the elements defined by
\[
b = \custpartn{1,2,3,4}{1,2,3,4}{\stline11\stline22\stline33\stline44} \COMMA
e=e_1=f_1 = \custpartn{1,2,3,4}{1,2,3,4}{\uarc12\darc12\stline33\stline44} \COMMA
e_2=f_2 = \custpartn{1,2,3,4}{1,2,3,4}{\uarc13\darc13\stline22\stline44} \AND
f = \custpartn{1,2,3,4}{1,2,3,4}{\uarc14\darc14\stline22\stline33} .
\]
These elements are all projections (so in particular $q=r=b$), and they satisfy conditions~\ref{LP1}--\ref{LP4} of Lemma \ref{lem:LP}.  However, we have
\[
ee_1bf_1f=
\custpartn{1,2,3,4}{1,2,3,4}{\uarc12\darc14\stline33\stline42}
\not=
\custpartn{1,2,3,4}{1,2,3,4}{\uarc12\darc14\stline32\stline43}
=ee_2bf_2f.
\]
\end{rem}

\section{\boldmath From chained projection groupoids to regular $*$-semigroups}\label{sect:GtoS}

In Section \ref{sect:G} we introduced the category $\CPG$ of chained projection groupoids, and constructed a functor $\bG:\RSS\to\CPG$ from the category of regular $*$-semigroups.  Our goal in this section is to go backwards, and construct a functor $\bS:\CPG\to\RSS$.  The bulk of the work is carried out in Section \ref{subsect:SPGve}, where we show how to build a regular $*$-semigroup $\bS(P,\G,\ve)$ from a chained projection groupoid $(P,\G,\ve)$; see Definition \ref{defn:pr} and Theorem \ref{thm:SGve}.  After showing that this does indeed yield a functor in Section \ref{subsect:Sfunctor}, we establish some basic properties of the semigroups $\bS(P,\G,\ve)$ in Section \ref{subsect:PoP}.  These will be used in Section \ref{sect:iso} to show that $\bG$ and $\bS$ are mutually inverse isomorphisms between the categories $\RSS$ and $\CPG$.

\subsection[The regular $*$-semigroup associated to a chained projection groupoid]{\boldmath The regular $*$-semigroup associated to a chained projection groupoid}\label{subsect:SPGve}

%For the duration of this section, we fix a chained projection groupoid $(P,\G,\ve)$, as in Definition~\ref{defn:CoPG}, and we continue to write $\C=\C(P)$, and so on.  

Our aim in this section is to construct a regular $*$-semigroup $S=\bS(P,\G,\ve)$ from a chained projection groupoid $(P,\G,\ve)$.  The underlying set of~$S$ will simply be $\G$, and the involution of $S$ will simply be inversion in $\G$:
\[
a^* = a^{-1} \qquad\text{for $a\in\G$.}
\]
The definition of the product in $S$, which we will denote by $\pr$, is more involved.  Its inspiration is drawn from the properties of regular $*$-semigroups discussed in Remark \ref{rem:GS}.  Coherency (cf.~Definition \ref{defn:CoPG}) is not required for the definition of $\pr$, or indeed for the first series of results in this section.

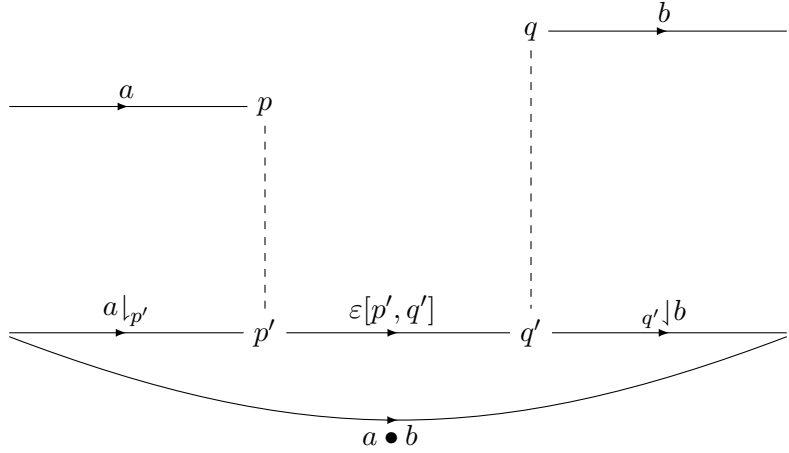
\begin{figure}[h]
\begin{center}
\begin{tikzpicture}[xscale=.7]
\tikzstyle{vertex}=[circle,draw=black, fill=white, inner sep = 0.07cm]
\node (p) at (5,3){$p$};
\node (q) at (10,4){$q$};
\node (p0) at (0,3){};
\node (q0) at (15,4){};
\node (p') at (5,0){$p'$};
\node (q') at (10,0){$q'$};
\node (p0') at (0,0){};
\node (q0') at (15,0){};
\draw[dashed] (p)--(p') 
(q)--(q') 
;
\draw[->-=0.5] (p0)--(p);
\draw[->-=0.5] (q)--(q0);
\draw[->-=0.5] (p0')--(p');
\draw[->-=0.5] (q')--(q0');
\draw[->-=0.5] (p')--(q');
\draw[->-=0.5] (p0') to [bend right = 15] (q0');
\node () at (2.4,3.2) {$a$};
\node () at (12.5,4.25) {$b$};
\node () at (2.4,.3) {$a\rest_{p'}$};
\node () at (12.5,.25) {${}_{q'}\corest b$};
\node () at (7.4,.3) {$\ve[p',q']$};
\node () at (7.35,-1.35) {$a\pr b$};
\end{tikzpicture}
\caption{Construction of the product $a\pr b$, as in Definition \ref{defn:pr}.}
\label{fig:pr}
\end{center}
\end{figure}

\begin{defn}\label{defn:pr}
Let $(P,\G,\ve)$ be a weak chained projection groupoid (cf.~Definition \ref{defn:ChPG}), and consider an arbitrary pair of morphisms~$a,b\in\G$.  Let $p=\br(a)$ and $q=\bd(b)$, and define the projections
\[
p' = q\th_p \AND q'=p\th_q.
\]
By Lemma \ref{lem:p'q'}, we have $p'\leq p$, $q'\leq q$ and $p'\F q'$.  In particular, the morphisms $ a\rest_{p'}$, $\ve[p',q']$ and ${}_{q'}\corest b$ exist, and we define
\begin{equation}\label{eq:pr}
a\pr b = a\rest_{p'} \circ \ve[p',q'] \circ {}_{q'}\corest b.
\end{equation}
This is illustrated in Figure \ref{fig:pr}.  We define $\bS(P,\G,\ve)$ to be the $(2,1)$-algebra with underlying set~$\G$, and with binary operation $\pr$ and unary operation ${}^*={}^{-1}$.
\end{defn}

Theorem \ref{thm:SGve} below says that $S=\bS(P,\G,\ve)$ is a regular $*$-semigroup when $(P,\G,\ve)$ belongs to $\CPG$, and not just to $\WCPG$.  To prove this, we must establish the identities:
\[
(a\pr b)\pr c = a\pr(b\pr c) \COMMA (a^*)^* = a = (a\pr a^*)\pr a \AND (a\pr b)^*=b^*\pr a^*.
\]
Associativity of $\pr$ is quite difficult to demonstrate, and is finally achieved in Lemma \ref{lem:abcabc} below.  The identities involving ${}^*$ are comparatively easy, however, and we verify these shortly in Lemma~\ref{lem:*}.  For its proof, we need Lemma \ref{lem:prcirc}, which shows that the (total) product $\pr$ extends the (partial) composition $\circ$, in the sense that these operations agree whenever the latter is defined.

For the rest of this section we fix a weak chained projection groupoid $(P,\G,\ve)$; we only assume coherency in Lemma \ref{lem:abcabc}, at the very last stage of the proof of associativity.

\begin{lemma}\label{lem:prcirc}
If $a,b\in\G$ are such that $\br(a)=\bd(b)$, then $a\pr b=a\circ b$.
\end{lemma}

\pf
In the notation of Definition \ref{defn:pr}, we have $p=q$, so also $p'=q'=p$.  But then
\[
a\pr b = a\rest_{p} \circ \ve[p,p] \circ {}_{p}\corest b = a\circ p\circ b = a\circ b.  \qedhere
\]
\epf

\begin{lemma}\label{lem:*}
For any $a,b\in\G$, we have
\[
(a^*)^* = a = (a\pr a^*)\pr a \AND (a\pr b)^*=b^*\pr a^*.
\]
\end{lemma}

\pf
Of course $(a^*)^*=(a^{-1})^{-1}=a$.  It follows from Lemma \ref{lem:prcirc} that $a\pr a^* = a\circ a^{-1} = \bd(a)$, and so $(a\pr a^*)\pr a = \bd(a)\pr a = \bd(a)\circ a = a$.

Next, let $p,q,p',q'$ be as in Definition \ref{defn:pr}.  Since $q=\br(b^*)$ and $p=\bd(a^*)$, we have
\[
a\pr b = a\rest_{p'} \circ \ve[p',q'] \circ {}_{q'}\corest b \AND b^*\pr a^* = b^*\rest_{q'} \circ \ve[q',p'] \circ {}_{p'}\corest a^* ,
\]
and so
\[
(a\pr b)^*= ({}_{q'}\corest b)^* \circ \ve[p',q']^* \circ (a\rest_{p'})^* = b^*\rest_{q'} \circ \ve[q',p'] \circ {}_{p'}\corest a^* = b^*\pr a^*.  \qedhere
\]
\epf

The symmetry/duality afforded by the identity $(a\pr b)^*=b^*\pr a^*$ will allow for some simplification in some of the proofs to follow.  This is the case with the proof of the next result, which identifies the domain and range of a product $a\pr b$.

\begin{lemma}\label{lem:drab}
For any $a,b\in\G$ we have
\[
\bd(a\pr b) = \bd(b)\Th_{a^{-1}} \AND \br(a\pr b) = \br(a)\Th_b.
\]
\end{lemma}

\pf
Let $p,q,p',q'$ be as in Definition \ref{defn:pr}, so that $a\pr b$ is as in \eqref{eq:pr}.  Then
\[
\br(a\pr b) = \br({}_{q'}\corest b) = q'\vt_b = p\th_q\vt_b = \br(a) \th_{\bd(b)}\vt_b = \br(a)\Th_b.
\]
The identity involving domains may be proved in similar fashion (using \eqref{eq:vta*}).  Alternatively, it follows from the range identity and duality:
\[
\bd(a\pr b) = \br((a\pr b)^*) = \br(b^*\pr a^*) = \br(b^*)\Th_{a^*} = \bd(b)\Th_{a^{-1}}.  \qedhere
\]
\epf

Our next result establishes an important identity, and can be thought of as a generalisation of Lemma \ref{lem:Thab} (in light of Lemma \ref{lem:prcirc}) in the special case of (weak) chained projection groupoids.

\begin{lemma}\label{lem:Thaprb}
For any $a,b\in\G$ we have $\Th_{a\pr b} = \Th_a\Th_b$.
\end{lemma}

\pf
Let $p,q,p',q'$ be as in Definition \ref{defn:pr}.  Then by the projection algebra axioms we have
\[
\th_{p'}\th_{q'} = \th_{q\th_p}\th_{p\th_q} =_4 \th_p\th_q\th_p\th_q\th_p\th_q =_5 \th_p\th_q. 
\]
We then calculate
\begin{align*}
\Th_{a\pr b} &= \Th_{a\rest_{p'}} \circ \Th_{\ve[p',q']} \circ \Th_{{}_{q'}\corest b} &&\text{by \eqref{eq:pr} and Lemma \ref{lem:Thab}}\\
&= \Th_a\th_{p'} \circ \th_{p'}\th_{q'} \circ \th_{q'}\Th_b &&\text{by \ref{G2c} and Lemmas \ref{lem:G2}\ref{G25} and \ref{lem:ve}\ref{ve2}}\\
&= \Th_a \circ \th_{p'}\th_{q'} \circ \Th_b &&\text{by \ref{P2}}\\
&= \Th_a \circ \th_p\th_q \circ \Th_b &&\text{by the above observation}\\
&= \Th_a \th_{\br(a)}\circ\th_{\bd(b)} \Th_b\\
&= \Th_a \Th_b &&\text{by Lemma \ref{lem:G2}\ref{G22}.}  \qedhere
\end{align*}
\epf

We are almost ready to show that $\pr$ is associative when $(P,\G,\ve)$ is coherent, i.e.~that $(a\pr b)\pr c$ and $a\pr(b\pr c)$ are equal, for all $a,b,c\in\G$.  Before we do so, we note that Lemmas \ref{lem:drab} and \ref{lem:Thaprb} can be used to show that these two terms have equal domains, and equal ranges, still without assuming coherency.  For example,
\begin{equation}\label{eq:rabc}
\br((a\pr b)\pr c) = \br(a\pr b)\Th_c = \br(a)\Th_b\Th_c = \br(a)\Th_{b\pr c} = \br(a\pr(b\pr c)).  
\end{equation}
A similar calculation gives
\begin{equation}\label{eq:dabc}
\bd((a\pr b)\pr c) = \bd(a\pr(b\pr c)) = \bd(c)\Th_{b^{-1}}\Th_{a^{-1}}.
\end{equation}

As a stepping stone to associativity, the next result gives expansions of $(a\pr b)\pr c$ and $a\pr(b\pr c)$, expressing each term as a composition of five morphisms.

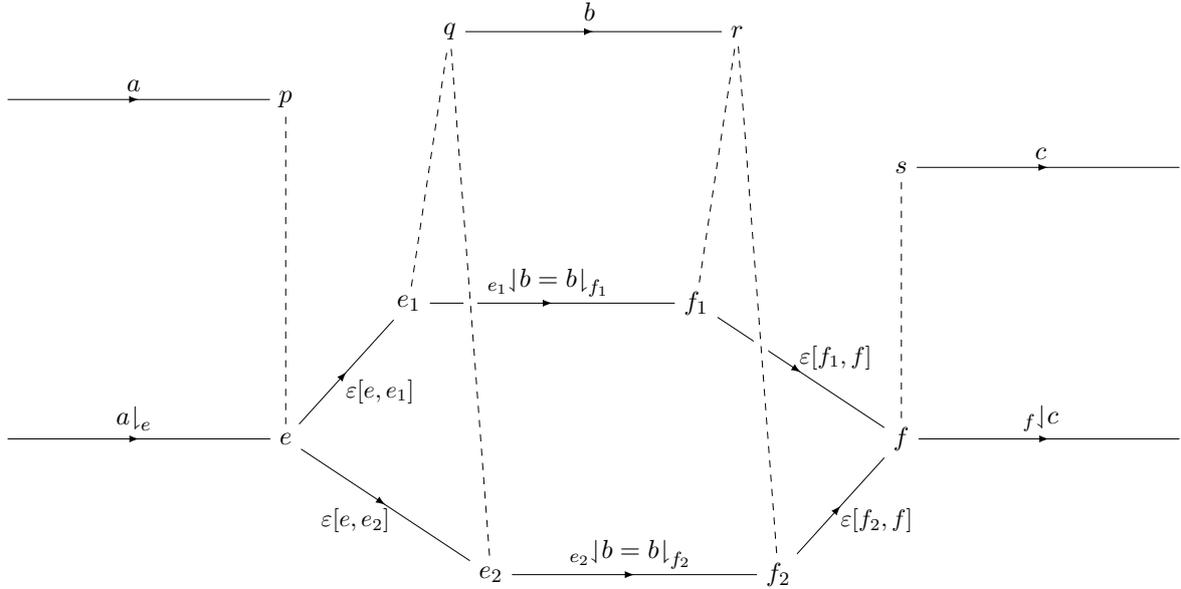
\begin{figure}[h]
\begin{center}
\scalebox{0.9}{
\begin{tikzpicture}[xscale=0.6]
\tikzstyle{vertex}=[circle,draw=black, fill=white, inner sep = 0.07cm]
\node (e) at (0,0){$e$};
\node (e1) at (3,2){$e_1$};
\node (e2) at (5,-2){$e_2$};
\node (f) at (15,0){$f$};
\node (f1) at (10,2){$f_1$};
\node (f2) at (12,-2){$f_2$};
\node (p) at (0,5){$p$};
\node (q) at (4,6){$q$};
\node (r) at (11,6){$r$};
\node (s) at (15,4){$s$};
\node (da) at (-7,5) {};
\node (rc) at (22,4) {};
\node (dae) at (-7,0) {};
\node (rfc) at (22,0) {};
\draw[->-=0.5] (q)--(r);
\draw[->-=0.5] (e)--(e1);
\draw[->-=0.5] (e1)--(f1);
\draw[->-=0.5] (f1)--(f);
\draw[->-=0.5] (e)--(e2);
\draw[->-=0.5] (e2)--(f2);
\draw[->-=0.5] (f2)--(f);
\draw[->-=0.5] (da)--(p);
\draw[->-=0.5] (s)--(rc);
\draw[->-=0.5] (dae)--(e);
\draw[->-=0.5] (f)--(rfc);
\draw[white,line width=2mm] (q)--(e2);
\draw[white,line width=2mm] (r)--(f2);
\draw[dashed] (e1)--(q)--(e2) (f1)--(r)--(f2) (p)--(e) (s)--(f);
\node () at (-3.7,5.2) {$a$};
\node () at (7.4,6.3) {$b$};
\node () at (18.4,4.2) {$c$};
\node () at (6.4,2.3) {${}_{e_1}\corest b = b\rest_{f_1}$};
\node () at (8.4,-1.7) {${}_{e_2}\corest b = b\rest_{f_2}$};
\node () at (2.3,0.7) {{\small $\ve[e,e_1]$}};
\node () at (1.7,-1.2) {{\small $\ve[e,e_2]$}};
\node () at (13.4,1.2) {{\small $\ve[f_1,f]$}};
\node () at (14.4,-1.2) {{\small $\ve[f_2,f]$}};
\node () at (-3.7,.3) {$a\rest_e$};
\node () at (18.4,.3) {${}_f\corest c$};
\end{tikzpicture}
}
\caption{Top:  three morphisms $a,b,c\in\G$.  Bottom: the upper and lower paths represent $(a\pr b)\pr c$ and $a\pr(b\pr c)$, respectively; cf.~Lemma \ref{lem:abc}.  Dashed lines indicate $\leq$ relationships.}
\label{fig:abc}
\end{center}
\end{figure}

\begin{lemma}\label{lem:abc}
If $a,b,c\in\G$, and if $p=\br(a)$, $q=\bd(b)$, $r=\br(b)$ and $s=\bd(c)$, then 
\ben
\item \label{abc2} $(a\pr b)\pr c = a\rest_e \circ \ve[e,e_1] \circ {}_{e_1}\corest b \circ \ve[f_1,f] \circ {}_f\corest c$,
where
\[
e = s\Th_{b^{-1}}\th_p \COMMA e_1 = e\th_q \COMMA f_1 = e\Th_b \AND f = p\Th_b\th_s,
\]
\item \label{abc1} $a\pr (b\pr c) = a\rest_e \circ \ve[e,e_2] \circ b\rest_{f_2} \circ \ve[f_2,f] \circ {}_f\corest c$,
where
\[
e = s\Th_{b^{-1}}\th_p \COMMA e_2 = f\Th_{b^{-1}} \COMMA f_2 = f\th_r \AND f = p\Th_b\th_s.
\]
\een
The above factorisations are illustrated in Figure \ref{fig:abc}.
\end{lemma}

\pf
\firstpfitem{\ref{abc2}}  Put $p'=q\th_p$ and $q'=p\th_q$, so that 
\[
a\pr b = a'\circ\ve'\circ b' \WHERE a' = a\rest_{p'} \COMMa \ve' = \ve[p',q'] \ANd b' = {}_{q'}\corest b.
\]
Keeping Lemma \ref{lem:drab} in mind, we now define
\[
t = \br(a\pr b) = p\Th_b.
\]
Then
\[
(a\pr b)\pr c = (a\pr b)\rest_{t'} \circ \ve[t',s'] \circ {}_{s'}\corest c \WHERE \text{$t'=s\th_t$ and $s'=t\th_s$.}
\]
We now focus on the term $(a\pr b)\rest_{t'}$.  By \eqref{eq:pab} we have
\[
(a\pr b)\rest_{t'} = (a'\circ\ve'\circ b')\rest_{t'} = a'\rest_u \circ \ve'\rest_v \circ b'\rest_{t'} \WHERE v = t'\vt_{b'}^{-1} \ANd u = t'\vt_{b'}^{-1}\vt_{\ve'}^{-1}.
\]
%\newpage\noindent
Now, 
\bit
\item $a'\rest_u = a\rest_{p'}\rest_u = a\rest_u$,
\item $\ve'\rest_v \leq \ve' = \ve[p',q']$, $\bd(\ve'\rest_v) = \br(a'\rest_u) = u$ and $\br(\ve'\rest_v)=v$, so $\ve'\rest_v = \ve[u,v]$ by Lemma \ref{lem:aepq}, and
\item $b'\rest_{t'}\leq b'\leq b$ and $\bd(b'\rest_{t'})=\br(\ve'\rest_v)=v$, so $b'\rest_{t'}={}_v\corest b$.  (For later use, since $b'\rest_{t'}$ and $b\rest_{t'}$ are both below $b$ and have codomain $t'$, it also follows that $b\rest_{t'}=b'\rest_{t'}={}_v\corest b$.)
\eit
Putting everything together, it follows that
\begin{align*}
(a\pr b)\pr c &= (a\pr b)\rest_{t'} \circ \ve[t',s'] \circ {}_{s'}\corest c \\
&= a'\rest_u \circ \ve'\rest_v \circ b'\rest_{t'} \circ \ve[t',s'] \circ {}_{s'}\corest c \\
&= a\rest_u \circ \ve[u,v] \circ {}_v\corest b \circ \ve[t',s'] \circ {}_{s'}\corest c .
\end{align*}
It therefore remains to check that
\bena
\bmc4
\item \label{dabc1} $u=e$, 
\item \label{dabc2} $v=e_1$, 
\item \label{dabc3} $t'=f_1$, 
\item \label{dabc4} $s'=f$,
\emc
\een
where $e,e_1,f_1,f$ are as defined in the statement of the lemma.  

\pfitem{\ref{dabc1}}  First note that $\bd(a\rest_u) = \bd((a\pr b)\pr c) = s\Th_{b^{-1}}\Th_{a^{-1}}$ by \eqref{eq:dabc}, so $a\rest_u = {}_{s\Th_{b^{-1}}\Th_{a^{-1}}}\corest a$.  But then
\[
u = \br(a\rest_u) = \br({}_{s\Th_{b^{-1}}\Th_{a^{-1}}}\corest a) = s\Th_{b^{-1}}\Th_{a^{-1}}\vt_a = s\Th_{b^{-1}}\th_p = e,
\]
where we used Lemma \ref{lem:G2}\ref{G23} in the second-last step.

\pfitem{\ref{dabc3}}  Here we use \ref{G2b} to calculate $t' = s\th_t = s\th_{p\Th_b} = s\Th_{b^{-1}}\th_p\Th_b = e\Th_b = f_1$.

\pfitem{\ref{dabc2}}  We noted above that ${}_v\corest b=b\rest_{t'}$.  Combining this with \eqref{eq:vta*}, the just-proved item \ref{dabc3}, and Lemma \ref{lem:G2}\ref{G23}, it follows that
\[
v = \bd({}_v\corest b) = \bd(b\rest_{t'}) = t'\vt_{b^{-1}} = f_1\vt_{b^{-1}} = (e\Th_b)\vt_{b^{-1}} = e\th_q = e_1.
\]

\pfitem{\ref{dabc4}}  Finally, $s'=t\th_s=p\Th_b\th_s=f$.

\aftercases  As noted above, this completes the proof of \ref{abc2}.

\pfitem{\ref{abc1}}  This can be proved in similar fashion to part \ref{abc2}, but in fact we can also deduce it from \ref{abc2}.  Indeed, we first use Lemma \ref{lem:*} to expand
\[
(a\pr(b\pr c))^{-1} = (b\pr c)^{-1}\pr a^{-1} = (c^{-1}\pr b^{-1})\pr a^{-1}.
\]
Note that $s=\br(c^{-1})$, $r=\bd(b^{-1})$, $q=\br(b^{-1})$ and $p=\bd(a^{-1})$.  It follows from part \ref{abc2} that 
\[
(c^{-1}\pr b^{-1})\pr a^{-1} = c^{-1}\rest_{e'} \circ \ve[e',e_1'] \circ {}_{e_1'}\corest b^{-1} \circ \ve[f_1',f'] \circ {}_{f'}\corest a^{-1},
\]
where
\[
e' = p\Th_b\th_s \COMMA e_1' = e'\th_r \COMMA f_1' = e'\Th_{b^{-1}} \AND f' = s\Th_{b^{-1}}\th_p.
\]
But then
\begin{align}
\nonumber a\pr (b\pr c) &= ((c^{-1}\pr b^{-1})\pr a^{-1})^{-1} \\
\nonumber &= (c^{-1}\rest_{e'} \circ \ve[e',e_1'] \circ {}_{e_1'}\corest b^{-1} \circ \ve[f_1',f'] \circ {}_{f'}\corest a^{-1})^{-1} \\
\nonumber &= ({}_{f'}\corest a^{-1})^{-1} \circ \ve[f_1',f']^{-1} \circ ({}_{e_1'}\corest b^{-1})^{-1} \circ \ve[e',e_1']^{-1} \circ (c^{-1}\rest_{e'})^{-1} \\
\label{eq:a(bc)1} &= a\rest_{f'} \circ \ve[f',f_1'] \circ b\rest_{e_1'} \circ \ve[e_1',e'] \circ {}_{e'}\corest c .
\end{align}
We note immediately that $f'=e$ and $e'=f$ (where $e,f$ are as in the statement of the lemma).  We also have
\[
e_1' = e'\th_r = f\th_r = f_2 \AND f_1' = e'\Th_{b^{-1}} = f\Th_{b^{-1}} = e_2.
\]
Thus, continuing from \eqref{eq:a(bc)1}, we have
\[
a\pr (b\pr c) = a\rest_e \circ \ve[e,e_2] \circ b\rest_{f_2} \circ \ve[f_2,f] \circ {}_f\corest c,
\]
as required.
\epf

\begin{lemma}\label{lem:abcabc}
If $(P,\G,\ve)$ is a chained projection groupoid, then $(a\pr b)\pr c = a\pr(b\pr c)$ for any $a,b,c\in\G$.
\end{lemma}

\pf
By Lemma \ref{lem:abc}, it suffices to show (in the notation of that lemma) that
\begin{equation}\label{eq:abcabc}
\ve[e,e_1] \circ {}_{e_1}\corest b \circ \ve[f_1,f] = \ve[e,e_2] \circ b\rest_{f_2} \circ \ve[f_2,f].
\end{equation}
Comparing the definitions of the $e_i,f_i$ from Lemma \ref{lem:abc} (in terms of $e,f,q,r,b$) with \eqref{eq:e1e2f1f2}, an application of \ref{G3} will give us \eqref{eq:abcabc} if we can show that $(e,f)$ is $b$-linked (cf.~Definitions \ref{defn:LP1} and \ref{defn:LP2}).  In other words, the proof of the lemma will be complete if we can show that
\[
f = e\Th_b\th_f \AND e = f\Th_{b^{-1}}\th_e.
\]
For the first, we have
\begin{align*}
e\Th_b\th_f &= s\Th_{b^{-1}}\th_p\Th_b\th_{p\Th_b\th_s} &&\text{by definition}\\
&= s\Th_{b^{-1}}\th_p\Th_b\th_s\th_{p\Th_b}\th_s &&\text{by \ref{P4}}\\
&= p\Th_b\th_s\Th_{b^{-1}}\th_{p}\Th_b\th_s &&\text{by \ref{G2b} and Lemma \ref{lem:pa'qap}}\\
&= s\Th_{b^{-1}}\th_{p}\Th_b\th_s &&\text{by Lemma \ref{lem:pa'qap}}\\
&= p\Th_b\th_s = f &&\text{by Lemma \ref{lem:pa'qap}.}
\end{align*}
The proof that $e = f\Th_{b^{-1}}\th_e$ is virtually identical.
\epf

Lemmas \ref{lem:*} and \ref{lem:abcabc} give the following:

\begin{thm}\label{thm:SGve}
If $(P,\G,\ve)$ is a chained projection groupoid, then $\bS(P,\G,\ve)$ is a regular $*$-semigroup.  \epfres
\end{thm}

\subsection{A functor}\label{subsect:Sfunctor}

We have just seen how to construct a regular $*$-semigroup $\bS(P,\G,\ve)$ from a chained projection groupoid $(P,\G,\ve)$.  As ever, this construction can be thought of as an object map from the category $\CPG$ to the category $\RSS$.  Since a morphism $(P,\G,\ve)\to(P',\G',\ve')$ in $\CPG$ is a (special kind of) map $\phi:\G\to\G'$, and since the semigroups $S=\bS(P,\G,\ve)$ and $S'=\bS(P',\G',\ve')$ have the same underlying sets as~$\G$ and~$\G'$ (respectively), we can think of $\phi$ as a map ${\bS(\phi) = \phi : S\to S'}$.  We now show that $\bS$ (interpreted this way) is a functor:

\begin{thm}\label{thm:Sfunctor}
$\bS$ is a functor $\CPG\to\RSS$.
\end{thm}

\pf
It remains to show that any morphism $\phi:(P,\G,\ve)\to(P',\G',\ve')$ in $\CPG$ is a $*$-morphism $S\to S'$, where $S=\bS(P,\G,\ve)$ and $S'=\bS(P',\G',\ve')$.

Since $\phi$ is a groupoid functor, we have
\[
(a^*)\phi = (a^{-1})\phi = (a\phi)^{-1} = (a\phi)^* \qquad\text{for all $a\in S$.}
\]
To make the rest of this proof notationally concise, we will write $\ol a = a\phi$ for $a\in S(=\G)$.  It remains to show that $\phi$ is a semigroup homomorphism, i.e.~that $\ol{a\pr b}=\ol a\pr \ol b$ for all $a,b\in S$.  (Here we write~$\pr$ for the product on both $S$ and $S'$.)  For this, write as usual
\[
p = \br(a) \COMMA q = \bd(b) \COMMA p' = q\th_p \AND q' = p\th_q,
\]
so that
\[
a\pr b = a\rest_{p'} \circ \ve[p',q'] \circ {}_{q'}\corest b.
\]
Since $\phi$ is a functor, we have
\[
\br(\ol a) = \br(\ol{a\circ p}) = \br(\ol a\circ \ol p)  = \br(\ol p) = \ol p  \ANDSIM \bd(\ol b) = \ol q.
\]
Thus,
\[
\ol a \pr \ol b = \ol a\rest_{\ol p'} \circ \ve'[\ol p',\ol q'] \circ {}_{\ol q'}\corest \ol b \WHERE \ol p' = \ol q\th'_{\ol p} \AND \ol q' = \ol p\th'_{\ol q}.
\]
Since $\phi$ is a projection groupoid morphism, its object map $v\phi:P\to P'$ is a projection algebra morphism, and so
%By \ref{F1} we have
\[
\ol{p'} = \ol{q\th_p} = \ol q\th'_{\ol p} = \ol p' \ANDSIM \ol{q'} = \ol q'.
\]
Since $\phi$ respects the evaluation maps, we have
\[
\ol{\ve[p',q']} = (\ve[p',q'])\phi = \ve'([p',q']\Phi) = \ve'[\ol{p'},\ol{q'}] = \ve'[\ol p',\ol q'].
\]
Since $\phi$ is an ordered functor, we have
\[
\ol{a\rest_{p'}} = \ol a\rest_{\ol{p'}} = \ol a\rest_{\ol p'} \ANDSIM \ol{{}_{q'}\corest b} = {}_{\ol q'}\corest\ol b.
\]
Putting everything together, we have
\[
\ol{a\pr b} = \ol{a\rest_{p'} \circ \ve[p',q'] \circ {}_{q'}\corest b} = \ol{a\rest_{p'}} \circ \ol{\ve[p',q']} \circ \ol{{}_{q'}\corest b} = \ol a\rest_{\ol p'} \circ \ve'[\ol p',\ol q'] \circ {}_{\ol q'}\corest \ol b = \ol a\pr\ol b.  \qedhere
\]
\epf

\subsection{Products of projections}\label{subsect:PoP}

For the duration of this section we fix a chained projection groupoid $(P,\G,\ve)$.  The next two results concern $\pr$ products involving projections in the regular $*$-semigroup $S=\bS(P,\G,\ve)$.  These will be of use later, but we also apply them here to describe the projections, idempotents and idempotent-generated subsemigroup of $S$ in Proposition \ref{prop:ES}.

\begin{lemma}\label{lem:apr}
If $a\in\G(p,q)$, and if $t\in P$, then
\ben
\item \label{apr1} $a\pr t = a\rest_{q'}\circ\ve[q',t']$, where $q'=t\th_q$ and $t'=q\th_t$,
\item \label{apr2} $t\pr a = \ve[t',p']\circ{}_{p'}\corest a$, where $p'=t\th_p$ and $t'=p\th_t$,
\item \label{apr3} $a\pr t = a\circ\ve[q,t]$ if $q\F t$,
\item \label{apr4} $t\pr a = \ve[t,p]\circ a$ if $p\F t$,
\item \label{apr5} $a\pr t = a\rest_t$ if $t\leq q$,
\item \label{apr6} $t\pr a = {}_t\corest a$ if $t\leq p$.
\een
\end{lemma}

\pf
By symmetry, it suffices to prove \ref{apr1}, \ref{apr3} and \ref{apr5}.

\pfitem{\ref{apr1}}  We have $a\pr t = a\rest_{q'}\circ\ve[q',t']\circ{}_{t'}\corest t = a\rest_{q'}\circ\ve[q',t']\circ t' = a\rest_{q'}\circ\ve[q',t']$.

\pfitem{\ref{apr3}}  This follows from part \ref{apr1}, as $q'=q$ and $t'=t$ if $q\F t$.

\pfitem{\ref{apr5}}  If $t\leq q$, then we also have $t\leqF q$ by Lemma \ref{lem:pqp}\ref{pqP2}, so that $t=t\th_q=q'$ and $t=q\th_t=t'$.  It then follows from part \ref{apr1} that
\[
a\pr t = a\rest_{q'}\circ\ve[q',t'] = a\rest_t\circ\ve[t,t] = a\rest_t\circ t = a\rest_t.  \qedhere
\]
\epf

\begin{lemma}\label{lem:pqvepq}
If $p,q\in P$, then
\ben
\item \label{pqvepq1} $p\pr q = \ve[p',q']$, where $p'=q\th_p$ and $q'=p\th_q$,
\item \label{pqvepq2} $p\pr q = \ve[p,q]$ if $p\F q$.
\item \label{pqvepq3} $p\pr q\pr p = q\th_p$.
\een
\end{lemma}

\pf
\firstpfitem{\ref{pqvepq1}}  By Lemma \ref{lem:apr}\ref{apr1} we have $p\pr q = p\rest_{p'}\circ\ve[p',q'] = p'\circ\ve[p',q'] = \ve[p',q']$.

\pfitem{\ref{pqvepq2}} This follows from part \ref{pqvepq1}, as $p=p'$ and $q=q'$ when $p\F q$.

\pfitem{\ref{pqvepq3}}  Here we have
\begin{align*}
p\pr q\pr p &= \ve[p',q']\pr p &&\text{by part \ref{pqvepq1}, where $p'=q\th_p$ and $q'=p\th_q$}\\
&= \ve[p',q']\rest_{q''}\circ\ve[q'',p''] &&\text{by Lemma \ref{lem:apr}\ref{apr1}, where $q''=p\th_{q'}$ and $p''=q'\th_p$.}  
\end{align*}
Using the projection algebra axioms, it is easy to show that in fact $p''=p'$ and $q''=q'$.  Thus, continuing from above, we have
\[
p\pr q\pr p
= \ve[p',q']\rest_{q''} \circ \ve[q'',p'']  
= \ve[p',q']\rest_{q'} \circ \ve[q',p'] 
= \ve[p',q'] \circ \ve[q',p']
= p' = q\th_p.  \qedhere
\]
\epf

Recall that the sets of idempotents and projections of a regular $*$-semigroup $S$ are denoted
\[
E(S) = \set{e\in S}{e^2=e} \AND P(S) = \set{p\in S}{p^2=p=p^*}.
\]
Although these sets are not equal in general, Lemma \ref{lem:PS1}\ref{PS12} says that they generate the same subsemigroup of $S$, i.e.~$\la E(S)\ra = \la P(S)\ra$.  We write $\E(S)$ for this idempotent-generated subsemigroup.  The results proved above allow us to give the following descriptions of $P(S)$, $E(S)$ and~$\E(S)$ in the case that $S=\bS(P,\G,\ve)$ arises from a chained projection groupoid $(P,\G,\ve)$.  %(It follows from Proposition \ref{prop:SoG} that every regular $*$-semigroup has this form.)

\begin{prop}\label{prop:ES}
Let $(P,\G,\ve)$ be a chained projection groupoid, and let $S=\bS(P,\G,\ve)$.  Then
\ben
\item \label{ES1} $P(S)=P$, 
\item \label{ES2} $E(S)=\set{\ve[p,q]}{(p,q)\in{\F}}$,
\item \label{ES3} $\E(S) = \im(\ve) = \set{\ve(\c)}{\c\in\C(P)}$, and consequently $S$ is idempotent-generated if and only if~$\ve$ is surjective.
\een
\end{prop}

\pf
\firstpfitem{\ref{ES1}}  By Lemmas \ref{lem:PS1}\ref{PS11} and \ref{lem:prcirc} we have
\[
P(S) = \set{a\pr a^*}{a\in S} = \set{a\circ a^{-1}}{a\in \G} = \set{\bd(a)}{a\in \G} = v\G = P.
\]

\pfitem{\ref{ES2}}  By Lemma \ref{lem:PS1}\ref{PS12}, and part \ref{ES1}, we have $E(S) = \set{p\pr q}{p,q\in P}$.  The claim then follows from Lemma \ref{lem:pqvepq}.

\pfitem{\ref{ES3}}  If $\c=[p_1,\ldots,p_k]\in\C(P)$, then 
\begin{align}
\nonumber \ve(\c) &= \ve[p_1,p_2]\circ\ve[p_2,p_3]\circ\cdots\circ\ve[p_{k-1},p_k] &&\text{by \eqref{eq:vec}}\\
\nonumber &= \ve[p_1,p_2]\pr\ve[p_2,p_3]\pr\cdots\pr\ve[p_{k-1},p_k] &&\text{by Lemma \ref{lem:prcirc}}\\
\nonumber &= (p_1\pr p_2)\pr(p_2\pr p_3)\pr\cdots\pr(p_{k-1}\pr p_k) &&\text{by Lemma \ref{lem:pqvepq}\ref{pqvepq2}}\\
\label{eq:ES} &= p_1\pr p_2\pr\cdots\pr p_k \in \E(S).
\end{align}

Conversely, fix some $a\in\E(S)$, so that $a = p_1\pr\cdots\pr p_k$ for some $p_1,\ldots,p_k\in P$.  We show that $a\in\im(\ve)$ by induction on $k$.  The $k=1$ case being trivial, we assume that $k\geq2$, and let $b=p_1\pr\cdots\pr p_{k-1}$.  By induction we have $b=\ve(\c)$ for some $\c\in\C(P)$, and we write $r=\br(b)=\br(\c)$.  Then with $r'=p_k\th_r$ and $p_k'=r\th_{p_k}$, it follows from Lemma \ref{lem:apr}\ref{apr1} and properties of $\ve$ from Definition \ref{defn:ve} that
\[
a = b\pr p_k = b\rest_{r'}\circ\ve[r',p_k'] = \ve(\c)\rest_{r'}\circ\ve[r',p_k'] = \ve(\c\rest_{r'})\circ\ve[r',p_k'] = \ve(\c\rest_{r'}\circ[r',p_k']) \in \im(\ve).  \qedhere
\]
\epf

\begin{rem}\label{rem:ES}
We will see in Section \ref{sect:iso} (see Proposition \ref{prop:SoG}) that every regular $*$-semigroup has the form $\bS(P,\G,\ve)$ for some chained projection groupoid $(P,\G,\ve)$.  Combining this with Proposition~\ref{prop:ES}\ref{ES3}, we obtain a proof of Proposition \ref{prop:ERSS}, stated in Section \ref{subsect:Sprelim}.
\end{rem}

\section{The category isomorphism}\label{sect:iso}

We are now ready to bring all of the above ideas together.  By Theorems \ref{thm:Gfunctor} and \ref{thm:Sfunctor}, we have two functors
\[
\bG:\RSS\to\CPG \AND \bS:\CPG\to\RSS,
\]
and we now wish to show that these are mutual inverses, thereby proving the following, which is the main result of the paper:

\begin{thm}\label{thm:iso}
The category $\RSS$ of regular $*$-semigroups (with $*$-semigroup homomorphisms) is isomorphic to the category $\CPG$ of chained projection groupoids (with chained projection functors).
\end{thm}

\pf
The functors $\bG$ and $\bS$ both act identically on morphisms, which in each category are certain structure-preserving mappings on objects (regular $*$-semigroups or chained projection groupoids).  It follows that $\bG$ and $\bS$ are mutually inverse at the level of morphisms.  For objects, it remains to check that
\[
\bG(\bS(P,\G,\ve)) = (P,\G,\ve) \AND \bS(\bG(S)) = S
\]
for any chained projection groupoid $(P,\G,\ve)$, and any regular $*$-semigroup $S$.  This is achieved in Propositions \ref{prop:GoS} and \ref{prop:SoG} below.
\epf

\begin{prop}\label{prop:GoS}
If $(P,\G,\ve)$ is a chained projection groupoid, then $\bG(\bS(P,\G,\ve))=(P,\G,\ve)$.
\end{prop}

\pf
Write
\[
S = \bS(P,\G,\ve) \AND (P',\G',\ve') = \bG(S) = (P(S),\G(S),\ve(S)).
\]
We must show that:
\ben
\item \label{GoS1} $P'=P$, meaning that the two projection algebras have the same underlying sets, and the same $\th$ operations,
\item \label{GoS2} $\G'=\G$, meaning that the two groupoids have the same underlying (object and morphism) sets, and all the same operations and relations (domain/range maps, composition, inversion and order), and
\item \label{GoS3} $\ve'=\ve$, meaning simply that the two maps are identical.
\een
\pfitem{\ref{GoS1}}  This follows quickly from Proposition~\ref{prop:ES}\ref{ES1} and Lemma \ref{lem:pqvepq}\ref{pqvepq3}.

\pfitem{\ref{GoS2}}  By construction, the underlying sets of $S$ and $\G'=\G(S)$ are simply $\G$, and the unary operations of $S$ and~$\G'$ are simply inversion in $\G$.  The projection set of $S$ is $P$ (by Proposition~\ref{prop:ES}\ref{ES1}), and the object set of $\G'=\G(S)$ is $P(S)=P$ (by construction).

Recall that the product $\pr$ in $S$ is given by Definition \ref{defn:pr}.  The groupoid $\G'$ is constructed from $S$, as outlined in Definition \ref{defn:GS}.  To avoid ambiguity, we will denote the composition on~$\G'$ by~$\ccirc$ in order to distinguish it from the original composition $\circ$ in $\G$.  For the same reason, we also denote the domain and range maps in $\G'$ by $\bd'$ and $\br'$.  So for $a,b\in\G'(=\G)$, we have
\[
\bd'(a) = a\pr a^{-1} \COMMA \br'(a) = a^{-1}\pr a \AND a\ccirc b=a\pr b \text{ when $\br'(a)=\bd'(b)$.}
\]
By Lemma \ref{lem:prcirc}, and since $\br(a)=\bd(a^{-1})$, it follows that
\begin{equation}\label{eq:d'r'}
\bd'(a) = a\pr a^{-1} = a \circ a^{-1} = \bd(a) \ANDSIM \br'(a) = \br(a).
\end{equation}
It also follows that
\[
\br'(a)=\bd'(b) \Iff \br(a)=\bd(b) \qquad\text{and then}\qquad a\ccirc b = a\pr b = a\circ b.
\]
Thus, the compositions in the categories $\G$ and $\G'$ coincide.

We now move on to the orders on $\G$ and $\G'$.  We use the symbols $\corest$ and $\rest$ to denote restrictions in $\G$ as usual, and $\corestt$ and $\restt$ for restrictions in $\G'$ (note the placement of the arrow heads).  Let $a\in\G$ and write $q=\bd(a)$.  By \eqref{eq:d'r'}, ${}_p\corest a$ is defined in $\G$ precisely when ${}_p\corestt a$ is defined in $\G'$, i.e.~when $p\leq q$.  And for such $p$, it follows from Definition \ref{defn:GS} and Lemma \ref{lem:apr}\ref{apr6} that
\[
{}_p\corestt a = p\pr a = {}_p\corest a .
\]
But then for $a,b\in\G$, and writing $\leq$ and $\leq'$ for the orders in $\G$ and $\G'$, we have
\[
a\leq b \Iff a={}_p\corest b\text{ for some $p\leq\bd(b)$} \Iff a={}_p\corestt b\text{ for some $p\leq\bd'(b)$} \Iff a\leq' b.
\]

\pfitem{\ref{GoS3}}  If $\c=[p_1,\ldots,p_k]\in\C$, then Definition \ref{defn:veS} and \eqref{eq:ES} give $\ve'(\c) = p_1\pr\cdots\pr p_k = \ve(\c)$.
%and we show that $\ve(\c)=\ve'(\c)$ by induction on $k$.  The $k=1$ case being clear, suppose $k\geq2$, and write $\d=[p_1,\ldots,p_{k-1}]$ and $a=\ve(\d)$.  Note that $\c=\d\circ[p_{k-1},p_k]$, and by induction we have
%\[
%a = \ve(\d) = \ve'(\d) = p_1\pr\cdots\pr p_{k-1}.
%\]
%It follows that $\ve'(\c)=a\pr p_k$.  Since $\br(a)=\br(\d)=p_{k-1}$ (as $\ve$ is a $v$-functor), and since $p_{k-1}\F p_k$ (as $\c\in\C$), it follows from Lemma \ref{lem:apr}\ref{apr3} that
%\[
%\ve'(\c) = a\pr p_k = a\circ\ve[p_{k-1},p_k] = \ve(\d)\circ\ve[p_{k-1},p_k] = \ve(\d\circ[p_{k-1},p_k]) = \ve(\c),
%\]
%as required.
\epf

\begin{prop}\label{prop:SoG}
If $S$ is a regular $*$-semigroup, then $\bS(\bG(S))=S$.
\end{prop}

\pf
Write
\[
(P,\G,\ve) = \bG(S) = (P(S),\G(S),\ve(S)) \AND S' = \bS(P,\G,\ve).
\]
Again, by construction the underlying sets of $S$ and $S'$ coincide, and so do the involutions.  It therefore remains to show that the products in $S$ and $S'$ coincide as well.  We denote the product in $S$ by juxtaposition, as usual.  The product in $S'$ is $\pr$, which is constructed from $(P,\G,\ve)$ as in Definition \ref{defn:pr}, while $(P,\G,\ve)$ is in turn constructed from $S$ as in Definition \ref{defn:PGveS}.

So let $a,b\in S$; we must show that $a\pr b=ab$.  Write $p=\br(a)$ and $q=\bd(b)$, so that
\[
a\pr b = a\rest_{p'} \circ \ve[p',q'] \circ {}_{q'}\corest b \WHERE p' = q\th_p \ANd q' = p\th_q.
\]
Following Definition \ref{defn:GS}, and using Lemma \ref{lem:PS1}, we have
\begin{align*}
p &= a^*a, & p' &= pqp, & a\rest_{p'} &= ap' = apqp, \\
q &= bb^*, & q' &= qpq, & {}_{q'}\corest b &= q'b = qpqb, & \ve[p',q'] = p'q' = pqpqpq = pq.
\end{align*}
Combining all of the above, and keeping Definition \ref{defn:GS} and Lemma \ref{lem:PS1} in mind, it follows that
\[
a\pr b = apqp \circ pq \circ qpqb = apqpqpqb = apqb = aa^*abb^*b = ab,
\]
as required.
\epf

As noted above, this completes the proof of Theorem \ref{thm:iso}.

\begin{rem}\label{rem:structure}
We can use Theorem \ref{thm:iso} to deduce a \emph{structure theorem} for regular $*$-semigroups, in the sense that it shows how the entire structure of such a semigroup can be entirely, and uniquely, determined by `simpler' objects.  Specifically, given a regular $*$-semigroup $S$, one constructs:
\bit
\item the projection algebra $P=P(S)$, as in Definition~\ref{defn:PS},
\item the ordered groupoid $\G=\G(S)$, as in Definition \ref{defn:GS}, and
\item the evaluation map $\ve=\ve(S)$, as in Definition \ref{defn:veS}.
\eit
Proposition \ref{prop:SoG} tells us that $S$ is completely determined by these three simpler objects, in the sense that $S=\bS(P,\G,\ve)$ can be constructed from $P$, $\G$ and $\ve$ in the manner described in Definition~\ref{defn:pr}.  Examples \ref{eg:AG} and \ref{eg:Rees} demonstrate the subtlety of the decomposition, showing how small changes to $P$, $\G$ or $\ve$ can result in very different semigroups.
\end{rem}

\begin{rem}
Theorem \ref{thm:iso} shows that the categories of regular $*$-semigroups and chained projection groupoids are in a sense the same.  It also follows from the definition of the functors~$\bS$ and~$\bG$ that the algebras of projections of regular $*$-semigroups are precisely the object sets of chained projection groupoids.  It follows from Proposition \ref{prop:PtoS} (or indeed from results of \cite{Imaoka1983,Jones2012,Paper2,Paper3}) that this is precisely the class of all (abstract) projection algebras.
\end{rem}

\begin{rem}
It is again worth considering an example.  Let $S=\M^0(P,P,G,M)$ be a Rees 0-matrix regular $*$-semigroup, as in Example \ref{eg:Rees}, and let $(P_0,\G,\ve)=\bG(S)$ be the associated chained projection groupoid.  The product $a\pr b$ of two non-zero elements $a=(p,g,q)$ and $b=(r,h,s)$ of $S$ is given by
\[
a\pr b = a\rest_{q'} \circ \ve[q',r'] \circ {}_{r'}\corest b,
\WHERE
q' = r\th_q \AND r' = q\th_r.
\]
By \eqref{eq:M0th}, we have
\[
q' = \begin{cases}
q &\text{if $m_{qr}\not=0$}\\
0 &\text{if $m_{qr}=0$}
\end{cases}
\AND
r' = \begin{cases}
r &\text{if $m_{qr}\not=0$}\\
0 &\text{if $m_{qr}=0$.}
\end{cases}
\]
When $m_{qr} = 0$, we have $a\pr b = a\rest_0\circ\ve[0,0]\circ{}_0\corest b = 0$.  Otherwise, we use \eqref{eq:Rees2} to calculate
\[
a\pr b = a\rest_q \circ \ve[q,r] \circ {}_r\corest b = (p,g,q)\circ(q,m_{qr},r)\circ(r,g,s) = (p,gm_{qr}h,s).
\]
Thus, in both cases we have $a\pr b=ab$, as expected.
\end{rem}

\section{Inverse semigroups and inductive groupoids}\label{sect:I}

We have noted on a number of occasions that any inverse semigroup is a regular $*$-semigroup (with $a^*=a^{-1}$).  In this final section, we look at how the theory developed above simplifies in the case of inverse semigroups.  In particular, we will see that our results allow us to deduce the celebrated Ehresmann--Schein--Nambooripad (ESN) Theorem, stated below as Theorem \ref{thm:ESN}, and discussed in more detail in Section \ref{sect:intro}.  
%This theorem was first explicitly formulated by Lawson in \cite[Theorem 4.1.8]{Lawson1998}, who named it as such in order to honour the contributions of the three stated mathematicians to the development of the result.  %We discussed this in some detail in Section \ref{sect:intro}, but see 
%See \cite[Chapter 4]{Lawson1998} and \cite{Hollings2012,HL2017} for more details.

In what follows, we write $\IS$ for the category of inverse semigroups.  Morphisms in $\IS$ are simply the semigroup homomorphisms.  (Any semigroup homomorphism between inverse semigroups automatically respects the involutions.)  We also write $\IG$ for the category of \emph{inductive groupoids}, i.e.~the ordered groupoids whose object set is a semilattice (under the order inherited from the containing groupoid).  Morphisms in $\IG$ are the \emph{inductive functors}, i.e.~the ordered groupoid functors whose object maps are semilattice morphisms.  The ESN Theorem is as follows:

\begin{thm}\label{thm:ESN}
The category $\IS$ of inverse semigroups (with semigroup homomorphisms) is isomorphic to the category $\IG$ of inductive groupoids (with inductive functors).
\end{thm}

We prove this by
%\bit
%\item 
first identifying (in Section \ref{subsect:I}) the image of the category $\IS$ under the isomorphism $\bG:\RSS\to\CPG$, and then 
%\item 
showing (in Section \ref{subsect:ESN}) that this image is isomorphic to~$\IG$.
%\eit

%We will give a proof of this theorem in Section \ref{subsect:ESN}, relying on our above results on regular $*$-semigroups.  En route to this, in Section \ref{subsect:I} we identify the image of the category $\IS$ under the functor $\bG:\RSS\to\CPG$; see Proposition \ref{prop:triv}.

%show how we are somewhat-inevitably led to the ESN Theorem by considering the simplifications that arise when we specialise our general theory to inverse semigroups.

\subsection{The chained projection groupoid associated to an inverse semigroup}\label{subsect:I}

In Theorem \ref{thm:iso} we showed that the functors
\[
\bG:\RSS\to\CPG \AND \bS:\CPG\to\RSS
\]
are mutually inverse isomorphisms between
\bit
\item the category $\RSS$ of regular $*$-semigroups, with $*$-semigroup homomorphisms, and
\item the category $\CPG$ of chained projection groupoids, with chained projection functors.
\eit
It follows that for any subcategory $\CC$ of $\RSS$, the functor $\bG$ restricts to an isomorphism from $\CC$ onto its image $\bG(\CC)$ in $\CPG$:
\begin{center}
\scalebox{0.8}{
\begin{tikzpicture}
\draw [fill=blue!20,rounded corners] (0,0)--(4,0)--(4,5)--(0,5)--cycle; \node (RSS) at (2,4) {$\RSS$};
\draw [fill=red!30,rounded corners] (1,1)--(3,1)--(3,3)--(1,3)--cycle; \node (IS) at (2,2) {$\CC$};
\begin{scope}[shift={(8,0)}]
\draw [fill=blue!20,rounded corners] (0,0)--(4,0)--(4,5)--(0,5)--cycle; \node (CPG) at (2,4) {$\CPG$};
\draw [fill=red!30,rounded corners] (1,1)--(3,1)--(3,3)--(1,3)--cycle; \node (???) at (2,2) {$\bG(\CC)$};
\end{scope}
\draw[-{latex}] (RSS) to [bend left = 10] (CPG);
\draw[-{latex}] (CPG) to [bend left = 10] (RSS);
\draw[dashed,{latex}-{latex}] (IS) to (???);
\node () at (6,4.7) {$\bG$};
\node () at (6,3.3) {$\bS$};
\end{tikzpicture}
}
\end{center}
In particular, the category $\IS$ of inverse semigroups is isomorphic to its image $\bG(\IS)$.  The next result identifies this image.

\newpage

\begin{prop}\label{prop:triv}
If $(P,\G,\ve)$ is a chained projection groupoid, then the following are equivalent:
\ben
\item \label{triv1} $q\th_p=p\th_q$ for all $p,q\in P$,
\item \label{triv2} ${\leqF}={\leq}$,
\item \label{triv3} ${\F}=\De_P$,
\item \label{triv4} $\C(P)=P$ (in which case $\ve$ is the inclusion  $\ve=\io:P\hookrightarrow\G$),
\item \label{triv5} $\im(\ve)=P$,
\item \label{triv6} $(P,\G,\ve)=\bG(S)$ for some inverse semigroup $S$.
\een
\end{prop}

\pf
\firstpfitem{\ref{triv1}$\implies$\ref{triv2}}  This is clear, upon inspecting \eqref{eq:leqP} and \eqref{eq:leqF}.

\pfitem{\ref{triv2}$\implies$\ref{triv3}}  We have ${\F}={\leqF}\cap{\geqF}={\leq}\cap{\geq}=\De_P$.

\pfitem{\ref{triv3}$\implies$\ref{triv4}}  The only $P$-paths have the form $(p,p,\ldots,p)\approx(p)\equiv p$.

\pfitem{\ref{triv4}$\implies$\ref{triv5}}  This is clear.

%\pfitem{\ref{triv5}$\implies$\ref{triv1}}  Let $p,q\in P$, and put $p'=q\th_p$ and $q'=p\th_q$; we must show that $p'=q'$.  Lemma \ref{lem:p'q'} gives $p'\F q'$, and so we have the $P$-chain $[p',q']$.  Since $\ve[p',q']\in\im(\ve)=P$ by assumption, we have $\ve[p',q']=\bd(\ve[p',q'])=p'$ by Lemma \ref{lem:ve}\ref{ve2}.  We obtain $\ve[p',q']=q'$ in similar fashion, and so $p'=q'$, as required.

\pfitem{\ref{triv5}$\implies$\ref{triv6}}  By Proposition \ref{prop:GoS} we have $(P,\G,\ve)=\bG(S)$, for $S=\bS(P,\G,\ve)$ as in Definition \ref{defn:pr}.  By Proposition \ref{prop:ES}, and by assumption, we have ${E(S)\sub\im(\ve)=P=P(S)}$.  It then follows that $E(S)=P(S)$, so $S$ is inverse by Lemma \ref{lem:inv}.

%We show that $S$ is inverse by showing that $E=E(S)$ is commutative.  By Proposition \ref{prop:ES}, and by assumption, we have $E\sub\im(\ve)=P=P(S)$, and so $E=P$.  But also $E=P^2$ by Lemma~\ref{lem:PS1}\ref{PS12}, so it follows that $P=P^2$, i.e.~that $P$ is a subsemigroup of $S$.  Now let $e,f\in E=P$.  Since $e$, $f$ and $e\pr f$ are all projections, we have $e\pr f = (e\pr f)^* = f^*\pr e^* = f\pr e$, as required.

%\aftercases We now tie these in with condition \ref{triv6}.
%
%\pfitem{\ref{triv6}$\implies$\ref{triv1}}  We saw this in \eqref{eq:effe}.

\pfitem{\ref{triv6}$\implies$\ref{triv1}}  This follows from Lemma \ref{lem:inv}.
%By Definition \ref{defn:PS}, and since idempotents commute, we have
%\[
%q\th_p = pqp = pq = qpq = p\th_q \qquad\text{for $p,q\in P$.}  \qedhere
%\]
\epf

\begin{rem}
%Items \ref{triv1}--\ref{triv4} concerned $P$ (from the chained projection groupoid $(P,\G,\ve)$), while items \ref{triv4} and \ref{triv5} concerned $\ve$ (\ref{triv4} concerned both), but none concerned $\G$.  
At the outset, one might have expected the following condition to be included in the proposition:
\bit
\item $\G$ is inductive (i.e., $(P,{\leq})$ is a semilattice).
\eit
However, as in Remark \ref{rem:inv}, this condition is necessary, but not sufficient, for $S$ to be inverse; see also Examples \ref{eg:AG} and \ref{eg:Rees}.  
%it follows from Example \ref{eg:AG} (or Example~\ref{eg:Rees}) that this is not a sufficient condition for $S$ to be inverse.  Indeed, there we saw that the groupoid $\G(S)$ is inductive for any adjacency semigroup $S=A(\Ga)$ over a reflexive, symmetric digraph $\Ga$.  But one can easily check that $S=A(\Ga)$ is inverse if and only if $\Ga$ consists only of loops at each vertex.
\end{rem}

%The next result is essentially a reformulation of Proposition \ref{prop:invGve}, but it seems to be worth recording.
%
%\begin{prop}
%Let $S$ be a regular $*$-semigroup, with projection algebra $P=P(S)$.  Then $S$ is inverse if and only if ${\F}=\De_P$.
%\end{prop}
%
%\pf
%Clearly ${\F}=\De_P$ is equivalent to $\C=P$, and the result then follows from Proposition~\ref{prop:invGve} (and Theorem \ref{thm:iso}).
%\epf

%Now that we understand the image $\bG(\IS)$ of $\IS$ in $\CPG$ (cf.~Proposition~\ref{prop:triv}), we are a step closer to proving Theorem \ref{thm:ESN}.  One might now go on to show that $\bG(\IS)$ is isomorphic to~$\IG$, the category of inductive groupoids.  
%%However, there are still some subtle points remaining to be dealt with, some so subtle as to be almost invisible at this point.  
%Rather than continue on this path, however, we take a more direct/canonical route in the next section, though we will make use of Proposition~\ref{prop:triv} on a number of occasions.

\begin{defn}
We call a chained projection groupoid \emph{trivial} if it satisfies any (and hence all) of the conditions of Proposition \ref{prop:triv}.  We denote by $\TCPG$ the full subcategory of $\CPG$ consisting of all trivial chained projection groupoids (and all chained projection functors between them).
\end{defn}

\subsection{The Ehresmann--Schein--Nambooripad Theorem}\label{subsect:ESN}

Theorem \ref{thm:iso} and Proposition \ref{prop:triv} show that the categories $\IS$ (of inverse semigroups) and $\TCPG$ (of trivial chained projection groupoids) are isomorphic.  Thus, we can complete the proof of the ESN Theorem by showing that the category $\IG$ of inductive groupoids is also isomorphic to~$\TCPG$.  We do this by constructing a pair of mutually inverse functors
\[
\bI:\TCPG\to\IG \AND \bT:\IG\to\TCPG.
\]
The first will simply be a forgetful functor:
\bit
\item On objects, $\bI(P,\G,\io) = \G$.
\item On morphisms, $\bI(\phi)=\phi$.
\eit

\begin{prop}\label{prop:F1functor}
$\bI$ is a functor $\TCPG\to\IG$.
\end{prop}

\pf
Let $(P,\G,\io)$ be a trivial chained projection groupoid.  Comparing Proposition \ref{prop:triv}\ref{triv1} with Lemma \ref{lem:meet}, it is clear that $P=v\G$ is a $\wedge$-semilattice, so $\G$ is inductive.  It remains to check that any chained projection functor $\phi:(P,\G,\io)\to(P',\G',\io')$ in $\TCPG$ is an inductive functor~$\G\to\G'$.  By definition, $\phi$ is an ordered groupoid functor $\G\to\G'$, and $v\phi$ is a projection algebra morphism $P\to P'$.  Again keeping Proposition \ref{prop:triv}\ref{triv1} and Lemma \ref{lem:meet} in mind, it follows that
\[
(p\wedge q)\phi = (q\th_p)\phi = (q\phi)\th_{p\phi} = p\phi \wedge q\phi \qquad\text{for any $p,q\in P$.} \qedhere
\]
\epf

The definition of the functor $\bT:\IG\to\TCPG$ is relatively straightforward as well, although there are a few details to verify before we can give it.

\begin{lemma}\label{lem:IGPG}
If $\G$ is an inductive groupoid, with object semilattice $v\G=P$, then
\ben
\item \label{IGPG1}  $P$ is a projection algebra with respect to the maps
\[
\th_p:P\to P \GIVENBY q\th_p = p\wedge q \qquad\text{for $p,q\in P$,}
\]
\item \label{IGPG2}  $(P,\G)$ is a projection groupoid.
\een
\end{lemma}

\pf
\firstpfitem{\ref{IGPG1}}  It is a routine matter to check that axioms \ref{P1}--\ref{P5} hold.

\pfitem{\ref{IGPG2}}  It is easy to see that $(P,\G)$ is a \emph{weak} projection groupoid; indeed, this follows quickly from the construction of the projection algebra operations, and the fact that $p\leq q \iff p=p\wedge q$ in a $\wedge$-semilattice.  It remains to verify \ref{G2}, and for this we check~\ref{G2d}.  To do so, fix some $a\in\G$, and write $p=\bd(a)$ and $q=\br(a)$.  We must show that $\vt_a:p^\da\to q^\da$ is a projection algebra morphism, i.e.~that
\begin{align*}
(s\th_t)\vt_a &= (s\vt_a) \th_{t\vt_a} &&\hspace{-3cm}\text{for all $s,t\leq p$.}
\intertext{By the definition of the $\th$ maps, this amounts to showing that}
(s\wedge t)\vt_a &= s\vt_a\wedge t\vt_a &&\hspace{-3cm}\text{for all $s,t\leq p$,}
\end{align*}
i.e.~that $\vt_a$ is a semilattice morphism.  But this follows quickly from basic order-theoretic facts:
\bit
\item $p^\da$ and $q^\da$ are both $\wedge$-semilattices, as order ideals of the $\wedge$-semilattice $P$, and
\item $\vt_a:p^\da\to q^\da$ and its inverse $\vt_a^{-1}=\vt_{a^{-1}}:q^\da\to p^\da$ are both order-preserving, by Lemma \ref{lem:vtavta*}.
\eit
It follows that $\vt_a$ preserves meets.
\epf

\begin{lemma}\label{lem:IGCPG}
If $\G$ is an inductive groupoid, with object semilattice $v\G=P$, then 
\ben
\item \label{IGCPG1}  $\C(P)=P$, where the projection algebra structure on $P$ is as in Lemma \ref{lem:IGPG},
\item \label{IGCPG2}  $(P,\G,\io)$ is a (trivial) chained projection groupoid, where $\io:P\hookrightarrow\G$ is the inclusion.
%\item \label{IGCPG3}  $\bbS(\G)=\bS(P,\G,\io)$ is an inverse semigroup.
\een
\end{lemma}

\pf
\firstpfitem{\ref{IGCPG1}}  For any $p,q\in P$ we have $q\th_p = p\wedge q = q\wedge p = p\th_q$.  As in the proof of Proposition~\ref{prop:triv}, this implies ${\leqF}={\leq}$, and then $\F=\De_P$, from which the claim follows.

\pfitem{\ref{IGCPG2}}  By Lemma \ref{lem:IGPG}\ref{IGPG2}, and since $\io:P\hookrightarrow\G$ is clearly an ordered $v$-functor, it remains to check that \ref{G3} holds, and this is essentially trivial.  Indeed, consider some $b$-linked pair $(e,f)$, where $b\in\G$, and let $e_1,e_2,f_1,f_2\in P$ be as in \eqref{eq:e1e2f1f2}.  Since ${\F}=\De_P$, it follows from Lemma \ref{lem:LP}\ref{LP3} that in fact $e_1=e_2=e$ and $f_1=f_2=f$.  But then
\[
\lam(e,b,f) = \io[e,e]\circ{}_e\corest b\circ\io[f,f] = e\circ{}_e\corest b\circ f = {}_e\corest b \ANDSIM \rho(e,b,f) = {}_e\corest b.  \qedhere
\]
%
%\pfitem{\ref{IGCPG3}}  The underlying sets of $\bbS(\G)$ and $\bS(P,\G,\io)$ are both $\G$.  To check that the operations~$\star$ and~$\pr$ coincide, let $a,b\in\G$.  Write $e=\br(a)$ and $f=\bd(b)$, and also set $g=e\wedge f$.  Following Definition~\ref{defn:pr} (applied to the chained projection groupoid $(P,\G,\io)$), we have
%\[
%a\pr b = a\rest_{e'} \circ \io[e',f'] \circ {}_{f'}\corest b \WHERE e'=f\th_e \ANd f'=e\th_f.
%\]
%But $e'=f\th_e=e\wedge f=g$, and similarly $f'=g$, so in fact
%\[
%a\pr b = a\rest_{g} \circ \io[g,g] \circ {}_{g}\corest b = a\rest_{g} \circ g \circ {}_{g}\corest b = a\rest_{g}\circ {}_{g}\corest b = a\star b.
%\]
%It remains to check that the regular $*$-semigroup $S=\bS(P,\G,\io)$ is inverse, and this follows from Proposition \ref{prop:triv} and part \ref{IGCPG1} of the current lemma, as $(P,\G,\io)=\bG(S)$.
\epf

We can now define the functor $\bT:\IG\to\TCPG$:
\bit
\item On objects, $\bT(\G) = (P,\G,\io)$ is the trivial chained projection groupoid from Lemma \ref{lem:IGCPG}.
\item On morphisms, $\bT(\phi)=\phi$.
\eit

\begin{prop}\label{prop:F2functor}
$\bT$ is a functor $\IG\to\TCPG$.
\end{prop}

\pf
It remains to check that any inductive functor $\phi:\G\to\G'$ is a chained projection functor $(P,\G,\io)\to(P',\G',\io')$, where $(P,\G,\io)=\bT(\G)$ and $(P',\G',\io')=\bT(\G')$.  
Following Definition~\ref{defn:ChPG}, we need to show that:
\ben
\item \label{F1} $\phi$ is a projection groupoid morphism, and
\item \label{F2} $\phi$ respects the evaluation maps.
\een

\pfitem{\ref{F1}}  Since $\phi$ is an ordered groupoid functor by assumption, we just need to check that $v\phi$ is a projection algebra morphism ${P\to P'}$.
This follows from the fact that $v\phi$ is a semilattice morphism (as $\phi$ is inductive), and keeping the rule $q\th_p=p\wedge q$ in mind.

\pfitem{\ref{F2}}  Since $\C(P)=P$ and $\C(P')=P'$, the functor $\Phi=\C(v\phi)$ at the top of the diagram in Definition \ref{defn:ChPG} is just the object map $\Phi=v\phi$.  The diagram then becomes:
\[
\begin{tikzcd}
P \arrow{rr}{v\phi} \arrow[hook,swap]{dd}{\io} & ~ & P' \arrow[hook]{dd}{\io'} \\%
~&~&~\\
\G \arrow{rr}{\phi}& ~ & \G',
\end{tikzcd}
\]
and this obviously commutes.
\epf

Now that we have defined the two functors $\bI$ and $\bT$, it is completely routine to check that they are mutually inverse.  As noted above, this completes the proof of Theorem \ref{thm:ESN}.

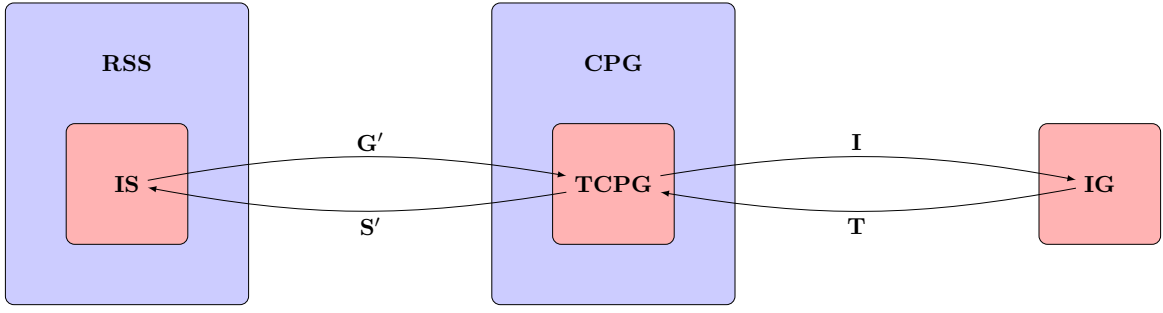
\begin{figure}[h]
\begin{center}
\scalebox{0.8}{
\begin{tikzpicture}
\draw [fill=blue!20,rounded corners] (0,0)--(4,0)--(4,5)--(0,5)--cycle; \node (RSS) at (2,4) {$\RSS$};
\draw [fill=red!30,rounded corners] (1,1)--(3,1)--(3,3)--(1,3)--cycle; \node (IS) at (2,2) {$\IS$};
\begin{scope}[shift={(8,0)}]
\draw [fill=blue!20,rounded corners] (0,0)--(4,0)--(4,5)--(0,5)--cycle; \node (CPG) at (2,4) {$\CPG$};
\draw [fill=red!30,rounded corners] (1,1)--(3,1)--(3,3)--(1,3)--cycle; \node (TCPG) at (2,2) {$\TCPG$};
\end{scope}
\begin{scope}[shift={(16,0)}]
\draw [fill=red!30,rounded corners] (1,1)--(3,1)--(3,3)--(1,3)--cycle; \node (IG) at (2,2) {$\IG$};
\end{scope}
%\draw[-{latex}] (RSS) to [bend left = 10] (CPG);
%\draw[-{latex}] (CPG) to [bend left = 10] (RSS);
%\node () at (6,4.7) {$\bG$};
%\node () at (6,3.3) {$\bS$};
\draw[-{latex}] (IS) to [bend left = 10] (TCPG);
\draw[-{latex}] (TCPG) to [bend left = 10] (IS);
\node () at (6,2.7) {$\bG'$};
\node () at (6,1.3) {$\bS'$};
\draw[-{latex}] (TCPG) to [bend left = 10] (IG);
\draw[-{latex}] (IG) to [bend left = 10] (TCPG);
\node () at (14,2.7) {$\bI$};
\node () at (14,1.3) {$\bT$};
\end{tikzpicture}
}
\caption{The categories and isomorphisms discussed in Remark \ref{rem:ESN}.}
\label{fig:isos}
\end{center}
\end{figure}

\begin{rem}\label{rem:ESN}
Restricting the isomorphisms $\bG:\RSS\to\CPG$ and $\bS:\CPG\to\RSS$ as appropriate, we obtain isomorphisms
\[
\bG':\IS\to\TCPG \AND \bS':\TCPG\to\IS,
\]
as shown in Figure \ref{fig:isos}.  Composing these with $\bI$ and $\bT$, we therefore obtain isomorphisms
\[
\bbG = \bI\circ\bG':\IS\to\IG \AND \bbS = \bS'\circ\bT:\IG\to\IS.
\]
The functors $\bbG$ and $\bbS$ both act identically on morphisms.  At the object level, we have
\[
\bbG(S)=\bI(P(S),\G(S),\ve(S))=\G(S) \qquad\text{for an inverse semigroup $S$,}
\]
so that $\bbG$ is simply the restriction to $\IS\sub\RSS$ of the functor $\G:\RSS\to\OG$ (cf.~Definition~\ref{defn:GS} and Proposition \ref{prop:calGfunctor}).  
On the other hand, $\bbS$ maps an inductive groupoid $\G$ to the inverse semigroup
\[
S = \bbS(\G) = \bS(P,\G,\io) ,
\]
where $(P,\G,\io)=\bT(\G)$ is the trivial chained projection groupoid from Lemma \ref{lem:IGCPG}\ref{IGCPG2}.  Since~$S$ is inverse, the product $\pr$ given in Definition \ref{defn:pr} takes on a simpler form here.  Indeed, let $a,b\in S(\equiv \G)$, write $p=\br(a)$ and $q=\bd(b)$, and also let $e=p\wedge q$.  Examining Proposition \ref{prop:triv}\ref{triv1} and Lemma \ref{lem:meet}, we note that $e=p'=q'$, where these projections are as in Definition \ref{defn:pr}.  It then follows that
\[
a\pr b 
%= a\restr_{p'}\circ\io[p',q']\circ{}_{q'}\corest b 
= a\rest_e\circ\io[e,e]\circ{}_e\corest b 
= a\rest_e\circ e\circ{}_e\corest b 
= a\rest_e\circ{}_e\corest b .
\]
In this way, we see that the functors $\bbG = \bI\circ\bG'$ and $\bbS = \bS'\circ\bT$ are precisely those used in the canonical proof of the ESN Theorem; see \cite[Chapter 4]{Lawson1998}.
\end{rem}

\begin{rem}
Alternatively, one could also obtain a direct (and quite transparent) proof of the ESN Theorem by specialising our proof of Theorem \ref{thm:iso} to the inverse case.  We will not give the full outline of this, but will comment briefly on the most involved step, which is to show that the product $\pr$ defined on the groupoid $\G$ is associative.  As in the previous remark, this product reduces (in the inductive/inverse case) to
\[
a\pr b = a\rest_{e} \circ {}_{e}\corest b \qquad\text{for $a,b\in\G$, where $e=\br(a)\wedge \bd(b)$.}
\]
It follows that for $a,b,c\in\G$, we have
\[
(a\pr b)\pr c = a\rest_f\circ b\rest_g\circ c\rest_h \AND a\pr(b\pr c) = a\rest_{f'}\circ b\rest_{g'}\circ c\rest_{h'}
\]
for some $f,g,h,f',g',h'\in P=v\G$.  We can then show that $(a\pr b)\pr c$ and $a\pr(b\pr c)$ are equal by showing that they have equal domains and equal ranges.  This follows from \eqref{eq:rabc} and \eqref{eq:dabc}, which were proved quite early in Section \ref{subsect:SPGve}.
\end{rem}

\begin{rem}
Our proof of the ESN Theorem involved showing that two categories of categorical structures ($\IG$ and $\TCPG$) were isomorphic.  This bears some resemblance to the approach of \cite{MV2022}, which gives a direct proof of the equivalence of Nambooripad's inductive groupoids \cite{Nambooripad1979} and cross-connections \cite{Nambooripad1994}.
\end{rem}

\section{Conclusion and discussion of future work}\label{sect:conclusion}

The main result of this paper, Theorem \ref{thm:iso}, states that the category $\RSS$ of regular $*$-semigroups (with $*$-semigroup homomorphisms) is isomorphic to the category $\CPG$ of chained projection groupoids (with chained projection functors).  At the object level, the functor $\bG:\RSS\to\CPG$ maps a regular $*$-semigroup $S$ to $\bG(S)=(P,\G,\ve)$, where
\bit
\item $P=P(S)$ is the projection algebra of $S$, a classical structure that goes back to Imaoka \cite{Imaoka1983},
\item $\G=\G(S)$ is the ordered groupoid associated to $S$, constructed in entirely the same way as inductive groupoids are from inverse semigroups \cite{Lawson1998}, and
\item $\ve=\ve(S)$ is a functor from the chain groupoid $\C(P)$ to $\G$, which provides the key link between the structures of $P$ and $\G$.
\eit
As explained in Remark \ref{rem:structure}, the existence of the inverse functor $\bS=\bG^{-1}:\CPG\to\RSS$ means that the structure of $S$ is entirely determined by that of these smaller and simpler parts, in the sense that $S=\bS(P,\G,\ve)$ can be reconstructed directly from them.
This has a number of important consequences, and leads to a number of interesting applications that will be the subject of future articles \cite{Paper2,Paper3}, and which we will briefly discuss here.

\subsection[Fundamental regular $*$-semigroups]{\boldmath Fundamental regular $*$-semigroups}\label{subsect:MP}

One of the important themes in early papers in the area was the classification of \emph{fundamental} regular $*$-semigroups, i.e.~those with no non-trivial idempotent-separating congruences.  See for example the papers of Imaoka and Yamada \cite{Yamada1981,Imaoka1980,Imaoka1983}, and the more recent work of Jones~\cite{Jones2012}.  Some of these papers drew inspiration from Munn's classical construction of fundamental \emph{inverse} semigroups~\cite{Munn1970}, which were built from order-isomorphisms between principal ideals of a semilattice.  The need for a regular $*$-analogue of a semilattice eventually led to the definition of (what we call) projection algebras in~\cite{Imaoka1983}.  

In the sequel \cite{Paper2}, we will use our groupoid approach to provide what we believe is a more transparent analysis of fundamental regular $*$-semigroups.  For one thing, we can use our~$\vt$ and~$\Th$ maps to describe the maximum idempotent-separating congruence $\mu_S$ of a regular $*$-semigroup~$S$,~\emph{viz.}:
\[
\mu_S = \set{(a,b)\in S\times S}{\vt_a=\vt_b} = \set{(a,b)\in S\times S}{\Th_a=\Th_b \text{ and } \Th_{a^*}=\Th_{b^*}}.
\]
We can also construct a maximum fundamental regular $*$-semigroup $\M_P$ with a given projection algebra $P$ by first building an appropriate chained projection groupoid $(P,\M,\ve)$, and then defining $\M_P = \bS(P,\M,\ve)$.  As with some previous studies (see especially \cite[Section~4]{Jones2012}), our construction involves projection algebra isomorphisms $p^\da\to q^\da$ between down-sets in~$P$.  The operations in $\M$ are simply function composition, inversion and restriction.  The fundamental semigroups constructed in the literature \cite{Yamada1981,Imaoka1980,Imaoka1983,Jones2012} all have somewhat-complicated products, but these are all easily understood as special instances of the $\pr$ product from Definition \ref{defn:pr}.  The morphisms $\ve[p,q]$ featuring in such $\pr$ products have the form $\ve[p,q]=\th_q|_{p^\da}$.

Given a regular $*$-semigroup $S$, with $P=P(S)$, each $\vt_a$ ($a\in S$) is a projection algebra isomorphism $\bd(a)^\da\to\br(a)^\da$ (cf.~\ref{G2d}), so we have a well-defined map $\phi:S\to\M_P:a\mt\vt_a$.  One of the main results of \cite{Paper2} is that $\phi$ is a $*$-semigroup homomorphism with $\ker(\phi)=\mu_S$, and whose image $\im(\phi)\cong S/\mu_S$ is the fundamental image of $S$.

Specialising to diagram monoids (cf.~Section~\ref{subsect:D}), we will show in \cite{Paper2} that when $P=P(\PP_n)$ is the projection algebra of a finite partition monoid $\PP_n$, the maximum fundamental regular $*$-semigroup $\M_P$ is (isomophic to) $\PP_n$ itself.  This is achieved by continuing the combinatorial analysis of the projection algebra $P(\PP_n)$ initiated in previous works such as \cite{EG2017}.

\subsection[Free (idempotent- and projection-generated) regular $*$-semigroups]{\boldmath Free (idempotent- and projection-generated) regular $*$-semigroups}\label{subsect:PGP}

In Section \ref{subsect:PA} we showed that the construction $S\mt P(S)$ of the projection algebra of a regular $*$-semigroup is the object part of a functor $P:\RSS\to\PA$.  It is natural to think of $P$ as a \emph{forgetful functor}.  Indeed, identifying $\RSS$ with $\CPG$ (via the isomorphism from Theorem~\ref{thm:iso}), we can think of $P$ as a functor $\CPG\to\PA$, under which $(P,\G,\ve)\mt P$.  
One of the main results of the sequel \cite{Paper3} is the following theorem:
\bit
\item The forgetful functor $P:\RSS\to\PA$ has a \emph{left adjoint}.
\eit
At the object level, the adjoint $\FF:\PA\to\RSS$ maps a projection algebra $P$ to what we will call the \emph{chain semigroup} $\C_P$.  This semigroup is constructed as $\C_P=\bS(P,\ol\C,\ve)$, where $\ol\C$ is a certain quotient of the chain groupoid $\C=\C(P)$, and $\ve$ is the canonical quotient map.  (The definition of $\ol\C$ is very natural, but too involved to give here.  Suffice it to say that it involves the concept of \emph{linked pairs} in $P$, which are regular $*$-analogues of Nambooripad's \emph{singular squares} in biordered sets \cite{Nambooripad1979}.)  There are two main consequences of the above theorem.

First, since the image of $\FF$ turns out to be a full subcategory of $\RSS$, we can identify $\PA$ with a \emph{coreflective} subcategory of $\RSS$.  This has an analogue in inverse semigroups, where it can be seen that the category of semilattices is coreflective in $\IS$.  We are unaware of a statement of this fact in the literature (but we note that \cite[Theorem 2.4.2]{Lawson1998} says that the category of groups is \emph{reflective} in $\IS$).  This lends extra support to the claim that projection algebras are the appropriate regular $*$-analogue of semilattices.

The second main consequence of the above theorem is that the semigroups $\C_P$ are therefore the \emph{$P$-free objects} in the category $\RSS$.%
\footnote{This notion of free-ness is different from the \emph{varietal} sense; free regular $*$-semigroups with respect to the forgetful functor $\RSS\to\Set$ are studied in \cite{Polak2001}.  We note, however, that the left adjoint $\Set\to\RSS$ to \emph{this} forgetful functor does not map onto a full subcategory of $\RSS$, so that $\Set$ is not (isomorphic to) a coreflective subcategory of $\RSS$.  (This is analogous to other varietal free constructions.  For example, there exist morphisms $X^+\to Y^+$ between free semigroups that do not arise from set maps $X\to Y$.)}
%
%:
%\ben
%\item \label{31} The semigroups $\C_P$ are therefore the \emph{$P$-free objects} in the category $\RSS$.
%\item \label{32} Since the image of $\FF$ is a full subcategory of $\RSS$, we can identify $\PA$ with a \emph{coreflective} subcategory of $\RSS$.
%\een
%Statement \ref{32} has an analogue in inverse semigroups, where it can be seen that the category of semilattices is coreflective in $\IS$.  We are unaware of a statement of this fact in the literature (but we note that \cite[Theorem 2.4.2]{Lawson1998} states that the category of groups is reflective in $\IS$).  This lends extra support to the claim that projection algebras are the appropriate regular $*$-analogue of semilattices.
%
%Statement \ref{31} 
This conveys a lot of deep information, which we will not fully unpack here (see \cite{Paper3} for details), but one of the main points is the following theorem:
\bit
\item If $P$ is a projection algebra and $S$ a regular $*$-semigroup, then for any projection algebra morphism $\phi:P\to P(S)$ there is a unique $*$-semigroup homomorphism $\Phi:\C_P\to S$ such that the following diagram commutes (where both `$\io$'s denote inclusion maps):
\[
\begin{tikzcd}
P \arrow{rr}{\phi} \arrow[swap,hookrightarrow]{dd}{\io} & ~ & P(S) \arrow[hookrightarrow]{dd}{\io} \\%
~&~&~\\
\C_P \arrow{rr}{\Phi}& ~ & S.
\end{tikzcd}
\]
\eit
As we will explain in \cite{Paper3}, the semigroups $\C_P$ are the regular $*$-analogues of the free (idempotent-generated) semigroups $\FIG(E)$ associated to a biordered set $E$, studied by Nambooripad, Easdown and others \cite{Nambooripad1979,NP1980,Easdown1985,GR2012,GR2012b,DGR2017,DG2014}.  The semigroup $\FIG(E)$ can be defined in terms of a presentation by generators and relations:
\[
\FIG(E) = \pres{X_E}{x_ex_f=x_{ef} \text{ for } (e,f)\in \bbB},
\]
%\bit
%\item $\FIG(E) = \pres{X_E}{x_ex_f=x_{ef} \text{ for } (e,f)\in \bbB}$,
%\eit
where here
%\bit
%\item 
$X_E=\set{x_e}{e\in E}$ is an alphabet in one-one correspondence with $E$, and
%\item 
$\bbB$ is the set of all \emph{basic pairs} in $E$, i.e.~all pairs $(e,f)\in E\times E$ for which $\{ef,fe\}\cap\{e,f\}\not=\es$.
%\eit
One of the main results in \cite{Paper3} is a presentation for the chain semigroup $\C_P$:
\[
\C_P\cong\FPG(P) = \pres{X_P}{x_p^2=x_p,\ (x_px_q)^2=x_px_q,\ x_px_qx_p=x_{q\th_p} \text{ for } p,q\in P},
\]
%\bit
%\item $\C_P\cong\FPG(P) = \pres{X_P}{x_p^2=x_p,\ (x_px_q)^2=x_px_q,\ x_px_qx_p=x_{q\th_p} \text{ for } p,q\in P}$,
%\eit
where again
%\bit
%\item 
$X_P=\set{x_p}{p\in P}$ is an alphabet in one-one correspondence with $P$.
%\eit

Another result of \cite{Paper3} shows that a projection algebra $P$ uniquely determines a biordered set $E=E(P)$, in the sense that any regular $*$-semigroup $S$ with projection algebra $P(S)=P$ also has biordered set $E(S)\cong E$.  The chain semigroup $\C_P\cong\FPG(P)$ is then exhibited
%(up to isomorphism) all regular $*$-semigroup with projection algebra $P$ have the same biordered set $E=E(P)$, and exhibits $\C_P$ 
as a natural quotient of $\FIG(E)$, \emph{viz.}:
\[
\C_P \cong \FIG^*(E) = \pres{X_E}{x_ex_f=x_{ef} \text{ for } (e,f)\in \bbB\cup(P\times P)}.
\]
%\bit
%\item $\C_P \cong \FIG^*(E) = \pres{X_E}{x_ex_f=x_{ef} \text{ for } (e,f)\in \bbB\cup(P\times P)}$.
%\eit
A natural family of questions therefore arises:
\bit
\item Given an important regular $*$-semigroup $S$, with projection algebra $P=P(S)$ and biordered set $E=E(S)$, how are the structures of $\C_P\cong\FPG(P)\cong\FIG^*(E)$ and $\FIG(E)$ related, and how are these related to $S$ itself?
\eit
We will show in \cite{Paper3} that when $P=P(\TL_n)$ is the projection algebra of a (finite) \emph{Temperley-Lieb monoid}, $\C_P$ is (isomorphic to) $\TL_n$ itself.

The situation for partition monoids $\PP_n$ is rather more complicated, and will be treated in a separate paper \cite{Paper4}, where we will show that when $P=P(\PP_n)$ and $E=E(\PP_n)$:
\bit
\item the maximal subgroup of $\FPG(P)\cong\FIG^*(E)$ corresponding to an idempotent of rank $r\leq n-2$ is isomorphic to $\S_r$, the symmetric group of degree $r$, and
\item the maximal subgroup of $\FIG(E)$ corresponding to an idempotent of rank $r\leq n-2$ is isomorphic to $\Z\times\S_r$, where $\Z$ is the infinite cyclic group.
\eit
Higher rank subgroups are free (of calculable degrees).  These results naturally complement existing studies of maximal subgroups in other families of monoids, such as transformation and linear monoids \cite{BMM2009,GR2012b,DG2014,YDG2015,Dolinka2013}.  As we explain in \cite{Paper4}, the additional free factor in the subgroups $\Z\times\S_r$ of~$\FIG(E)$ comes from the fact that the biordered set $E=E(\PP_n)$ is isomorphic to that of a certain `twisted' monoid mapping onto $\PP_n$.
The above results are obtained by utilising a general (group) presentation for a maximal subgroup of an arbitrary free projection-generated regular $*$-semigroup, obtained in \cite{Paper3} via an application of the Reidemeister--Schreier rewriting machinery from \cite{Ruskuc1999}.

%Inspired by the main result of \cite{GR2012b}, we also believe the following question is very interesting:
%\bit
%\item Is every group (isomorphic to) a maximal subgroup of some free projection-generated regular $*$-semigroup?
%\eit

\subsection{Open problems}\label{subsect:OP}

Throughout the text we have mentioned a number of open problems.  For convenience, we list these here, as well as some others:
\bit
\item To what extent is a projection algebra determined by its associated relations ($\leq$, $\leqF$ and~$\F$)?
\item Is every group (isomorphic to) a maximal subgroup of some free projection-generated regular $*$-semigroup (cf.~\cite{DR2013,GR2012})?
\item Does the groupoid-approach developed here have consequences for the representation theory of regular $*$-semigroups (cf.~\cite{Steinberg2008,Steinberg2010,Steinberg2006, Stein2016,Stein2017,Stein2019,Stein2020,MS2021,Stein2022})?
\item Can the theory developed above be applied to other categories of `projection-based' semigroups (cf.~\cite{Jones2012,Wang2022})?
\item Is there a parallel theory for (non-projection-based) generalisations of regular $*$-semigroups, e.g.~for involution semigroups satisfying $aa^*aa^*=aa^*$ in place of $aa^*a=a$?
\item Does the isomorphism $\RSS\cong\CPG$ lead to an efficient algorithm for enumerating finite regular $*$-semigroups (cf.~\cite{Malandro2019})?
\item Can one develop a symmetrical theory of cross-connections for regular *-semigroups (cf.~\cite{MV2022,MV2019,Nambooripad1994,Grillet1974a,Grillet1974b,Grillet1974c,Grillet1974d})?
\eit

\footnotesize
\def\bibspacing{-1.1pt}
\bibliography{biblio}
\bibliographystyle{abbrv}

\end{document}